\documentclass[a4paper,reqno]{amsart}
\usepackage[english]{babel}

\usepackage{amsmath,amstext,amssymb,amsthm}
\usepackage{mathrsfs}
\usepackage{bbm}

\usepackage{enumitem}
\setlist[enumerate]{leftmargin=*}
\setlist[itemize]{leftmargin=*}

\def\Xint#1{\mathchoice
{\XXint\displaystyle\textstyle{#1}}%
{\XXint\textstyle\scriptstyle{#1}}%
{\XXint\scriptstyle\scriptscriptstyle{#1}}%
{\XXint\scriptscriptstyle\scriptscriptstyle{#1}}%
\!\int}
\def\XXint#1#2#3{{\setbox0=\hbox{$#1{#2#3}{\int}$ }
\vcenter{\hbox{$#2#3$ }}\kern-.6\wd0}}

\def\dashint{\Xint-}

\usepackage{hyperref}

\newcommand{\arXivlink}[1]{\href{https://arxiv.org/abs/#1}{\texttt{arXiv:#1}}}

\newcommand*\RR{\mathbb{R}}
\newcommand*\CC{\mathbb{C}}
\newcommand*\NN{\mathbb{N}}
\newcommand*\ZZ{\mathbb{Z}}

\newcommand*{\admR}{\mathscr{R}}
\newcommand*{\cP}{\mathscr{P}}
\newcommand*{\cS}{\mathscr{S}}
\newcommand*\adm{\mathscr{R}_{\mathrm{ad}}}
\newcommand*\chil{\mathscr{C}}
\newcommand*\BMO{\mathrm{BMO}}

\newcommand*\ext{\mathcal{X}}

\newcommand{\dnsA}{\mathcal{A}}

\newcommand*{\Rnon}{\RR_+}
\newcommand*{\Rpos}{\mathring{\RR}_+}
\newcommand*{\Rnoz}{\RR^\ast}

\newcommand*{\Npos}{\NN_+}

\newcommand*{\vfX}{Y}

\newcommand*{\even}{\mathrm{e}}
\newcommand*{\odd}{\mathrm{o}}

\newcommand*{\loc}{\mathrm{loc}}
\newcommand*{\sloc}{\mathrm{sloc}}

\newcommand*{\ind}{\mathbbm{1}}

\newcommand*{\modu}{\mathrm{m}}
\newcommand*{\Hank}{\mathrm{H}}

\newcommand{\sobolev}[2]{L^{#2}_{#1}}
\newcommand{\hormander}[1]{\mathcal{H}_2^{#1}}

\DeclareMathOperator{\sign}{sign}
\DeclareMathOperator*{\esssup}{ess\,sup}

\newcommand*\supp{\mathop{\mathrm{supp}}}
\newcommand*\arccosh{\mathop{\mathrm{arccosh}}}

\DeclareMathOperator*{\absconv}{abs\,conv}

\newcommand*\id{\mathrm{id}}

\newcommand*{\tc}{\mathrel{\,:\,}}
\newcommand*{\defeq}{\mathrel{:=}}

\newcommand{\Beta}{\mathrm{B}}

\DeclareSymbolFont{t1letters}{T1}{cmr}{m}{it}
\DeclareMathSymbol{\DD}{0}{t1letters}{208}

\newcommand*{\dist}{\varrho}

\newcommand*\Sz{\mathcal{S}}

\newcommand*\LinBnd{\mathcal{L}}

\newcommand*{\RBd}{\mathfrak{R}}

\newcommand*\fctE{\mathcal{E}}

\newcommand*\cT{\mathcal{T}}

\newcommand*\dd{\mathrm{d}}
\newcommand*\dmu{\dd\mu_\nu}
\newcommand*\dmub{\dd\bar{\mu}_\nu}
\newcommand*\dpr{\dd\varpi_\nu}
\newcommand*\dv{\dd v}
\newcommand*\du{\dd u}

\newcommand*\ddlam{\frac{\dd\lambda}{\lambda}}

\newcommand*\bdx{\mathbf{x}}
\newcommand*\bdy{\mathbf{y}}
\newcommand*\bdz{\mathbf{z}}
\newcommand*\bdzero{\mathbf{0}}

\newcommand*\Dom{\mathfrak{D}}

\newcommand*\Riesz{\mathcal{R}}

\newcommand*\iker{\mathcal{K}}
\newcommand*\iheat{\mathcal{H}}

\newcommand*\aux{\mathfrak{r}}
\newcommand*\HaId{\widetilde{\Hank}_\nu}
\newcommand*\BigO{O}
\newcommand*\rade{\varepsilon}
\newcommand*\rD{\mathrm{D}}

\newcommand{\IT}{\mathrm{IT}}

\theoremstyle{plain}
\newtheorem{thm}{Theorem}[section]
\newtheorem{lm}[thm]{Lemma}
\newtheorem{prop}[thm]{Proposition}
\newtheorem{cor}[thm]{Corollary}
\theoremstyle{definition}
\newtheorem{defin}[thm]{Definition}
\newtheorem{rem}[thm]{Remark}

\numberwithin{equation}{section}

\begin{document}
\title[Singular integrals on $ax+b$ hypergroups]{Singular integrals on $ax+b$ hypergroups and an operator-valued spectral multiplier theorem}

\author[A. Martini]{Alessio Martini}
\address[A. Martini]{Dipartimento di Scienze Matematiche ``G.L. Lagrange'' \\ Politecnico di Torino \\ Corso Duca degli Abruzzi 24 \\ 10129 Torino \\ Italy}
\email{alessio.martini@polito.it}

\author[P. Plewa]{Pawe\l{} Plewa}
\address[P. Plewa]{Dipartimento di Scienze Matematiche ``G.L. Lagrange'',
	Politecnico di Torino\\ Corso Duca degli Abruzzi 24 \\ 10129 Torino \\ Italy
	and Department of Pure and Applied Mathematics, 
	Wroc{\l}aw University of Science and Technology\\ wyb. Wys\-pia{\'n}\-skie\-go  27\\ 50–-370 Wroc{\l}aw\\ Poland }
\email{pawel.plewa@pwr.edu.pl}

\begin{abstract}
Let $L_\nu = -\partial_x^2-(\nu-1)x^{-1} \partial_x$ be the Bessel operator on the half-line $X_\nu = [0,\infty)$ with measure $x^{\nu-1} \,\dd x$.
In this work we study singular integral operators associated with the Laplacian $\Delta_\nu = -\partial_u^2 + e^{2u} L_\nu$ on the product $G_\nu$ of $X_\nu$ and the real line with measure $\dd u$. For any $\nu \geq 1$, the Laplacian $\Delta_\nu$ is left-invariant with respect to a noncommutative hypergroup structure on $G_\nu$, which can be thought of as a fractional-dimension counterpart to $ax+b$ groups. In particular, equipped with the Riemannian distance associated with $\Delta_\nu$, the metric-measure space $G_\nu$ has exponential volume growth.
We prove a sharp $L^p$ spectral multiplier theorem of Mihlin--H\"ormander type for $\Delta_\nu$, as well as the $L^p$-boundedness for $p \in (1,\infty)$ of the associated first-order Riesz transforms. To this purpose, we develop a Calder\'on--Zygmund theory \`a la Hebisch--Steger adapted to the nondoubling structure of $G_\nu$, and establish large-time gradient heat kernel estimates for $\Delta_\nu$. In addition, the Riesz transform bounds for $p > 2$ hinge on an operator-valued spectral multiplier theorem, which we prove in greater generality and may be of independent interest.
\end{abstract}



\thanks{The authors gratefully acknowledge the financial support of Compagnia di San Paolo. The first-named author is a member of the Gruppo Nazionale per l'Analisi Matematica, la Probabilit\`a e le loro Applicazioni (GNAMPA) of the Istituto Nazionale di Alta Matematica (INdAM). The second-named author is also supported by the Foundation for Polish Science (START 057.2023).
}

\maketitle
	
\section{Introduction}

\subsection{Statement of the results}

For any real parameter $\nu \geq 1$,
let $X_\nu$ be the measure space $(\Rnon,\mu_\nu)$, where $\Rnon = [0,\infty)$ and $\dmu(x) = x^{\nu-1} \,\dd x$.
The Bessel operator $L_\nu$ is the singular second-order differential operator on $X_\nu$ given by
\begin{equation}\label{eq:bessel_op}
L_\nu = -\partial_x^2 - \frac{\nu-1}{x} \partial_x.
\end{equation}
When $\nu$ is an integer, $L_\nu$ is the radial part of the Laplace operator on $\RR^\nu$, and of course the latter operator is translation-invariant. An analogous invariance property can be stated for the Bessel operator $L_\nu$ for any $\nu \geq 1$, by introducing an appropriate notion of generalised translations on $X_\nu$, known as Hankel translations. The Hankel translations and the corresponding Hankel convolution determine a commutative hypergroup structure on $X_\nu$, with which $X_\nu$ is also known as a Bessel--Kingman hypergroup. Beside being translation-invariant on $X_\nu$, the Bessel operator $L_\nu$ is also homogeneous with respect to the natural dilations on $\Rnon$, which are automorphisms of the hypergroup $X_\nu$.

We now consider the semidirect product hypergroup $G_\nu = X_\nu\rtimes \RR$, where $\RR$ acts on $X_\nu$ via dilations. We equip $G_\nu$ with the right Haar measure $\dmub(x,u) = \dmu(x) \,\dd u$ and define the natural left-invariant Laplacian $\Delta_\nu$ on $G_\nu$ by
\begin{equation}\label{eq:Deltanu}
\Delta_\nu = -\partial_u^2+e^{2u}L_\nu.
\end{equation}
We shall think of $\Delta_\nu$ as a self-adjoint operator on $L^2(G_\nu)$ with Neumann boundary conditions. We refer the reader to Section \ref{s:prelim} for additional details on the hypergroup structures of $G_\nu$ and $X_\nu$. Here we just point out that, as a measure space, $G_\nu$ is simply the product of $X_\nu$ and $\RR$, and the hypergroup translations on $G_\nu$ play a role in the analysis of the Laplacian $\Delta_\nu$ which is analogous to that of the Hankel translations for the Bessel operator $L_\nu$. Crucially, due to the presence of the scaling factor $e^{2u}$, the two summands in \eqref{eq:Deltanu} do not commute; correspondingly, $G_\nu$ is equipped with a noncommutative hypergroup structure, which is not just the direct product of those of the factors.

The operator $\Delta_\nu$ is a fractional-dimension analogue of a distinguished group-invariant Laplacian on $ax+b$ groups, which has been studied in multiple works in the literature (see, e.g., \cite{GaSj,HeSt,Ma23,SjVa}). Indeed, when $\nu$ is an integer, we can consider the semidirect product group $\RR^\nu \rtimes \RR$ and the left-invariant Laplacian $\Delta_{\RR^\nu}$ thereon given by
\[
\Delta_{\RR^\nu} = -\partial_u^2 + e^{2u} L_{\RR^\nu},
\]
where $L_{\RR^\nu} = -\sum_{j=1}^\nu \partial_{x_j}^2$ is the usual (positive-definite) Laplace operator on $\RR^\nu$; then $\Delta_\nu$ coincides with the restriction of $\Delta_{\RR^\nu}$ to $\RR^\nu$-radial functions.

We remark that the operator $\Delta_\nu$ on $G_\nu$ can be written in divergence form as 
\begin{equation}\label{eq:divform}
\Delta_\nu = \vfX_0^+ \vfX_0 + \vfX_1^+ \vfX_1, 
\end{equation}
where $\vfX_0$ and $\vfX_1$ are the vector fields on $G_\nu$ given by
\begin{equation}\label{eq:vfs}
\vfX_0=\partial_u, \qquad \vfX_1=e^u \partial_x
\end{equation}
and $\vfX_0^+$ and $\vfX_1^+$ are their formal adjoints with respect to $\dmub$. Thus, if we equip $G_\nu$ with the Riemannian structure that makes the fields $\vfX_0,\vfX_1$ an orthonormal frame, then we can think of $\Delta_\nu$ of a weighted Laplacian on this Riemannian manifold.

Notice that the vector fields \eqref{eq:vfs} are independent of $\nu$, so the same is true for the Riemannian structure, which coincides, in suitable coordinates, with that of a half-plane in the real hyperbolic plane. What depends on $\nu$ is the measure $\dmub$, which determines the adjoints $\vfX_0^+,\vfX_1^+$ and thus the Laplacian $\Delta_\nu$ as in \eqref{eq:divform}. We remark that, equipped with the Riemannian metric and the measure $\dmub$, the space $G_\nu$ has exponential volume growth for any $\nu \geq 1$.

In this work, we study $L^p$-boundedness properties of singular integral operators on $G_\nu$ naturally related with the Laplacian $\Delta_\nu$ and the underlying geometry. Due to the exponential volume growth of $G_\nu$, the standard singular integral theory for doubling metric measure spaces does not apply here; nevertheless, as we shall see, the Calder\'on--Zygmund theory developed in \cite{HeSt} and the corresponding Hardy space theory of \cite{Va-PhD} can be fruitfully applied to this setting. The present work thus confirms the broad applicability of the singular integral theory of \cite{HeSt,Va-PhD}, even beyond the settings of Lie groups and flow trees to which it has been applied so far (see, e.g., \cite{DLMV,He18,LMSTV,MOV,MSTV} and references therein).

Our first main result concerns the functional calculus for $\Delta_\nu$. As the Laplacian $\Delta_\nu$ is self-adjoint on $L^2(G_\nu)$, by the spectral theorem we have a bounded Borel functional calculus for $\Delta_\nu$, whereby to any bounded Borel function $F : \Rnon \to \CC$ there corresponds an $L^2(G_\nu)$-bounded operator $F(\Delta_\nu)$. What we prove here is an $L^p$ spectral multiplier theorem of Mihlin--H\"ormander type for $\Delta_\nu$, giving sufficient conditions, in terms of size and smoothness of the function $F$, so that the operator $F(\Delta_\nu)$ extends to an $L^p$-bounded operator for other values of $p \neq 2$.
Interestingly enough, despite the exponential volume growth of $G_\nu$, our result shows that the operator $\Delta_\nu$ has a differentiable $L^p$-functional calculus, i.e., it suffices to control a finite number of derivatives of $F$ on $\Rnon$ to ensure the $L^p$-boundedness of $F(\Delta_\nu)$. 

Let us fix a cutoff $\chi \in C_c^\infty(0,\infty)$ such that $\ind_{[1/2,2]} \leq \chi \leq \ind_{[1/4,4]}$.
Let $\sobolev{s}{2}(\RR)$ denote the $L^2$ Sobolev space of order $s$ on $\RR$. Moreover, let us write $A_+$ for the positive part $\max\{A,0\}$ of a real number $A$.

\begin{thm}\label{thm:main_mult}
Suppose that $s_0>3/2$ and $s_\infty>(\nu+1)/2$. If a bounded Borel function $F : [0,\infty)\to\CC$ satisfies
\begin{equation}\label{eq:ass_mult}
\sup_{t > 0} {(1+t)^{-(2-\nu)_+/2}} \|F(t\cdot)\chi\|_{\sobolev{s_0}{2}(\RR)} < \infty, \qquad  \sup_{t \geq 1} \|F(t\cdot)\chi\|_{\sobolev{s_\infty}{2}(\RR)} < \infty,
\end{equation}
then $F(\sqrt{\Delta_\nu})$ extends to an operator of weak type $(1,1)$ and bounded on $L^p(G_\nu)$ for $p \in (1,\infty)$, bounded from $H^1(G_\nu)$ to $L^1(G_\nu)$ and from $L^\infty(G_\nu)$ to $\BMO(G_\nu)$.
\end{thm}

We refer to Section \ref{s:CZ} below for further details on the definition of the Hardy space $H^1(G_\nu)$ and the bounded mean oscillation space $\BMO(G_\nu)$.  We point out that the regularity thresholds $3/2$ and $(\nu+1)/2$ in Theorem \ref{thm:main_mult} are sharp, i.e., they cannot be replaced by smaller quantities; we discuss this in Section \ref{ss:sharpness} below. Interestingly enough, when $\nu>1$, both thresholds are strictly larger than $1$, i.e., half the topological dimension of the manifold $G_\nu$.

The smoothness condition \eqref{eq:ass_mult} on the spectral multiplier $F$ can be thought of as a refinement of the usual scale-invariant smoothness condition of Mihlin--H\"ormander type \cite{Ho,Mi}:
\begin{equation}\label{eq:ass_mult_MH}
\|F\|_{L^2_{s,\sloc}} \defeq \sup_{t > 0} \|F(t\cdot) \chi\|_{\sobolev{2}{s}} < \infty
\end{equation}
for some $s > \max\{3,\nu+1\}/2$. The form of \eqref{eq:ass_mult} is analogous to that used in other contexts where there is a mismatch between the ``local dimension'' and the ``dimension at infinity'' of the environment space (see, e.g., \cite{A}, \cite[Theorem 2]{Si2} and \cite[Theorem 8.2]{DOSi}). In the case $\nu \geq 2$, one may assume $s_\infty \geq s_0$ and $(2-\nu)_+=0$, so \eqref{eq:ass_mult} is effectively equivalent to
\[
\sup_{0 < t \leq 1} \|F(t\cdot)\chi\|_{\sobolev{s_0}{2}(\RR)} < \infty, \qquad  \sup_{t \geq 1} \|F(t\cdot)\chi\|_{\sobolev{s_\infty}{2}(\RR)} < \infty;
\]
the latter can be thought of as a variant of \eqref{eq:ass_mult_MH}, where different smoothness requirements are imposed on the local part and the part at infinity of the spectral multiplier $F$. In the case $\nu < 2$, instead, one may assume $s_0 \geq s_\infty$, and \eqref{eq:ass_mult} can be thought of as a refinement of the condition \eqref{eq:ass_mult_MH} with $s=s_0$, due to the presence of the decaying factor $(1+t)^{-(2-\nu)/2}$ for large $t$.

Our second main result concerns the first-order Riesz transforms
\begin{equation}\label{eq:def_Riesz}
\Riesz_j = \vfX_j \Delta_\nu^{-1/2},\qquad j=0,1,
\end{equation}
associated with the Laplacian $\Delta_\nu$. From the expression \eqref{eq:divform} it is immediate that both Riesz transforms are bounded on $L^2(G_\nu)$; what we prove here is the $L^p$-boundedness for $p \in (1,\infty)$ of these singular integral operators.

\begin{thm}\label{thm:main_riesz}
The Riesz transforms $\Riesz_j$, $j=0,1$, are bounded on $L^p(G_\nu)$, $p\in(1,\infty)$, from $H^1(G_\nu)$ to $L^1(G_\nu)$ and are of weak type $(1,1)$. 
\end{thm}

While we state Theorem \ref{thm:main_riesz} as a single result, the bounds for $p > 2$ require substantially different proofs from those for $p<2$. Notice that the Riesz transforms in \eqref{eq:def_Riesz} are not skew-adjoint, so the case $p>2$ cannot be simply reduced to the case $p<2$ by duality.

When the parameter $\nu$ is an integer, both Theorems \ref{thm:main_mult} and \ref{thm:main_riesz} above correspond to analogous results for the Laplacian $\Delta_{\RR^\nu}$ on the $ax+b$ group $\RR^\nu \rtimes \RR$ proved elsewhere. Specifically, the $L^p$ spectral multiplier bounds in Theorem \ref{thm:main_mult} and the Riesz transform bounds for $p \leq 2$ in Theorem \ref{thm:main_riesz} correspond to \cite[Theorem 2.4]{HeSt}, while the Riesz transform bounds for $p> 2$ correspond to \cite[Theorem 1.1]{Ma23}. As a matter of fact, when $\nu$ is an integer, the results above could be deduced from the known results for $\Delta_{\RR^\nu}$, by restricting to appropriate classes of $\RR^\nu$-radial functions. The novelty in Theorems \ref{thm:main_mult} and \ref{thm:main_riesz} is the fact that $\nu$ is allowed to be fractional, which creates a number of additional difficulties in their proofs.

The results of \cite{HeSt} for the elliptic Laplacian $\Delta_{\RR^\nu}$ on the group $\RR^\nu \rtimes \RR$ were generalised in \cite{He99,MOV,MaVa} to the case of the semidirect extension $N \rtimes \RR$ of a noncommutative Carnot group $N$, equipped with the left-invariant sub-Laplacian $\Delta_N = -\partial_u^2 + e^{2u} L_N$, where $L_N$ is the homogeneous sub-Laplacian on $N$. In the recent work \cite{MaPl}, we proved a sharpening of the multiplier theorem of \cite{MOV} for certain classes of Carnot groups $N$. The analysis of the Laplacian $\Delta_\nu$ on the hypergroup $G_\nu$ for fractional $\nu$ turned out to be a crucial ingredient in the proofs in \cite{MaPl}. This is a further reason of interest for the operator $\Delta_\nu$, which is studied in its own right in the present paper.

The analysis developed here has actually a close connection to that of the sub-Laplacian $\Delta_N$ on the Carnot group extension $N \rtimes \RR$. Indeed, notice that the $L^p$-boundedness of the Riesz transforms associated to $\Delta_N$ was proved in \cite{MaVa} for $p \leq 2$, but the case $p > 2$ has so far remained open for nonabelian $N$. However, the approach that we develop here for $\Delta_\nu$, based on the operator-valued spectral multiplier theorem proved in Section \ref{s:opval} below, turns out to be fruitfully applicable to the case of $\Delta_N$ as well. We shall substantiate this in an upcoming paper \cite{MaPl2}, which constitutes an additional motivation for the present work.

We also point out that discrete analogues of the above results were proved in \cite{HeSt,LMSTV,MSTV} in the setting of flow trees. The latter setting can be thought of as capturing the ``coarse geometry'' of the semidirect products $\RR^\nu \rtimes \RR$ and $G_\nu$ equipped with the Laplacians $\Delta_{\RR^\nu}$ and $\Delta_\nu$. Indeed, the analogues of Theorem \ref{thm:main_mult} for flow trees only feature the threshold $3/2$ (corresponding to the ``dimension at infinity'') and there is no counterpart to the threshold $(\nu+1)/2$ there.

\subsection{Proof strategy}

The proofs of Theorem \ref{thm:main_mult} and the $p < 2$ part of Theorem \ref{thm:main_riesz} are based on two fundamental ingredients.

First of all, in Section \ref{s:CZ} we show that the Calder\'on--Zygmund theory of \cite{HeSt,Va-PhD} applies to the nondoubling space $G_\nu$.
To this purpose, we first work in the general framework of a metric measure spaces equipped with a family of ``admissible sets'' satisfying a number of geometric axioms (see Definition \ref{def:1'} below), and prove that under these assumptions any integrable function admits a Calder\'on--Zygmund decomposition in the sense of \cite{HeSt}.
We then show that the hypergroup $G_\nu$ with the measure $\dmub$ and the Riemannian distance falls under this framework, by constructing a suitable family of admissible sets in $G_\nu$. 

Secondly, we show in Section \ref{s:heatestimates} that the heat kernel associated with $\Delta_\nu$ and its gradient satisfy suitable weighted $L^1$-estimates, which are valid both for small and large times.
To this purpose, much as in \cite{MOV,MaVa}, we exploit the relation between the heat semigroups on $G_\nu$ and $X_\nu$ given by a general result of \cite{G}. Compared with the case of Lie groups studied in \cite{MOV,MaVa}, an additional technical complication here comes from the fact that hypergroup translations are particular ``averaging operators'' and are not just given by a change of variables. Correspondingly, the relation between the integral kernel and the convolution kernel of a left-invariant operator is more complicated, and some care is needed in order to deduce gradient estimates for the integral kernel from analogous estimates for the convolution kernel (see Lemma \ref{lm:l1est_int_conv} below).

Once the Calder\'on--Zygmund theory and the heat kernel estimates on $G_\nu$ are established, an adaptation of the arguments from \cite{HeSt,MOV,MaVa}, also exploiting the finite propagation speed property for the wave equation associated with $\Delta_\nu$, yields the multiplier theorem and the Riesz transform bound for $p <2$. The corresponding proofs are discussed in Section \ref{s:multriesz}, where we also show that the smoothness thresholds $3/2$ and $(\nu+1)/2$ in Theorem \ref{thm:main_mult} are sharp and cannot be lowered.

The proof of the $p>2$ part of Theorem \ref{thm:main_riesz} requires different arguments. Much as in \cite{Ma23}, our analysis relies, on the one hand, on precise asymptotics for the convolution kernels of the Riesz transforms and their adjoints at zero and at infinity, which we discuss in Section \ref{s:convRiesz}; on the other hand, in order to deal with the parts at infinity of the Riesz transforms, we shall make use of an operator-valued multiplier theorem. Again, the fact that we are working on the hypergroup $G_\nu$, instead of the Lie group $\RR^\nu \rtimes \RR$ considered in \cite{Ma23}, brings additional challenges.

Indeed, the approach in \cite{Ma23} for the study of the parts at infinity of the Riesz transform was based on the observation that $L^p(\RR^\nu \rtimes \RR) \cong L^p(\RR^\nu;L^p(\RR))$ and that, moreover, left-invariant operators on $\RR^\nu \rtimes \RR$ are also translation-invariant operators with respect to Euclidean translations in the $\RR^\nu$-variables. As a consequence, those operators can be considered as operator-valued Fourier multipliers on $\RR^\nu$, whose operator-valued symbols act on functions on $\RR$. An operator-valued multiplier theorem for the Euclidean Fourier transform already available in the literature \cite{StWe,We} could be then applied to prove the $L^p$-boundedness of the parts at infinity of the Riesz transforms on $\RR^\nu \rtimes \RR$.

In the case of the hypergroup $G_\nu = X_\nu \rtimes \RR$, it is still true that $L^p(G_\nu) \cong L^p(X_\nu;L^p(\RR))$ and that convolution operators on $G_\nu$ are translation-invariant in the $X_\nu$-variable, so they can be considered as operator-valued multipliers for the Hankel transform on $X_\nu$ (i.e., the ``radial part'' of the Fourier transform, which makes sense for fractional $\nu$ as well). However, an operator-valued multiplier theorem for the Hankel transform, analogous to that of \cite{StWe,We} for the Fourier transform, appears not to be available in the literature.

As a matter of fact, the proof in \cite{StWe} of the operator-valued Fourier multiplier theorem for integer $\nu > 1$ exploits in an essential way the multivariate structure of $\RR^\nu$, since it is based on Littlewood--Paley decompositions where each component has Fourier support in a dyadic rectangle (see also the discussion and the references in \cite{Hy07}). Thus, even if one starts with a radial Fourier multiplier on $\RR^\nu$, the proof forces one to work with non-radial functions, because the mentioned Littlewood--Paley decomposition destroys radiality. For this reason, the known proof of the operator-valued multiplier theorem for the Fourier transform on $\RR^\nu$ does not extend straightforwardly to the case of the Hankel transform on $X_\nu$ for fractional $\nu$.

Nevertheless, in Section \ref{s:opval} we manage to prove an operator-valued $L^p$ spectral multiplier theorem, which we can use here in place of the results of \cite{StWe,We}, and may be of independent interest. Specifically, the result in Section \ref{s:opval} is a conditional result, which can be applied to any self-adjoint operator $L$ on $L^2(X)$ for some measure space $X$. For such an operator $L$, by means of the spectral theorem we can readily define an operator-valued functional calculus, whereby to any bounded and weakly measurable symbol $M : \RR \to \LinBnd(L^2(Y))$, where $Y$ is another measure space, there corresponds a bounded operator $M(L)$ on $L^2(X \times Y) \cong L^2(X;L^2(Y))$. Under the assumption that $L$ admits a (scalar-valued) H\"ormander functional calculus on $L^p(X)$ for some $p  \in (1,\infty)$, i.e.,
\[
\|F(L)\|_{L^p(X) \to L^p(X)} \lesssim \|F\|_{L^2_{s,\sloc}}
\]
for some finite $s > 1/2$, we prove in Theorem \ref{thm:5} below that $L$ also has an operator-valued functional calculus on $L^p(X\times Y) \cong L^p(X;L^p(Y))$, 
where the condition for an operator-valued spectral multiplier $M$ to define an $L^p$-bounded operator $M(L)$ is a smoothness condition of Mihlin type expressed, much as in \cite{StWe,We}, in terms of R-boundedness properties:
\[
\|M(L)\|_{L^p(X \times Y) \to L^p(X \times Y)} \lesssim \RBd_{L^p(Y)}(\{ \xi^j \partial_\xi^j M(\xi) \tc \xi \in \Rnoz, \ j=0,\dots,N \}),
\]
where $N > s+3/2$ and $\RBd_{L^p(Y)}(\cT)$ denotes the R-bound of a family of operators $\cT$ on $L^p(Y)$, see Section \ref{s:opval} for more details.

By applying our abstract result from Section \ref{s:opval} to the Bessel operator $L_\nu$, we obtain an analogue of \cite{StWe,We} for the Hankel transform. We point out the result that we obtain for the Hankel transform is somewhat more restrictive in scope than that for the Fourier transform in \cite{StWe,We}: namely, the result in \cite{StWe,We} considers operator-valued Fourier multiplier operators on $L^p(\RR^\nu;Z)$ for a broad class of Banach spaces $Z$, while we only work with $L^p(X_\nu;L^p(Y))$. Nevertheless, our result is enough for proving the required $L^p$-bounds for $p>2$ for the Riesz transform $\Riesz_0$ on $G_\nu$, as we show in Section \ref{s:adjriesz}.

The Riesz transform $\Riesz_1$, however, does not directly fall under this framework. Indeed, while the Laplacian $\Delta_\nu$ and the vector field $\vfX_0 = \partial_u$ are left-invariant on $G_\nu$, the same is not true for the vector field $\vfX_1 = e^u \partial_x$, and as a consequence the Riesz transform $\Riesz_1$ is not a convolution operator on $G_\nu$, nor is it invariant under $X_\nu$-translations. This problem is already visible at the level of the hypergroup $X_\nu$, where the Bessel operator $L_\nu$ is invariant under the hypergroup translations, but the same is not true for the first-order differential operator $\partial_x$.

A workaround to this problem can be obtained by an application of Dunkl theory for rank-one root systems. Namely, by identifying functions on the half-line $X_\nu = [0,\infty)$ with even functions on the real line $X_\nu^\rD$ equipped with the measure $\dmu^\rD(x) = 2^{-1} |x|^{\nu-1} \,\dd x$, one obtains that the Hankel transform on $X_\nu$ is the restriction to even functions of the so-called Dunkl transform on $X_\nu^\rD$. Moreover, the Dunkl operator $D_\nu$ on $X_\nu^\rD$, given by
\[
D_\nu f(x) = f'(x) + \frac{\nu-1}{2} \frac{f(x)-f(-x)}{x},
\]
when restricted to even functions $f$ coincides with the differential operator $\partial_x$; notice though that $D_\nu$ maps even functions to odd functions. Crucially, the Dunkl operator $D_\nu$ is a multiplier for the Dunkl transform; for additional details we refer to Section \ref{ss:dunkl} below. In this way, we can think of the Riesz transform $\Riesz_1 = \vfX_1 \Delta_\nu^{-1/2}$ on $L^2(G_\nu) \cong L^2(X_\nu;L^2(\RR))$ as an appropriate restriction of an operator-valued Dunkl multiplier operator on $L^2(X_\nu^\rD;L^2(\RR))$. Now, by applying the abstract multiplier theorem from Section \ref{s:opval} to the self-adjoint operator $-i D_\nu$, we obtain an operator-valued multiplier theorem for the Dunkl transform, which eventually yields the required $L^p$-bounds for the Riesz transform $\Riesz_1$ for $p>2$.

We point out that the operator-valued Dunkl multiplier theorem that we obtain here is different from that of \cite{DelKri17}, which only applies to multipliers of the form $M(\xi) = F(\xi) I$ for a scalar function $F$, but proves boundedness in $L^p(X_\nu^\rD;Z)$ for a wide class of Banach spaces $Z$. Again, our result only considers the case of $Z=L^p(Y)$, but crucially allows us to deal with truly operator-valued multipliers $M$ and is therefore applicable to the Riesz transforms on $G_\nu$.

\subsection*{Notation}
We write $\Rnon$ and $\Rpos$ for the sets $[0,\infty)$ and $(0,\infty)$ of nonnegative and positive real numbers. We also write $\NN$ for the set of natural numbers (including $0$) and $\Npos$ for $\NN \setminus \{0\}$.

We denote by $A_+$ the positive part $\max\{A,0\}$ of a real number $A$.
We write $\lfloor x \rfloor = \max\{ k \in\ZZ \tc k \leq x \}$ and $\{x\} = x-\lfloor x \rfloor$ for the integer and fractional parts of $x \in \RR$.

We denote by $\Sz(\Rnon)$ and $\Sz_\even(\Rnon)$ the spaces of the restrictions to $\Rnon$ of all Schwartz functions on $\RR$ and all even Schwartz functions respectively.

We write $\ind_S$ for the characteristic function of a set $S$, while $\#S$ denotes the number of elements of $S$.

\section{Preliminaries}\label{s:prelim}

\subsection{The hypergroups \texorpdfstring{$X_\nu$}{Xnu} and \texorpdfstring{$G_\nu$}{Gnu}}

Let $\nu\geq 1$. We recall some basic facts concerning the Bessel--Kingman hypergroup $X_\nu$ and its semidirect extension $G_\nu$, which are discussed in greater detail in \cite{MaPl}, including further references to the literature. We refer to \cite{BH,J} for a comprehensive discussion about hypergroups.

$X_\nu$ is a commutative hypergroup, with identity element $0$ and Haar measure $\dmu(x) = x^{\nu-1} \,dx$. We denote by $*_\nu$ the convolution on $X_\nu$, also known as Hankel convolution, given by
\begin{equation}\label{eq:Hconv}
\phi \ast_\nu \psi(x) = \int_{X_\nu} \phi(y) \, \tau^{[x]}_\nu \psi(y) \,\dmu(y) = \int_{X_\nu} \tau^{[x]}_\nu \phi(y) \, \psi(y) \,\dmu(y)
\end{equation}
for suitable functions $\phi,\psi$ on $X_\nu$. Here
$\tau_\nu^{[x]}$ is the Hankel translation on $X_\nu$, i.e.,
\begin{equation}\label{eq:hankel_transl}
\tau_\nu^{[x]}\phi(y) = \tau_\nu^{[y]}\phi(x) = \int_{-1}^1 \phi(\langle x,y\rangle_\omega) \,\dpr(\omega),
\end{equation}
where
\begin{equation*}
\langle x,y\rangle_\omega\defeq\sqrt{x^2+y^2-2\omega xy}
\end{equation*}
and $\dpr$ is the probability measure on $[-1,1]$ given by
\[
\dpr(\omega) = \begin{cases}
\Beta(\frac{\nu-1}{2},\frac{1}{2})^{-1} (1-\omega^2)^{(\nu-3)/2} \,\dd\omega &\text{if } \nu > 1,\\
(\dd\delta_{-1}(\omega) + \dd\delta_{1}(\omega))/2 &\text{if } \nu = 1,
\end{cases}
\]
while $\Beta$ is the Beta function, and $\delta_{-1}$ and $\delta_1$ denote the Dirac delta measures at $-1$ and $1$.

Let $j^\nu$ be the even entire function defined by
\begin{equation}\label{eq:jnu}
j^\nu(x) = 2^{\nu/2-1} \Gamma(\nu/2) \, x^{1-\nu/2} J_{\nu/2-1}(x), \quad x \in \CC.
\end{equation}
The \emph{Hankel transform} of $f\in L^1(X_\nu)$ is given by
\begin{equation}\label{eq:H_transf}
\Hank_\nu f(\xi) = \int_\RR  f(x)\, j^{\nu}_{\xi}(x) \, \dmu(x), \qquad\xi \in X_\nu,
\end{equation}
where
\[
j^{\nu}_{\xi}(x) = j^{\nu}_{x}(\xi) = j^\nu(x \xi), \qquad \xi,x \in X_\nu.
\]
It can be proved that $\Hank_\nu$ maps $L^1(X_\nu)$ into $C_0(X_\nu)$ and extends to a multiple of an isometry on $L^2(X_\nu)$, i.e., the following Plancherel formula holds (see, e.g., \cite[eq.\ (3.11)]{MaPl}):
\begin{equation}\label{eq:plancherel_Hnu}
\| \Hank_\nu f \|_{L^2(X_\nu)} = 2^{\nu/2-1} \Gamma(\nu/2) \, \| f \|_{L^2(X_\nu)}.
\end{equation}
Moreover, the Hankel transform turns Hankel convolution into multiplication, and intertwines Hankel translations and the Bessel operator with multiplication operators:
\begin{equation}\label{eq:hankel_conv_transf}
\begin{aligned}
\Hank_\nu(f *_\nu g)(\xi) &= \Hank_\nu f(\xi) \Hank_\nu g(\xi), \\
\Hank_\nu(\tau^{[x]}_\nu f)(\xi) &= j^\nu_x (\xi) \Hank_\nu f(\xi), \\
\Hank_\nu (L_\nu f)(\xi) &= \xi^2 \Hank_\nu f(\xi).
\end{aligned}
\end{equation}

The semidirect product $G_\nu = X_\nu \rtimes \RR$ is instead a noncommutative hypergroup with identity element $\bdzero = (0,0)$. 
The involution of an element $\bdx = (x,u)\in G_\nu$ is given by $\bdx^- = (e^{-u}x,-u)$. The right Haar measure on $G_\nu$ is  $\dmub(\bdx) = \dmu(x) \,\du$ and the modular function is
\begin{equation}\label{eq:modular}
\modu_\nu(x,u)\defeq e^{-\nu u}.
\end{equation}
Unless stated otherwise, all function spaces on $G_\nu$ are defined with respect to the right Haar measure.
Left and right translations on $G_\nu$ are given by
\begin{equation}\label{eq:transl}
\ell_{(y,v)} f(x,u) = r_{(x,u)}f(y,v) = \tau_\nu^{[e^{v}x]} f(y,u+v);
\end{equation}
here the Hankel translation $\tau_\nu^{[e^{v}x]}$ is meant to act on the first variable of $f$.
We denote by $\diamond_\nu$ the hypergroup convolution on $G_\nu$, given by
\begin{equation}\label{eq:diamond}
f\diamond_\nu g(\bdx)
=\int_{G_\nu} f(\bdy^{-}) \, r_{\bdx} g(\bdy)\,\dmub(\bdy)
=\int_{G_\nu} \ell_{\bdx} f(\bdy^-) \, g(\bdy) \,\dmub(\bdy)
\end{equation}
for suitable functions $f,g$ on $G_\nu$. The convolution $\diamond_\nu$ satisfies Young's inequality: for any exponents $p,q,r\in[1,\infty]$ satisfying $1+1/r=1/p+1/q$, there holds
\begin{equation}\label{eq:young}
\Vert (f \modu_\nu^{-1/q'}) \diamond_\nu g \Vert_{L^r(G_\nu)} \leq \Vert f\Vert_{L^p(G_\nu)} \Vert g\Vert_{L^q(G_\nu)},
\end{equation}
where $1/q+1/q'=1$. An analogous inequality holds for the Hankel convolution $\ast_\nu$ on $X_\nu$.

For a function $\phi : X_\nu\to\CC$, we denote its $L^1$-isometric dilations by 
\begin{equation}\label{eq:dilation}
\phi_{(\lambda)}(x)=\lambda^{-\nu} \phi(\lambda^{-1}x), \qquad \lambda>0.
\end{equation}
In the sequel, we shall make use of the relation between $\diamond_\nu$-convolution and $\ast_\nu$-convolution, namely,
\begin{equation}\label{eq:39}
f\diamond_\nu g(x,u) = \int_\RR (f^{v}\ast_\nu g_{(e^{v})}^{u-v})(x)\,\dd v,
\end{equation}
where $f^u(x)\defeq f(x,u)$. This relation encodes the fact that $G_\nu$ is a semidirect product extension of $X_\nu$ relative to the action of $\RR$ on $X_\nu$ via dilations \cite{HeyKa,Wi}.

Notice that, as a measure space, $G_\nu$ is the product of the measure spaces $X_\nu$ and $\RR$, and therefore $L^2(G_\nu)$ is the Hilbert tensor product of $L^2(X_\nu)$ and $L^2(\RR)$.
Thus, if $\HaId \defeq \Hank_\nu \otimes \id$ denotes the ``partial Hankel transform'' operator (i.e., the Hankel transform acting on the first variable of a function on $G_\nu$), then $\HaId$ is a multiple of an isometry of $L^2(G_\nu)$. Moreover, from \eqref{eq:39} and \eqref{eq:hankel_conv_transf} it follows that
\begin{equation}\label{eq:HaId_conv}
\HaId(f \diamond_\nu g)(\xi,u) = \int_\RR \HaId f(\xi,v) \, \HaId g(e^{v}\xi,u-v)\,\dd v,
\end{equation}
where we used the fact that $\Hank_\nu(g_{(\lambda)})(\xi)=\Hank_\nu g(\lambda \xi)$.

We introduce the $L^1(G_\nu)$-isometric involution $f\mapsto f^\ast$ given by
\begin{equation}\label{eq:involution}
f^\ast(\bdx)=\modu_\nu(\bdx) \overline{f(\bdx^-)}.
\end{equation}
We record here some useful relations between convolution and involution on $G_\nu$.

\begin{lm}\label{lm:conv_inv}
The following identities hold.
\begin{enumerate}[label=(\roman*)]
\item\label{en:magicdiamond} For any $f \in L^1_\loc(G_\nu)$ there holds
\begin{equation*}
\modu_\nu(\bdy) r_{\bdx} f(\bdy^-)= \modu_\nu(\bdx) \overline{r_{\bdy}f^\ast(\bdx^-)}
\end{equation*}		
for almost all $\bdx,\bdy \in G_\nu$.
\item\label{en:conv_inv} For all $g \in L^1_\loc(G_\nu)$ and all compactly supported $f,h \in L^\infty(G_\nu)$,
\begin{equation}\label{eq:inv_conv}
(f \diamond_\nu g)^* = g^* \diamond_\nu f^*, \qquad \langle f \diamond_\nu g , h \rangle_{L^2(G_\nu)} = \langle f , h \diamond_\nu g^* \rangle_{L^2(G_\nu)}.
\end{equation}
\end{enumerate}
\end{lm}
\begin{proof}
\ref{en:magicdiamond}.
Let $\bdx=(x,u)$ and $\bdy=(y,v)$. By the definition  \eqref{eq:transl} of the right translations,
\[\begin{split}
		\modu_\nu(\bdx) \overline{r_{\bdy}f^\ast(\bdx^- )} 
		&=e^{-\nu u} \overline{\tau_\nu^{[e^{-u}y]}f^\ast(e^{-u}x,v-u)} \\
		&= e^{-\nu u} \int_{-1}^1 \overline{f^\ast(e^{-u}\langle x,y\rangle_\omega,v-u)} \,\dpr(\omega)\\
		&= e^{-\nu v} \int_{-1}^1 f(e^{-v}\langle x,y\rangle_\omega,u-v) \,\dpr(\omega)\\
		&= e^{-\nu v} \tau_\nu^{[e^{-v}x]}f(e^{-v}y,u-v)\\
		&=\modu_\nu(\bdy) r_{\bdx}f(\bdy^-).
\end{split}\]

\ref{en:conv_inv}. Notice that, by \eqref{eq:involution}, \eqref{eq:transl} and \eqref{eq:diamond},
\[
(f \diamond_\nu g)^*(\bdx) = \modu_\nu(\bdx) \, \overline{f \diamond_\nu g(\bdx^-)} = \modu_\nu(\bdx) \int_{G_\nu} \overline{r_{\bdy^-} f(\bdx^-)} \, \overline{g(\bdy)} \,\dmub(\bdy).
\]
Thus, by part \ref{en:magicdiamond}, and again \eqref{eq:involution} and \eqref{eq:diamond},
\begin{multline*}
(f \diamond_\nu g)^*(\bdx) = \int_{G_\nu} \modu_\nu(\bdy^-) r_{\bdx} f^*(\bdy) \, \overline{g(\bdy)} \,\dmub(\bdy) \\
= \int_{G_\nu} r_{\bdx} f^*(\bdy) \, g^*(\bdy^-) \,\dmub(\bdy) = g^* \diamond_\nu f^*(\bdx),
\end{multline*}
as claimed.

Similarly,
\[\begin{split}
\langle f \diamond_\nu g , h \rangle_{L^2(G_\nu)} 
&= \int_{G_\nu} \int_{G_\nu} f(\bdy^{-}) \, r_{\bdx} g(\bdy)\,\dmub(\bdy) \, \overline{h(\bdx)} \,\dmub(\bdx) \\
&= \int_{G_\nu} \int_{G_\nu} f(\bdy) \, r_{\bdx} g(\bdy^-) \, \modu_\nu(\bdy) \,\dmub(\bdy) \, \overline{h(\bdx)} \,\dmub(\bdx) \\
&= \int_{G_\nu} \int_{G_\nu} f(\bdy) \, \overline{r_{\bdy} g^*(\bdx^-)} \, \modu_\nu(\bdx) \,\dmub(\bdy) \, \overline{h(\bdx)} \,\dmub(\bdx) \\
&= \int_{G_\nu} f(\bdy) \overline{\int_{G_\nu} r_{\bdy} g^*(\bdx) \, h(\bdx^-) \,\dmub(\bdx)} \,\dmub(\bdy) = \langle f, h \diamond_\nu g^* \rangle_{L^2(G_\nu)},
\end{split}\]
where part \ref{en:magicdiamond}, equations \eqref{eq:involution} and \eqref{eq:diamond}, and Fubini's theorem were used.
\end{proof}

\subsection{Riemannian structure and Laplacian on \texorpdfstring{$G_\nu$}{Gnu}}

Recall from \cite[eq.\ (4.9)]{MaPl} that we can write $\Delta_\nu$ in divergence form, i.e., $\Delta_\nu = \nabla_\nu^+ \nabla_\nu$, where $\nabla_\nu$ is the gradient associated with $\Delta_\nu$, given by
\[
\nabla_{\nu} =  (\vfX_0,\vfX_1),
\]
and $\nabla^+$ is its formal adjoint with respect to $\dmub$. Here $\vfX_0$ and $\vfX_1$ are the vector fields of \eqref{eq:vfs}. Thus
\begin{equation*}
\vert\nabla_\nu f\vert^2 = |\vfX_0 f|^2+  |\vfX_1 f|^2.
\end{equation*}
Sometimes we use the upper script to emphasize on which variables the gradient acts, e.g., we may write $\nabla_\nu^{\bdx}$ for the gradient in $\bdx$. We point out that the vector fields $\vfX_0$ and $\vfX_1$ and the gradient $\nabla_\nu$ are independent of $\nu$, while the adjoint $\nabla_\nu^+$ and the Laplacian $\Delta_\nu$ depend on $\nu$.

$G_\nu$ is naturally equipped with a Riemannian structure that makes the vector fields $\vfX_0,\vfX_1$ an orthonormal frame. As discussed in \cite[Section 4.4]{MaPl}, this structure does not depend on $\nu$ and is the restriction to $G_\nu = \Rnon \times \RR$ of a suitable realisation of the hyperbolic plane Riemannian structure on $\RR \times \RR$; moreover $G_\nu$ is geodesically convex in $\RR \times \RR$, and therefore any pair of points of $G_\nu$ is joined by a single length-minimizing geodesic. An explicit formula for the Riemannian distance $\dist$ on $G_\nu$ is also available:
\begin{equation}\label{eq:dist}
\dist((x,u),(x',u'))=\arccosh\left(\cosh(u-u')+\frac{|x-x'|^2}{2e^{u+u'}} \right).
\end{equation}
We shall write $B_{G_\nu}(\bdx,r)$ and $\overline{B}_{G_\nu}(\bdx,r)$ for the open and closed $\dist$-balls centred at $\bdx \in G_\nu$ of radius $r > 0$.
We also define $|\bdx|_\dist \defeq \dist(\bdx,\bdzero)$. Notice that $|\bdx^-|_\dist = |\bdx|_\dist$ for all $\bdx \in G_\nu$.
We shall call \emph{radial} any function on $G_\nu$ that depends only on the distance $|\cdot|_\dist$ of its argument from the identity element $\bdzero$.

The next lemma serves as an analogue of the mean value theorem on $G_\nu$.

\begin{lm}\label{lm:10}
Fix $p\in[1,\infty]$ and $\bdx,\bdx'\in G_\nu$ and let $\gamma$ be the geodesic between them. If $\iker : G_\nu\times G_\nu\to\CC$ is continuously differentiable with respect to the first variable, then
\[
\left\Vert \iker(\bdx,\cdot) - \iker(\bdx',\cdot)\right\Vert_{L^p(G_\nu)}
\leq \dist(\bdx,\bdx') \sup_{\bdy\in\gamma} \left\Vert \nabla_\nu^{\bdy} \iker(\bdy,\cdot)\right\Vert_{L^p(G_\nu)}.
\]
\end{lm}
\begin{proof}
Fix $\bdx,\bdx'\in G_\nu$ and let $\gamma=(\gamma_1,\gamma_2) : [0,L]\to G_\nu$ be the geodesic between $\bdx$ and $\bdx'$ parameterized according to the arc length so that $L=\dist(\bdx,\bdx')$. Thus,
\begin{equation*}
	(e^{-\gamma_2(t)} \gamma_1'(t))^2 + (\gamma'_2(t))^2\leq 1,\qquad t \in [0,L].
\end{equation*}
Hence, for all $\bdz \in G_\nu$, the chain rule and the Cauchy--Schwarz inequality yield
\[
\begin{split}
| \partial_t \iker(\gamma(t),\bdz) | 
&= \left| e^{-\gamma_2(t)}\gamma_1'(t) \vfX_1^{\bdy} \iker(\bdy,\bdz)|_{\bdy=\gamma(t)}  +  \gamma_2'(t) \vfX_0^\bdy \iker(\bdy,\bdz)|_{\bdy=\gamma(t)} \right|\\
&\leq \left|\nabla_\nu^{\bdy} \iker(\bdy,\bdz)|_{\bdy=\gamma(t)} \right|.
\end{split}
\]
So, by the fundamental theorem of calculus and Minkowski's integral inequality,
\begin{multline*}
	\left\Vert \iker(\bdx,\cdot) - \iker(\bdx,\cdot)\right\Vert_{L^p(G_\nu)}\\
	=\left\Vert \int_0^L \partial_t \iker(\gamma(t),\cdot)\,\dd t\right\Vert_{L^p(G_\nu)}
	\leq L \sup_{\bdy\in\gamma}\left\Vert|\nabla^{\bdy}_\nu \iker(\bdy,\cdot)|\right\Vert_{L^p(G_\nu)},
\end{multline*}
as desired.
\end{proof}

We now collect some useful observations relating distance, derivatives and translations on $X_\nu$ and $G_\nu$. Much as before, by $C^1_\even(X_\nu)$ we denote the space of restrictions to $X_\nu = [0,\infty)$ of even $C^1$ functions on $\RR$; similarly, by $C^1_\even(G_\nu)$ we denote the space of restrictions to $G_\nu = X_\nu \times \RR$ of $C^1$ functions on $\RR^2$ which are even in the first variable.

\begin{lm}\label{lm:11}
Let $x\in X_\nu$ and $\bdx,\bdy\in G_\nu$.
\begin{enumerate}[label=(\roman*)]	
\item\label{lm:11(1)} Let $g \in C(G_\nu)$. Then
\[
 \supp (r_\bdx g) \subseteq \{ \bdz \in G_\nu \tc \dist(\supp g,\bdz) \leq |\bdx|_\dist\}.
\]
\item\label{lm:11(2)} Let $f\in C^1_\even(X_\nu)$ and $g\in C^1_\even(G_\nu)$. Then $(x,y) \mapsto \tau_\nu^{[x]} f(y)$ is in $C^1(X_\nu \times X_\nu)$, $(\bdx,\bdy) \mapsto r_{\bdx} f(\bdy)$ is in $C^1(G_\nu \times G_\nu)$, and moreover
\begin{equation}\label{eq:17}
|\partial_x \tau_{\nu}^{[x]} f | \leq \tau_\nu^{[x]} |f'|\qquad \text{and} \qquad
| \nabla_{\nu}^{\bdx} r_{\bdx} g | \leq r_{\bdx} |\nabla_\nu g|
\end{equation}
pointwise.
\item\label{lm:11(3)} Let $G : \Rnon \to \Rnon$ be an increasing function and define $g = G( \dist(\cdot,\bdy))$. Then
\begin{equation*}
\ell_{\bdy} g(\bdx) \leq G(|\bdx|_\dist).
\end{equation*}
\end{enumerate}
\end{lm}
\begin{proof}
\ref{lm:11(1)}. Let $\bdx = (x,u)$ and $\bdz = (z,w)$. From \eqref{eq:transl} and \eqref{eq:hankel_transl}, we see that $r_{\bdx} g(\bdz)$ only depends on the values of $g$ on the set
\[
[|z-e^w x|,z+e^w x] \times \{u+w\}.
\]
and vanishes if $g$ vanishes there. Thus, if $r_{\bdx} g(\bdz) \neq 0$, then
\[
([|z-e^w x|,z+e^w x] \times \{u+w\}) \cap \supp g \neq \emptyset.
\]
On the other hand, if $\bdy = (y,v) \in [|z-e^w x|,z+e^w x] \times \{u+w\}$, then $v = u+w$ and $|y-z| \leq e^w x$; so, by \eqref{eq:dist},
\[
\dist(\bdy,\bdz) = \arccosh(\cosh(v-w) + |y-z|^2/(2e^{v+w})) \leq |\bdx|_\dist.
\]
Thus, if $r_{\bdx} g(\bdz) \neq 0$, then $\dist(\supp g,\bdz) \leq |\bdx|_\dist$, as desired.

\ref{lm:11(2)}.
Let us first discuss the result for $f \in C^1_\even(X_\nu)$.
By assumption, $f$ extends to an even $C^1$ function on $\RR$, which we shall denote by the same symbol $f$.
If $\nu = 1$, then \eqref{eq:hankel_transl} reduces to
\[
\tau_1^{[x]} f(y) = \frac{1}{2} (f(x+y)+f(x-y)),
\]
whence both the $C^1$ regularity of $(x,y) \mapsto \tau_1^{[x]} f(y)$ and the desired estimate follow. If $\nu > 1$, then, $\langle x,y\rangle_\omega > 0$ for all $(x,y) \neq (0,0)$ and $\omega \in (-1,1)$; by \eqref{eq:hankel_transl} and differentiation under the integral sign, 
also using the estimates
\begin{equation}\label{eq:bd_xycos}
|x-y\omega|, |y-x\omega| \leq\langle x,y\rangle_\omega,
\end{equation}
we then obtain that
\begin{align*}
\partial_x \tau_{\nu}^{[x]} f(y) & = \int_{-1}^1 f'(\langle x,y\rangle_\omega) \frac{x-y \omega}{\langle x,y\rangle_\omega} \,\dpr(\omega),\\
\partial_y \tau_{\nu}^{[x]} f(y) & = \int_{-1}^1 f'(\langle x,y\rangle_\omega) \frac{y-x \omega}{\langle x,y\rangle_\omega} \,\dpr(\omega)
\end{align*}
for all $(x,y) \neq (0,0)$, thus showing the $C^1$ regularity of $(x,y) \mapsto \tau_\nu^{[x]} f(y)$ on $X_\nu \times X_\nu \setminus \{(0,0)\}$. As $f'(0) = 0$, from the above formulas we also see that
\[
\lim_{(x,y) \to (0,0)} \partial_x \tau_{\nu}^{[x]} f(y) = \lim_{(x,y) \to (0,0)} \partial_y \tau_{\nu}^{[x]} f(y) = 0,
\]
thus proving that $(x,y) \mapsto \tau_\nu^{[x]} f(y)$ is actually $C^1$ on the whole $X_\nu \times X_\nu$. Moreover, the above formulas and \eqref{eq:bd_xycos} show that
\begin{equation*}
|\partial_x \tau_{\nu}^{[x]} f(y) | 
= \left|\int_{-1}^1 f'(\langle x,y\rangle_\omega) \frac{x-y \omega}{\langle x,y\rangle_\omega} \,\dpr(\omega)\right|
\leq \tau_\nu^{[x]}|f'|(y),
\end{equation*}
as desired.

The $C^1$ regularity of $(\bdx,\bdy) \mapsto r_\bdx g(\bdy)$ on $G_\nu \times G_\nu$ follows from the definition \eqref{eq:transl} of the right translations and the just proved $C^1$ regularity result for Hankel translations. 		
In order to justify the second inequality in \eqref{eq:17} observe that, by the definition of the right translations \eqref{eq:transl},
\begin{equation}\label{eq:dutransl}
| \partial_u r_{(x,u)}g(y,v)| 
= | \tau_\nu^{[e^v x]} (\partial_u g)(y,u+v)| = |r_{(x,u)} \vfX_0 g(y,v)| \leq r_{(x,u)} |\vfX_0 g(y,v)|
\end{equation}
and, also by the first bound in \eqref{eq:17},
\[\begin{split}
| e^u \partial_x r_{(x,u)}g(y,v)| 
= | e^{u+v} (\partial_z \tau_\nu^{[z]}g)(y,u+v)|_{z=e^v x} |
&\leq e^{u+v} \tau_\nu^{[e^v x]} | \partial_y g|(y,u+v)\\
&= r_{(x,u)} | \vfX_1 g(y,v)|.
\end{split}\]

Therefore,
\begin{equation}\label{eq:18}
\begin{split}
|\nabla_\nu^{(x,u)} r_{(x,u)}g  | 
&= \sup_{|a_1|^2+|a_2|^2\leq 1} | a_1 \partial_u r_{(x,u)} g + a_2e^{u}\partial_x r_{(x,u)} g|\\
&\leq \sup_{|a_1|^2+|a_2|^2\leq 1} r_{(x,u)} ( |a_1| |\vfX_0 g| + |a_2|  |\vfX_1 g|)\\
&\leq r_{(x,u)} |\nabla_\nu g|,
\end{split}
\end{equation}
where in the last inequality we used the monotonicity of $r_{(x,u)}$. This proves \eqref{eq:17}.
		
\ref{lm:11(3)}.
Fix $\bdy=(y,v)\in G_\nu $. By the definition of the left translations \eqref{eq:transl} and the formula \eqref{eq:dist} for the distance on $G_\nu$,
\[\begin{split}
\ell_{(y,v)} g(x,u) &= \tau_\nu^{[y]}g(e^v x,u+v)\\
&=\int_{-1}^1 G\left(\arccosh\left[\cosh u +\frac{(\langle e^v x,y\rangle_\omega -y)^2}{2e^{u+2v}} \right]\right) \,\dpr(\omega).
\end{split}\]
Notice that
\begin{equation*}
|\langle e^v x,y\rangle_\omega -y|\leq e^v x,\qquad \omega\in[-1,1].
\end{equation*}
Thus, by the monotonicity of $G$,
\begin{equation*}
\ell_{(y,v)} g(x,u) \leq \int_{-1}^1 G\left(\arccosh\left[\cosh u+\frac{e^{2v}x^2}{2e^{u+2v}} \right]\right) \,\dpr(\omega) = G(|(x,u)|_\dist),
\end{equation*}
as desired.
\end{proof}

We now recall the Plancherel theorem for the functional calculus of $\Delta_\nu$ \cite[Proposition~4.12 and Corollary~4.16]{MaPl}. 
For $\lambda > 0$, $a,b \in \RR$, we write
\begin{equation}\label{eq:doublepower}
\lambda^{[a,b]} \defeq \begin{cases}
\lambda^a &\text{if }\lambda \leq 1,\\
\lambda^b &\text{if }\lambda \geq 1.
\end{cases}
\end{equation}

\begin{prop}\label{prop:6}
If $F : \RR_+\to\CC$ is a bounded Borel function that belongs to $L^2(\RR_+,\lambda^{[3/2,(\nu+1)/2]}\,\ddlam)$, then $F(\Delta_\nu)$ is a right $\diamond_\nu$-convolution operator, i.e.,
\[
F(\Delta_\nu) f = f \diamond_\nu K_{F(\Delta_\nu)},	
\]
whose convolution kernel $K_{F(\Delta_\nu)}$ is in $L^2(G_\nu)$. Moreover, there exists a regular Borel measure $\sigma_\nu$ on $\Rnon$, called the \emph{Plancherel measure} associated with $\Delta_\nu$, such that
\begin{equation}\label{eq:11}
	 \|  K_{F(\Delta_\nu)} \|_{L^2(G_\nu)}^2
	= \int_0^\infty |F|^2\, \dd\sigma_\nu 
	\lesssim_\nu \int_0^\infty |F(\lambda)|^2 \,\lambda^{[3/2,(\nu+1)/2]}\, \ddlam
\end{equation}
for all such $F$.
In addition, $\modu_\nu^{-1/2} K_{F(\Delta_\nu)}$ is a radial function.
\end{prop}

From the above Plancherel formula we can easily derive a sort of \emph{Riemann--Lebesgue Lemma} for the mapping $F\mapsto K_{F(\Delta_\nu)}$.

\begin{lm}\label{lm:1}
Let $F\in L^1(\RR_+,\lambda^{[3/2,(\nu+1)/2]}\,\ddlam)$ be bounded. Then $\modu_\nu^{-1/2} K_{F(\Delta_\nu)} \in C_0(G_\nu)$ and
\begin{equation*}
\Vert \modu_\nu^{-1/2} K_{F(\Delta_\nu)}\Vert_{\infty} \leq \int_0^\infty |F| \,\dd\sigma_\nu \lesssim_\nu \int_0^\infty |F(\lambda)| \, \lambda^{[3/2,(\nu+1)/2]}\,\ddlam.
\end{equation*}
\end{lm}
\begin{proof}
Let us decompose $F=F_1F_2$ in such a way that $|F_1|=|F_2|=|F|^{1/2}$. Since $F_1,F_2\in L^2(\RR_+,\lambda^{[3/2,(\nu+1)/2]}\,\ddlam)$, by Proposition \ref{prop:6} the functions $\modu_\nu^{-1/2}K_{F_i(\Delta_\nu)}$, $i=1,2$, are radial on $G_\nu$, and
\begin{equation*}
	\Vert \modu_\nu^{-1/2}K_{F_i(\Delta_\nu)}\Vert_{L^2}^2=\Vert K_{F_i(\Delta_\nu)}\Vert_{L^2}^2 = \int_0^\infty |F_i|^2 \,\dd\sigma_\nu = \int_0^\infty |F| \,\dd\sigma_\nu
\end{equation*}
(the first equality follows from radiality, cf.\ \cite[Proposition~4.15]{MaPl}). In addition, as $F(\Delta_\nu) = F_1(\Delta_\nu) F_2(\Delta_\nu)$,
\[\begin{split}
\modu_\nu^{-1/2} K_{F(\Delta_\nu)} 
&= \modu_\nu^{-1/2}(K_{F_1(\Delta_\nu)}\diamond_\nu K_{F_2(\Delta_\nu)}) \\
&= (\modu_\nu^{-1/2}K_{F_1(\Delta_\nu)})\diamond_\nu (\modu_\nu^{-1/2}K_{F_2(\Delta_\nu)})
\end{split}\]
by the properties of the modular function \cite[Section 5.3]{J}.
Hence, by Young's convolution inequality \eqref{eq:young} and \cite[(5.5P)]{J}, $\modu_\nu^{-1/2} K_{F(\Delta_\nu)} \in C_0(G_\nu)$ and
\[\begin{split}
\Vert \modu_\nu^{-1/2}K_{F(\Delta_\nu)}\Vert_{\infty}
&\leq \Vert K_{F_1(\Delta_\nu)}\Vert_{L^2} \Vert \modu_\nu^{-1/2}K_{F_2(\Delta_\nu)}\Vert_{L^2} = \int_0^\infty |F| \,\dd\sigma_\nu.
\end{split}\]
The bound \eqref{eq:11} on the Plancherel measure $\sigma_\nu$ concludes the proof.
\end{proof}

We record here a few integration formulas for functions of the form $\modu_\nu^{1/2} f$, where $f$ is radial, against certain weights on $G_\nu$; the formulas use the notation \eqref{eq:doublepower}. The proof is fully analogous to that of \cite[Lemma 3.3]{Ma23} and is omitted.

\begin{lm}\label{lm:13}
Let $k\in\Rpos$ and let $f : \RR \to [0,\infty)$ be measurable. Then
\begin{align*}
	\int_{G_\nu} \modu_\nu^{1/2}(\bdx) \, f(r)\,\dmub(\bdx) 
	&\simeq_{k,\nu} \int_0^\infty f(s) \, s^{[\nu,1]} \, e^{\nu s/2}\,\dd s,\\
	\int_{G_\nu} x^k \, \modu_\nu^{1/2}(\bdx) \, f(r)\,\dmub(\bdx) 
	&\simeq_{k,\nu} \int_0^\infty f(s) \, s^{[\nu+k,0]} \, e^{(2k+\nu)s/2}\,\dd s,\\
	\int_{G_\nu} \ind_{\{u\leq 1\}} \, x^k \, \modu_\nu^{1/2}(\bdx) \, f(r) \,\dmub(\bdx) 
	&\simeq_{k,\nu} \int_0^\infty f(s) \, s^{[\nu+k,0]} \, e^{(k+\nu)s/2}\,\dd s,\\
	\int_{G_\nu} \left|\sinh u\right|^k  \modu_\nu^{1/2}(\bdx) \, f(r)\,\dmub(\bdx) 
	&\simeq_{k,\nu} \int_0^\infty f(s) \, s^{[\nu+k,0]} \, e^{(2k+\nu)s/2}\,\dd s,\\
	\int_{G_\nu} \ind_{\{|u|\leq 1\}} \, |u|^k \, \modu_\nu^{1/2}(\bdx) \, f(r)\,\dmub(\bdx) 
	&\simeq_{k,\nu} \int_0^\infty f(s) \, s^{[\nu+k,0]} \, e^{\nu s/2}\,\dd s,
\end{align*}
where $\bdx = (x,u)$ and $r = |\bdx|_\dist$.
\end{lm}

Notice that, if $T$ is a right $\diamond_\nu$-convolution operator on $L^2(G_\nu)$ with convolution kernel $K_T$, then $T$ is an integral operator with integral kernel $\iker_T$ given by
\begin{equation}\label{eq:26}
\iker_T( \bdx,\bdy) = \modu_\nu(\bdy) r_{\bdx} K_T(\bdy^-) = \modu_\nu(\bdx) \overline{r_{\bdy}K_T^\ast(\bdx^-)};
\end{equation}
the former equality follows from \eqref{eq:diamond} via the change of variables $\bdy \mapsto \bdy^-$ (see also \cite[eq.\ (4.38)]{MaPl}), while the latter follows from Lemma \ref{lm:conv_inv}\ref{en:magicdiamond}.

We now show that, in the case of $\diamond_\nu$-convolution operators, certain weighted $L^1$-estimates for the integral kernels can be reduced to corresponding estimates for the convolution kernels.

\begin{lm}\label{lm:l1est_int_conv}
Assume that the functions $K_T : G_\nu \to \CC$ and $\iker_T : G_\nu \times G_\nu \to \CC$ are related by \eqref{eq:26}. Then, for any increasing function $w : \Rnon \to \Rnon$,
\begin{align*}
\sup_{\bdy\in G_\nu} \int_{G_\nu} |\iker_T(\bdx,\bdy)| \, w(\dist(\bdx,\bdy)) \,\dmub(\bdx) 
&\leq \int_{G_\nu} |K_T(\bdx)| \, w(|\bdx|_\dist) \,\dmub(\bdx),\\
\sup_{\bdy\in G_\nu} \int_{G_\nu} |\nabla^\bdy_\nu \iker_T(\bdx,\bdy)| \, w(\dist(\bdx,\bdy)) \,\dmub(\bdx)
&\leq \int_{G_\nu} |\nabla_\nu K_T^*(\bdx)| \, w(|\bdx|_\dist) \,\dmub(\bdx), \\
\sup_{\bdy\in G_\nu} \int_{G_\nu} |\nabla^\bdx_\nu \iker_T(\bdx,\bdy)| \, w(\dist(\bdx,\bdy)) \,\dmub(\bdx) 
&\leq \int_{G_\nu} |\nabla_\nu K_T(\bdx)| \, w(|\bdx|_\dist) \,\dmub(\bdx).
\end{align*}
\end{lm}
\begin{proof}
Fix $\bdy\in G_\nu$ and denote $g(\bdx)\defeq w(\dist(\bdx,\bdy))$.
Then, by \eqref{eq:26}, \eqref{eq:diamond} and Lemma \ref{lm:11}\ref{lm:11(3)},
\begin{multline}\label{eq:dist_conv_est}
\int_{G_\nu} |\iker_T(\bdx,\bdy)| \, w(\dist(\bdx,\bdy)) \,\dmub(\bdx) 
	= \int_{G_\nu} | r_{\bdy} K^\ast_T(\bdx)| \, g(\bdx^-) \,\dmub(\bdx)
	\leq g\diamond_\nu |K_T^\ast|(\bdy) \\
	= \int_{G_\nu} | K^\ast_T(\bdx)| \, \ell_{\bdy}g(\bdx^-)\,\dmub(\bdx)
	\leq \int_{G_\nu} | K_T(\bdx)| \, w(|\bdx|_\dist) \,\dmub(\bdx).
\end{multline}
Moreover, by \eqref{eq:26} and Lemma \ref{lm:11}\ref{lm:11(2)} we deduce that
\begin{multline*}
\int_{G_\nu} |\nabla_\nu^\bdy \iker_T(\bdx,\bdy)| \, w(\dist(\bdx,\bdy)) \,\dmub(\bdx) =
	\int_{G_\nu} | \nabla_\nu^\bdy r_{\bdy}K^\ast_T(\bdx)| \, g(\bdx^{-}) \,\dmub(\bdx)  \\
	\leq  \int_{G_\nu}  r_{\bdy} | \nabla_\nu K^\ast_T|(\bdx) \, g(\bdx^{-}) \,\dmub(\bdx)
	\leq \int_{G_\nu} |\nabla_\nu K^\ast_T(\bdx)| \, w(|\bdx|_\dist) \,\dmub(\bdx), 
\end{multline*}
where the last inequality is proved much as in \eqref{eq:dist_conv_est}. Similarly,
\begin{multline*}
\int_{G_\nu} |\nabla_\nu^\bdx \iker_T(\bdx,\bdy)| \, w(\dist(\bdx,\bdy)) \,\dmub(\bdx) =
	\modu_\nu(\bdy) \int_{G_\nu} | \nabla_\nu^\bdx r_{\bdx} K_T(\bdy^{-})| \, g(\bdx) \,\dmub(\bdx)  \\
	\leq \modu_\nu(\bdy) \int_{G_\nu}  r_{\bdx} | \nabla_\nu K_T| (\bdy^{-}) \, g(\bdx) \,\dmub(\bdx)
	= \int_{G_\nu} r_{\bdy} | (\nabla_\nu K_T)^*| (\bdx) \, g(\bdx^{-}) \,\dmub(\bdx),
\end{multline*}
where the last equality is due to Lemma \ref{lm:conv_inv}\ref{en:conv_inv} and the change of variables $\bdx \mapsto \bdx^{-}$; we can then proceed much as in  \eqref{eq:dist_conv_est} to bound the last integral by
\[
\int_{G_\nu} | \nabla_\nu K_T(\bdx)| \, w(|\bdx|_\dist) \,\dmub(\bdx),
\]
as desired.
\end{proof}

We conclude this section with a technical result dealing with certain singular $\diamond_\nu$-convolution kernels. The rough idea is that, even if a function on $G_\nu$ is singular (not locally integrable) at $\bdzero$, one can still use it as a $\diamond_\nu$-convolution kernel, at least in the sense of the ``off-diagonal kernels'' of singular integral operators.

\begin{prop}\label{prop:SIO_diamond}
Let $K \in L^1_\loc(G_\nu) \cap C^1_\even(G_\nu \setminus \{\bdzero\})$. Then, for all $f \in C_c(G_\nu)$, the function $f \diamond_\nu K$ is continuosly differentiable on $G_\nu \setminus \supp f$ and
\begin{equation}\label{eq:SIO_der}
\vfX_0(f \diamond_\nu K)(\bdx) = \int_{G_\nu} \ell_{\bdx} f(\bdy^-) \, \vfX_0 K(\bdy) \,\dmub(\bdy) \quad\text{for all } \bdx \notin \supp f.
\end{equation}
Moreover, for all $f,g \in C_c(G_\nu)$ with disjoint supports,
\begin{equation}\label{eq:SIO_der_adj}
\begin{split}
\langle \vfX_0(f \diamond_\nu K), g \rangle_{L^2(G_\nu)} 
&= \int_{G_\nu} \overline{g(\bdx)} \int_{G_\nu} \ell_{\bdx} f(\bdy^{-}) \, \vfX_0 K(\bdy) \,\dmub(\bdy) \,\dmub(\bdx) \\
&= \int_{G_\nu} f(\bdx) \overline{\int_{G_\nu} \ell_{\bdx} g(\bdy^{-}) \, (\vfX_0 K)^*(\bdy) \,\dmub(\bdy)} \,\dmub(\bdx).
\end{split}
\end{equation}
\end{prop}
\begin{proof}
Notice first that, by \eqref{eq:diamond}, 
\[
f \diamond_\nu K(\bdx) = \int_{G_\nu} \ell_{\bdx} f(\bdy^-) \, K(\bdy) \,\dmub(\bdy) .
\]
Moreover, by Lemma \ref{lm:11}\ref{lm:11(1)},
\[
\supp (r_{\bdy} f) \subseteq \{ \bdz \in G_\nu \tc \dist(\supp f,\bdz) \leq |\bdy|_\dist \};
\]
thus for any fixed $\varepsilon > 0$ we see that, if $\dist(\supp f,\bdx) > \varepsilon$, then $\ell_\bdx f(\bdy) = r_{\bdy} f(\bdx) = 0$ whenever $|\bdy|_\dist \leq \varepsilon$. Therefore, if $\Omega_\varepsilon = \{ \bdx \tc \dist(\supp f,\bdx) > \varepsilon \}$, then
\[
f \diamond_\nu K(\bdx) = \int_{G_\nu} \ell_{\bdx} f(\bdy^-) \, K_\varepsilon(\bdy) \,\dmub(\bdy) = \int_{G_\nu} f(\bdy^-) \, \ell_{\bdy} K_\varepsilon(\bdx) \,\dmub(\bdy) \quad\forall \bdx \in \Omega_\varepsilon,
\]
where $K_\varepsilon \in C^1_\even(G_\nu)$ is any function that coincides with $K$ on $G_\nu \setminus \overline{B}_{G_\nu}(\bdzero,\varepsilon)$ and vanishes in a neighbourhood of $\bdzero$.

As $(\bdx,\bdy) \mapsto \ell_{\bdy} K_\varepsilon(\bdx) = r_{\bdx} K_\varepsilon(\bdy)$ is of class $C^1$ on $X_\nu \times X_\nu$ (see Lemma \ref{lm:11}\ref{lm:11(2)}), from the above expression and differentiation under the integral sign we deduce that $f \diamond_\nu K$ is continuously differentiable on $\Omega_\varepsilon$. Moreover, as $\vfX_0$ commutes with left translations, i.e.,
\begin{equation}\label{eq:X0_li}
\vfX_0 \ell_\bdx g = \ell_\bdx \vfX_0 g \qquad \forall \bdx \in G_\nu, \ g \in C^1(G_\nu),
\end{equation}
(see \eqref{eq:vfs} and \eqref{eq:transl}), we see that, for all $\bdx \in \Omega_\varepsilon$,
\[
\vfX_0(f \diamond_\nu K)(\bdx) 
= \int_{G_\nu} f(\bdy^-) \, \ell_{\bdy} \vfX_0 K_\varepsilon(\bdx) \,\dmub(\bdy)
= \int_{G_\nu} \ell_{\bdx} f(\bdy^-) \, \vfX_0 K_\varepsilon(\bdy) \,\dmub(\bdy).
\]
As $\vfX_0$ is a differential operator, $\vfX_0 K_\varepsilon$ coincides with $\vfX_0 K$ off $\overline{B}_{G_\nu}(\bdzero,\varepsilon)$. This proves the identity \eqref{eq:SIO_der} for all $\bdx \in \Omega_\varepsilon$; as $\varepsilon > 0$ was arbitrary we obtain \eqref{eq:SIO_der} for all $\bdx \notin \supp f$.

Now, if $g \in C_c(G_\nu)$ has disjoint support from $\supp f$, then we can find $\varepsilon>0$ such that $\supp g \subseteq \Omega_{2\varepsilon} \subseteq \Omega_\varepsilon$. In particular, the previous expression for $\vfX_0(f \diamond_\nu K)$ on $\Omega_\varepsilon$ applies and we obtain
\[\begin{split}
\langle \vfX_0(f \diamond_\nu K), g \rangle_{L^2(G_\nu)} 
&= \int_{G_\nu} \overline{g(\bdx)} \int_{G_\nu} \ell_{\bdx} f(\bdy^{-}) \, \vfX_0 K_\varepsilon(\bdy) \,\dmub(\bdy) \,\dmub(\bdx) \\
&= \langle f \diamond_\nu (\vfX_0 K_\varepsilon), g \rangle_{L^2(G_\nu)}.
\end{split}\]
As $\vfX_0 K_\varepsilon \in L^1_\loc(G_\nu)$, while $f,g \in C_c(G_\nu)$, by Lemma \ref{lm:conv_inv}\ref{en:conv_inv} we deduce
\[\begin{split}
\langle \vfX_0(f \diamond_\nu K), g \rangle_{L^2(G_\nu)}
&= \langle f , g \diamond_\nu (\vfX_0 K_\varepsilon)^* \rangle_{L^2(G_\nu)} \\
&= \int_{G_\nu} f(\bdx) \overline{\int_{G_\nu} \ell_{\bdx} g(\bdy^{-}) \, (\vfX_0 K_\varepsilon)^*(\bdy) \,\dmub(\bdy)} \,\dmub(\bdx).
\end{split}\]
On the other hand, from the fact that $\supp g \subseteq \Omega_{2\varepsilon} = \{ \bdx \tc \dist(\supp f,\bdx) > 2\varepsilon\}$, we deduce that $\supp f \subseteq \{ \bdx \tc \dist(\supp g,\bdx) > \varepsilon\}$, and moreover, arguing as before, $\ell_{\bdx} g(\bdy^{-})$ vanishes if $|\bdy| \leq \varepsilon$ and $\dist(\supp g,\bdx) > \varepsilon$. As $(\vfX_0 K_\varepsilon)^* = (\vfX_0 K)^*$ off $\overline{B}_{G_\nu}(\bdzero,\varepsilon)$, we can replace $(\vfX_0 K_\varepsilon)^*$ with $(\vfX_0 K)^*$ in the last integral and complete the proof of \eqref{eq:SIO_der_adj}.
\end{proof}

In the above proof the left-invariance property \eqref{eq:X0_li} of $\vfX_0$ is crucial, and indeed the same result does not hold for $\vfX_1$. This is why, while the previous proposition will be enough to discuss the Riesz transform $\Riesz_0 = \vfX_0 \Delta^{-1/2}$ as a singular integral operator, for the Riesz transform $\Riesz_1 = \vfX_1 \Delta^{-1/2}$ we will need a different approach, based on the Dunkl theory discussed in Section \ref{ss:dunkl} below.

\subsection{Kernel formulas}
Recall from Proposition \ref{prop:6} that, for any bounded function $F\in L^2(\Rnon,\lambda^{[3/2,(\nu+1)/2]}\ddlam)$, the operator $F(\Delta_\nu)$ is a right $\diamond_\nu$-convolution operator with convolution kernel $K_{F(\Delta_\nu)} \in L^2(G_\nu)$, and moreover the function $\modu_\nu^{-1/2} K_{F(\Delta_\nu)}$ is radial on $G_\nu$. Thus, we can write
\begin{equation}\label{eq:H_notat}
K_{F(\Delta_\nu)}(\bdx) = \frac{2^{(2-\nu)/2}}{\Gamma(\nu/2)} \modu_\nu^{1/2}(\bdx) H_{F(\Delta_\nu)}(|\bdx|_\dist)
\end{equation}
for a suitable $H_{F(\Delta_\nu)} : \Rnon \to \CC$; the choice of the normalization constant in \eqref{eq:H_notat} makes the following computations neater.

In this section we shall derive an essentially explicit formula for $H_{F(\Delta_\nu)}$, in terms of the Euclidean Fourier transform of the function $\xi \mapsto F(\xi^2)$ on $\RR$. The starting point for our discussion is the case where $F(\lambda) = e^{-t\lambda}$, corresponding to the heat semigroup for $\Delta_\nu$.

From \cite[Proposition 4.4]{MaPl} we know that the heat semigroup $\{e^{-t\Delta_\nu} \}_{t>0}$ associated with $\Delta_\nu$ is a family of $\diamond_\nu$-convolution operators, and the corresponding kernels are given by
\begin{equation}\label{eq:6}
K_{e^{-t\Delta_\nu}}(x,u)=\frac{2}{\Gamma(\nu/2)} \int_0^\infty \Psi_t(\xi) \, (2\xi e^u)^{-\nu/2} \exp\left(-\frac{\cosh u}\xi-\frac{x^2}{2\xi e^u}\right)\,\dd\xi,
\end{equation}
where
\begin{equation*}
\Psi_t(\xi)
=\frac{e^{\frac{\pi^2}{4t}}}{\xi^2 \sqrt{4\pi^3 t}} \int_0^\infty \sinh \theta \, \sin\left(\frac{\pi\theta}{2t}\right) \exp\left(-\frac{\theta^2}{4t}-\frac{\cosh \theta}{\xi}\right)\,\dd\theta,\qquad \xi>0.
\end{equation*}
The formula \eqref{eq:6} follows from the relation between the heat semigroups $e^{-t\Delta_\nu}$ and $e^{-tL_\nu}$ on $G_\nu$ and $X_\nu$ (see \cite[Theorem~2.1]{G}).
Indeed, \eqref{eq:6} can be rewritten as
\begin{equation}\label{eq:7}
K_{e^{-t\Delta_\nu}}(x,u)=\int_0^\infty \Psi_t(\xi) \exp\left(-\frac{\cosh u}\xi \right) K_{e^{-\xi e^u L_\nu/2}}(x)\,\dd\xi,
\end{equation}
where $K_{e^{-tL_\nu}}$ is the $\ast_\nu$-convolution kernel of $e^{-tL_\nu}$, given by
\begin{equation}\label{eq:13}
K_{e^{-tL_\nu}}(x) = S_\nu \, (4\pi t)^{-\nu/2} e^{-\frac{x^2}{4t}},
\end{equation}
while
\[
S_\nu \defeq \frac{2\pi^{\nu/2}}{\Gamma(\nu/2)}
\]
(see, e.g., \cite[Lemma 3.3]{MaPl}).
If $\nu$ is an integer, then $S_\nu$ is the surface measure of a unit sphere in $\RR^\nu$, and the remaining part of \eqref{eq:13} is the radial profile of the heat kernel on $\RR^\nu$, as one should expect. 

Correspondingly, the constant $S_\nu$ appears as a multiplying factor when we compare our formulas for convolution kernels on $X_\nu$ and $G_\nu$ with the classical ones available in the literature for $\RR^\nu$ and $\RR^\nu\rtimes \RR$ in the case $\nu\in\Npos$.
Actually, due to the further normalisation constant in \eqref{eq:H_notat}, when directly comparing expressions for $H_{F(\Delta_\nu)}$ and  $K_{F(\Delta_{\RR^\nu})}$, a simpler multiplying factor appears, viz.,
\begin{equation}\label{eq:tS_nu}
\tilde S_\nu \defeq \frac{\Gamma(\nu/2)}{2^{(2-\nu)/2}} S_\nu = (2\pi)^{\nu/2}.
\end{equation}

For example, for $\nu=2$, from the explicit formula for the heat kernel for $\Delta_{\RR^2}$ on $\RR^2 \rtimes \RR$ (see, e.g, \cite[p.~1123]{SjVa}) and \eqref{eq:H_notat} we deduce that
\begin{equation}\label{eq:3}
H_{e^{-t\Delta_2}}(r) = \frac{1}{\sqrt{16 \pi t^3}} \frac{r}{\sinh r} \exp\left(-\frac{r^2}{4t}\right).
\end{equation}

On the other hand, for an arbitrary $\nu \geq 1$, by \eqref{eq:6} and \eqref{eq:dist},
\begin{equation}\label{eq:K_nu}
H_{e^{-t\Delta_\nu}}(r)=\int_0^\infty \Psi_t(\xi)\,\xi^{-\nu/2} \exp\left(-\frac{\cosh r}\xi\right)\,\dd\xi.
\end{equation}
From this formula we shall now derive a number of relations between heat kernels corresponding to different values of the dimensional parameter $\nu$; in the case of integer $\nu$, analogous formulas can be found, e.g., in \cite[eqs.\ (2.2)-(2.3)]{AnOs}.

\begin{prop}\label{prop:heat_dim}
Let $\nu \geq 1$.
\begin{enumerate}[label=(\roman*)]
\item\label{en:heat_dim_der} For any $k\in\NN$,
\begin{equation}\label{eq:19}
\left(\frac{-1}{\sinh r} \partial_r\right)^k H_{e^{-t\Delta_\nu}}(r) = H_{e^{-t\Delta_{\nu+2k}}}(r).
\end{equation}
\item\label{en:heat_dim_int} If $\beta \in (0,\infty)$, then
\begin{equation*}
	H_{e^{-t\Delta_\nu}}(r) = \frac{1}{\Gamma(\beta)} \int_r^\infty \frac{\sinh x}{(\cosh x-\cosh r)^{1-\beta}} H_{e^{-t\Delta_{\nu+2\beta}}}(x)\,\dd x.
\end{equation*}
\item\label{en:heat_dim_intder} There holds
\begin{equation*}\begin{split}
H_{e^{-t\Delta_\nu}}(r) 
&=\frac{1}{\Gamma(1-\{\frac{\nu}{2}\})} \int_r^\infty \frac{\sinh x}{(\cosh x-\cosh r)^{\{\frac\nu2 \}}} \left(\frac{-1}{\sinh x}\partial_x\right)^{\lfloor \frac{\nu}{2}\rfloor} H_{e^{-t\Delta_2}}(x)\,\dd x\\
&=\frac{1}{\Gamma(1-\{\frac{\nu}{2}\})} \left(\frac{-1}{\sinh r}\partial_r\right)^{\lfloor \frac{\nu}{2}\rfloor} \int_r^\infty \frac{\sinh x}{(\cosh x-\cosh r)^{\{\frac\nu2 \}}} H_{e^{-t\Delta_2}}(x)\,\dd x.
\end{split}
\end{equation*}
\end{enumerate}
\end{prop}
\begin{proof}
\ref{en:heat_dim_der}. This follows immediately from \eqref{eq:K_nu} and differentiation under the integral sign.

\ref{en:heat_dim_int}. 	We use \eqref{eq:K_nu} and obtain
\begin{multline*}
\int_r^\infty \frac{ H_{e^{-t\Delta_{\nu+2\beta}}}(x)  \sinh x}{(\cosh x-\cosh r)^{1-\beta}} \,\dd x  \\
= \int_0^\infty \Psi_t(\xi) \, \xi^{-(\nu+2\beta)/2} \int_r^\infty \frac{e^{-\frac{\cosh x}{\xi}} \sinh x }{(\cosh x-\cosh r)^{1-\beta}}  \,\dd x \,\dd\xi.
\end{multline*}
Notice that $|\Psi_t(\xi)| \lesssim_t \xi^{-2}$, see \cite[Theorem 2.1]{G}, so all the integrals are absolutely convergent.
In the inner integral we substitute $y=(\cosh x-\cosh r)/\xi $ and get
\[\begin{split}
\int_r^\infty \frac{ H_{e^{-t\Delta_{\nu+2\beta}}}(x) \sinh x}{(\cosh x-\cosh r)^{1-\beta}} \,\dd x  
&= \int_0^\infty \Psi_t(\xi) \, \xi^{-(\nu+2\beta)/2} \int_0^\infty (y\xi)^\beta e^{-y-\frac{\cosh r}{\xi}}   \frac{\dd y}{y} \,\dd\xi\\
&= \int_0^\infty \Psi_t(\xi) \, \xi^{-\nu/2} e^{-\frac{\cosh r}{\xi}}\int_0^\infty y^\beta e^{-y}   \frac{\dd y}{y} \,\dd\xi\\
&=\Gamma(\beta) H_{e^{-\Delta_{\nu}}}(r),
\end{split}\]
as desired.

\ref{en:heat_dim_intder}.	By part \ref{en:heat_dim_int} with $\beta=1-\{\frac{\nu}{2} \}$ and formula \eqref{eq:19} we get
\[\begin{split}
H_{e^{-t\Delta_\nu}}(r) 
&= \frac{1}{\Gamma(1-\{\frac\nu2 \})} \int_r^\infty \frac{\sinh x}{(\cosh x-\cosh r)^{\{\frac\nu2 \}}} H_{e^{-t\Delta_{2\lfloor\frac{\nu}{2}\rfloor+2}}}(x)\,\dd x\\
&=\frac{1}{\Gamma(1-\{\frac{\nu}{2}\})} \int_r^\infty \frac{\sinh x}{(\cosh x-\cosh r)^{\{\frac\nu2 \}}} \left(\frac{-1}{\sinh x}\partial_x\right)^{\lfloor \frac{\nu}{2}\rfloor}  H_{e^{-t\Delta_2}}(x)\,\dd x.
\end{split}\]

In order to conclude, it is enough to repeatedly apply the formula
\begin{multline}\label{eq:intbyparts}
	\left(\frac{1}{\sinh r} \partial_r\right) \int_r^\infty \frac{\sinh x}{(\cosh x- \cosh r)^\alpha} g(x)\,\dd x\\
	= \int_r^\infty \frac{\sinh x}{(\cosh x- \cosh r)^\alpha} \left(\frac{1}{\sinh x} \partial_x\right)g(x)\,\dd x,
\end{multline}
which holds for any $\alpha \in (0,1)$ and any $g : \Rnon \to \CC$ which decays sufficiently fast together with its derivative. To justify \eqref{eq:intbyparts},
we use the substitution $y=\cosh x-\cosh r$ in both integrals, and notice that $\frac{1}{\sinh x}\partial_x$ changes to $\partial_y$. Thus \eqref{eq:intbyparts} takes the form
\begin{equation*}
	\left(\frac{1}{\sinh r} \partial_r\right) \int_0^\infty y^{-\alpha} h(y+\cosh r) \,\dd y = \int_0^\infty y^{-\alpha} h'(y+\cosh r)\,\dd y,
\end{equation*} 
where $h(s)=g(\arccosh(s))$, and the latter formula follows by taking the derivative under the integral sign.
\end{proof}

Notice that if we let $\Delta_0$ be the classical Laplacian on $\RR$ and
\begin{equation*}
H_{e^{-t\Delta_0}}(r)\defeq \frac{1}{\sqrt{4\pi t}} e^{-\frac{r^2}{4t}}
\end{equation*} 
be the heat kernel on $\RR$, then from \eqref{eq:3} we see that
\begin{equation}\label{eq:heat_20}
H_{e^{-t\Delta_2}}(r) = \frac{-1}{\sinh r} \partial_r H_{e^{-t\Delta_0}}(r);
\end{equation}
in other words, with this notation, the relation \eqref{eq:19} also holds when $\nu = 0$.
From Proposition \ref{prop:heat_dim} we shall now derive a relation between $H_{F(\Delta_\nu)}$, $\nu\geq 1$, and $H_{F(\Delta_0)}$, where
\begin{equation}\label{eq:invfourier}
H_{F(\Delta_0)}(x) = \frac{1}{2\pi} \int_{\RR} F(\xi^2) \, e^{i\xi x}\,\dd\xi
\end{equation}
is the convolution kernel of $F(\Delta_0)$ on $\RR$.

\begin{cor}\label{cor:5}
Let $\nu\geq 1$ and let $F\in \Sz(\Rnon)$. Then 
\begin{multline*}
H_{F(\Delta_\nu)}(r) \\
= \frac{1}{\Gamma(1-\{\frac{\nu}{2}\})} \int_r^\infty \frac{\sinh x}{(\cosh x-\cosh r)^{\{\frac\nu2 \}}} \left(\frac{-1}{\sinh x}\partial_x\right)^{\lfloor\frac{\nu}{2}\rfloor +1} \, H_{F(\Delta_{0})}(x) \,\dd x
\end{multline*}
for all $r \in \Rnon$.
\end{cor}
\begin{proof}
With the change of variables $G(\lambda)=F(\lambda^2)$, we can reformulate the statement as
\begin{equation}\label{eq:rel_kernels_dim}
H_{G(\sqrt{\Delta_\nu})}(r) 
= \int_r^\infty \frac{x}{(\cosh x-\cosh r)^{\{\frac\nu2 \}}} E_\nu H_{G(\sqrt{\Delta_{0}})}(x)\,\dd x
\end{equation}
for all $G \in \Sz_\even(\Rnon)$, where $E_\nu$ is the differential operator given by
\[
E_\nu = \frac{-1}{\Gamma(1-\{\frac{\nu}{2}\})} x^{-1}\partial_x \left(\frac{-1}{\sinh x}\partial_x\right)^{\lfloor\frac{\nu}{2}\rfloor}.
\]

Notice now that, for any fixed $r \in \Rnon$, both sides of \eqref{eq:rel_kernels_dim} are continuous linear functions of $G \in \Sz_\even(\Rnon)$ with respect to the Schwartz class topology. For the left-hand side, this is clear from Lemma \ref{lm:1}. As for the right-hand side, this follows from the fact that $G \mapsto E_\nu H_{G(\sqrt{\Delta_0})}$ is a continuous operator on $\Sz_\even(\Rnon)$: indeed it is the composition of the inverse Fourier transform on $\RR$ (see \eqref{eq:invfourier}) and the differential operator $E_\nu$, which are clearly bounded on $\Sz_\even(\Rnon)$ (for the boundedness of the latter, see also \cite{Stk3}).

In light of these continuity properties, it is enough to check the identity \eqref{eq:rel_kernels_dim} on a class of functions whose linear span is dense in $\Sz_\even(\Rnon)$. One such class is that of gaussians, i.e., functions $G(\lambda) = e^{-t\lambda^2}$ for $t>0$ (see, e.g., \cite[Lemma~2]{J-McD}). On the other hand, in the case of gaussians, \eqref{eq:rel_kernels_dim} follows immediately by combining Proposition \ref{prop:heat_dim}\ref{en:heat_dim_intder} and \eqref{eq:heat_20}.
\end{proof}

\begin{rem}
The relation between $H_{F(\Delta_\nu)}$ and $H_{F(\Delta_{0})}$ in Corollary \ref{cor:5} could be interpreted in term of a suitable (inverse) Abel transform, in the sense of \cite[Section 5.3]{Koo}. A discrete analogue in the setting of flow trees is also stated in \cite[Proposition 5.4]{MSTV}.
\end{rem}

\subsection{Dunkl setting}\label{ss:dunkl}

We briefly recall some facts from the Dunkl theory for rank-one root systems. For a comprehensive treatment of the Dunkl theory we refer the reader to \cite{Du,dJ,Ro'99,Ro'03}. As we are interested in the Dunkl Laplacian and Dunkl transform on $\RR$ as ``extensions'' of the Bessel operator and Hankel transform on $X_\nu$, we shall use the same dimension parameter $\nu \geq 1$ as before to index them, instead of the multiplicity parameter $k = (\nu-1)/2$ or the Bessel function order $\nu/2-1$ used in other sources in the literature.

We shall denote by $X^\rD_\nu$ the space $\RR$ equipped with the measure $\dmu^{\rD}(x)=2^{-1} |x|^{\nu-1}\,\dd x$; notice that integrating an even function on $\RR$ with respect to $\dmu^{\rD}$ is the same as integrating it on $\Rnon$ with respect to $\dmu$. In this section, in order not to burden the notation, we shall often identify even functions on $X^\rD_\nu$ with their restrictions to $X_\nu$. In this way, operators such as $L_\nu$, $\Hank_\nu$ or $\tau_\nu^{[x]}$ discussed earlier can be also thought of as acting on even functions on $X^\rD_\nu$.

The \emph{Dunkl transform} of $f\in L^1(X^\rD_\nu)$ is given by
\begin{equation*}
\Hank_\nu^\rD f(\xi) = \int_\RR  f(x)\, \overline{j^{\nu,\rD}_{\xi}(x)} \, \dmu^{\rD}(x), \qquad\xi \in \RR
\end{equation*}
(see \cite[Example~2.29]{RoVo} and \cite[Remark 3.7.2]{Ro'95}), where 
\[
j^{\nu,\rD}_x(y) = j^{\nu,\rD}(xy), \qquad j^{\nu,\rD}(x)= j^\nu(x)+i\frac{x}{\nu} j^{\nu+2}(x)
\]
and $j^\nu$ is as in \eqref{eq:jnu}.
Notice that
\[
j^{\nu,\rD}(-x) = \overline{j^{\nu,\rD}(x)}, \qquad |j^{\nu,\rD}(x)| \leq 1 \qquad \forall x \in \RR
\]
\cite[Corollary 2.2]{Ro'95}. In the case $\nu=1$, we simply have $j^{1,\rD}(x) = e^{ix}$, so the Dunkl transform reduces to the Fourier transform on $\RR$, up to a normalisation factor.
The function
\[
E_\nu(x,y) = j^{\nu,\rD}(-ixy), \qquad x,y \in \CC,
\]
is also known as the \emph{Dunkl kernel} and satisfies the identity
\begin{equation}\label{eq:Dkerndiff}
D_\nu E_\nu(\cdot,y) = y E_\nu(\cdot,y)
\end{equation}
\cite[Theorem 2.27]{RoVo}, where $D_\nu$ is the \emph{Dunkl operator}
\begin{equation}\label{eq:Dop}
D_\nu f(x) = f'(x) +\frac{\nu-1}{2} \frac{f(x)-f(-x)}{x}. 
\end{equation}

For a function $f$ on $X_\nu^\rD$, we shall write $f=f_\even+f_\odd$, where $f_\even$ and $f_\odd$ are the even and odd parts of $f$, respectively. By comparing the definitions of $\Hank_\nu$ and $\Hank_\nu^\rD$, it follows that
\begin{equation}\label{eq:58}
\Hank_\nu^\rD f = \Hank_\nu f_\even -\frac{i}{\nu} U^{-1} \Hank_{\nu+2} U f_\odd,
\end{equation}
where
\begin{equation}\label{eq:opU}
Uf(x)=f(x)/x.
\end{equation}
Thus, we can think of the Dunkl transform as an extension of the Hankel transform, and deduce a number of properties of the Dunkl transform from those of the Hankel transform. For example, $\Hank_\nu^\rD$ is an isomorphism from the Schwartz class $\Sz(\RR)$ onto itself, and also from $L^2(X_\nu^\rD)$ onto itself; moreover, from the Plancherel formula \eqref{eq:plancherel_Hnu} for the Hankel transform one deduces an analogous result for the Dunkl transform:
\[
\| \Hank_\nu^\rD f \|_{L^2(X^\rD_\nu)} = 2^{\nu/2-1} \Gamma(\nu/2) \, \| f \|_{L^2(X^\rD_\nu)}.
\]

The \emph{Dunkl convolution} is defined by
\begin{equation}\label{eq:Dconv}
f\ast_\nu^\rD g (x) = \int_\RR \tau_{\nu,\rD}^{[x]} f(-y) \, g(y)\,\dmu^\rD(y) = \int_\RR  f(y) \, \tau_{\nu,\rD}^{[x]}g(-y)\,\dmu^\rD(y),
\end{equation}
where the Dunkl translations $\tau_{\nu,\rD}^{[x]}$ are given by (cf.\ \cite[Lemma 3.2]{Ro'95})
\begin{equation}\label{eq:Dtau}
\begin{split}
\tau_{\nu,\rD}^{[x]} f(y) 
&= \int_{-1}^1 \left( f_\even(\langle |x|,|y| \rangle_\omega) + \frac{x+y}{\langle |x|,|y|\rangle_\omega} f_\odd(\langle |x|,|y|\rangle_\omega) \right) \, (1-\omega \sign(xy) ) \,\dpr(\omega) \\
&= \int_{-1}^1 \left( f_\even(\langle x,y\rangle_\omega) + \frac{x+y}{\langle x,y\rangle_\omega} f_\odd(\langle x,y\rangle_\omega) \right) \, (1-\omega) \,\dpr(\omega),
\end{split}
\end{equation}
and $\dpr$ and $\langle x,y\rangle_\omega$ are as in \eqref{eq:hankel_transl}.
Notice that
\[
\tau_{\nu,\rD}^{[x]} f(y) = \tau_{\nu,\rD}^{[y]} f(x), \qquad (\tau_{\nu,\rD}^{[x]} f)\check{} = \tau_{\nu,\rD}^{[-x]} \check{f},
\]
where $\check f(y) \defeq f(-y)$; in particular, the equality of the two integrals in \eqref{eq:Dconv} corresponds to the fact that the adjoint of $\tau_{\nu,\rD}^{[x]}$ in $L^2(\RR,\dmu^\rD)$ is $\tau_{\nu,\rD}^{[-x]}$ (cf.\ \cite[\S 1.2]{Ro'95}).
For all $x,y \in \RR$ and $\omega \in [-1,1]$, there holds
\begin{equation}\label{eq:est_Dtrk}
0 \leq 1- \omega \leq 2, \quad | (x+y)(1-\omega) | \leq 2 \langle x,y\rangle_\omega
\end{equation}
(see \cite[pp.~298--299]{Ro'95}); so, by comparing \eqref{eq:Dtau} and \eqref{eq:hankel_transl},
\begin{equation}\label{eq:66}
|\tau_{\nu,\rD}^{[x]}f| \leq 2 \tau_\nu^{[|x|]}(|f_\even|+|f_\odd|),
\end{equation}
and therefore from \eqref{eq:Dconv} and \eqref{eq:Hconv} it follows that
\[
|f*_\nu^\rD g| \leq 2 (|f_\even|+|f_\odd|) *_\nu (|g_\even|+|g_\odd|).
\]
As a matter of fact, from \eqref{eq:Dconv} and \eqref{eq:Dtau} it also follows that
\begin{equation}\label{eq:DHconv_even}
f *_\nu^\rD g = f *_\nu g 
\end{equation}
whenever $f$ and $g$ are even; in this sense, the Dunkl convolution $*_\nu^\rD$ extends the Hankel convolution $*_\nu$. 

We can relate Dunkl convolution and Dunkl transform by observing that, by \cite[Theorem 2.4]{Ro'95},
\[
\tau_{\nu,\rD}^{[y]} j^{\nu,\rD}_x(z) = j^{\nu,\rD}_x(y) \, j^{\nu,\rD}_x(z),
\]
whence
\begin{equation}\label{eq:Dtr_conv}
\Hank_\nu^\rD (\tau_{\nu,\rD}^{[x]} f) = j^{\nu,\rD}_x f, \qquad
\Hank_\nu^\rD (f\ast_\nu^{\rD} g) = \Hank_\nu^\rD(f)\Hank_\nu^\rD(g)
\end{equation}
(cf.\ \cite[eq.\ (2.19)]{RoVo}).

The Dunkl operator $D_\nu$ of \eqref{eq:Dop} is formally skewadjoint with respect to $\dmu^{\rD}$, i.e., $D_\nu^+ = -D_\nu$. Moreover, from \eqref{eq:Dkerndiff} it follows that
\[
D_\nu j^{\nu,\rD}_x = ix j^{\nu,\rD}_x,
\]
whence
\begin{equation}\label{eq:Dtr_op}
\Hank_\nu^\rD (D_\nu f)(\xi) = i\xi \Hank_\nu^\rD f(\xi), \qquad \xi \in \RR
\end{equation}
(cf.\ \cite[Lemma 2.37]{RoVo}); from \eqref{eq:Dtr_conv} we then see that
\begin{equation}\label{eq:61}
D_\nu \tau_{\nu,\rD}^{[x]} f = \tau_{\nu,\rD}^{[x]} D_\nu f, \qquad 
D_\nu (f\ast_\nu^{\rD} g) = (D_\nu f)\ast_\nu^{\rD} g = f\ast_\nu^{\rD} (D_\nu g).
\end{equation}

Let $G_\nu^\rD$ be the space $\RR \times \RR$ equipped with the measure $\dmub^{\rD}(x,u) = \dmu^\rD(x) \, \du$.
For a function $f$ on $G_\nu$ we denote by $f_\even$ and $f_\odd$ the parts of $f$ which are even and odd with respect to the first variable. We shall always use this meaning of ``even'' and ``odd'' on $G_\nu^\rD$, and, much as before, in this section we shall identify even functions on $G_\nu^\rD$ with functions on $G_\nu$.

In analogy to \eqref{eq:39} we define
\[
f\diamond_\nu^{\rD} g(x,u) = \int_\RR (f^{v}\ast_\nu^{\rD} g_{(e^{v})}^{u-v})(x)\,\dd v,\qquad x,u\in\RR.
\]
Hence, much as in \eqref{eq:HaId_conv}, by \eqref{eq:Dtr_conv} and the fact that $\Hank_\nu^\rD(g_{(\lambda)})(\xi) = \Hank_\nu^\rD g(\lambda \xi)$ we deduce that
\begin{equation}\label{eq:Dtr_diamond}
\HaId^\rD (f\diamond_\nu^{\rD} g)(x,u) = \int_\RR \HaId^\rD f (\xi,v) \, \HaId^\rD g(e^v\xi,u-v)\,\dd v,
\end{equation}
where we used the notation $\HaId^\rD \defeq \Hank_\nu^\rD\otimes \id$ for the ``partial Dunkl transform''. Notice also that, by \eqref{eq:DHconv_even},
\begin{equation}\label{eq:DHconv_even_lift}
f \diamond^\rD_\nu g = f \diamond_\nu g
\end{equation}
whenever $f$ and $g$ are even.

Let us also define left and right translations on $G_\nu^\rD$:
\begin{equation}\label{eq:D-transl}
\ell_{(y,v)}^\rD f(x,u) = r_{(x,u)}^\rD f(y,v) =\tau_{\nu,\rD}^{[e^{v}x]}f(y,u+v).
\end{equation}
Thus, the $\diamond_\nu^\rD$-convolution takes a form that is similar to \eqref{eq:diamond}:
\begin{equation}\label{eq:Ddiamond2}
f\diamond_\nu^\rD g(\bdx)=\int_{G_\nu^\rD}   f(\bdy^{-1})r^\rD_{\bdx}g(\bdy)\,\dmub^\rD(\bdy) =\int_{G_\nu^\rD}   \ell^\rD_{\bdx}f(\bdy^{-1}) g(\bdy) \,\dmub^\rD(\bdy),
\end{equation}
where $(x,u)^{-1} \defeq (-e^{-u}x,-u)$. Notice that \eqref{eq:66} immediately gives
\begin{equation}\label{eq:72}
|r^{\rD}_{(x,u)} f|\leq 2r_{(|x|,u)} (|f_\even|+|f_\odd|),\qquad |\ell^{\rD}_{(x,u)} f|\leq 2\ell_{(|x|,u)} (|f_\even|+|f_\odd|),
\end{equation}
whence also
\[
|f\diamond_\nu^\rD g| \leq 2 (|f_\even| + |f_\odd|) \diamond_\nu (|g_\even|+|g_\odd|).
\]
From \eqref{eq:young} we then deduce Young-type inequalities for the $\diamond_\nu^\rD$-convolution: for all $p,q,r \in [1,\infty]$ with $1/p+1/q=1+1/r$,
\[
\Vert (f \modu_\nu^{-1/q'}) \diamond_\nu^\rD g\Vert_{L^r(G_\nu^\rD)} \leq 8 \Vert f \Vert_{L^p(G_\nu^\rD)} \Vert g \Vert_{L^q(G_\nu^\rD)}
\]
where $\modu_\nu(x,u) \defeq e^{-\nu u}$ for all $(x,u) \in G_\nu^\rD$.

We shall also use the $L^1$-isometric involution $f \mapsto f^\bullet$ given by
\[
f^{\bullet}(\bdx)\defeq \modu_\nu(\bdx)\overline{f(\bdx^{-1})}.
\]
The relation between involution and convolution on $G_\nu$ is given in the following statement, which is an analogue of Lemma \ref{lm:conv_inv}.

\begin{lm}\label{lm:conv_inv_D}
The following hold.
\begin{enumerate}[label=(\roman*)]
\item\label{en:conv_inv_D_magic} For any $f \in L^1_\loc(G_\nu^\rD)$ there holds
\begin{equation*}
\modu_\nu(\bdy) r^\rD_{\bdx} f(\bdy^{-1})= \modu_\nu(\bdx) \overline{r^\rD_{\bdy}f^{\bullet}(\bdx^{-1})}
\end{equation*}
for almost all $\bdx,\bdy \in G_\nu$.
\item\label{en:conv_inv_D_adj} For all $g \in L^1_\loc(G_\nu^\rD)$ and all compactly supported $f,h \in L^\infty(G_\nu^\rD)$,
\begin{equation}\label{eq:inv_conv_D}
(f \diamond_\nu^\rD g)^\bullet = g^\bullet \diamond_\nu^\rD f^\bullet, \qquad \langle f \diamond_\nu^\rD g , h \rangle_{L^2(G_\nu^\rD)} = \langle f , h \diamond_\nu^\rD g^\bullet \rangle_{L^2(G_\nu^\rD)}.
\end{equation}
\end{enumerate}
\end{lm}
\begin{proof}
\ref{en:conv_inv_D_magic}.
Let $\bdx=(x,u),\bdy=(y,v)\in G_\nu^\rD$. By \eqref{eq:D-transl} and \eqref{eq:Dtau},
\[\begin{split}
		&\modu_\nu(\bdx) \overline{r^\rD_{\bdy}f^\bullet(\bdx^{-1})} \\
		&=e^{-\nu u} \overline{\tau_{\nu,\rD}^{[e^{-u}y]} f^\bullet(-e^{-u}x,v-u)} \\
		&= e^{-\nu u} \int_{-1}^1 \left[ \overline{(f^\bullet)_\even (e^{-u}\langle -x,y\rangle_\omega,v-u)} + \frac{y-x}{\langle -x,y \rangle_\omega} \overline{(f^\bullet)_\odd (e^{-u}\langle -x,y\rangle_\omega,v-u)} \right] \\
		&\qquad\times (1-\omega)\,\dpr(\omega)\\
		&= e^{-\nu v} \int_{-1}^1 \left[ f_\even(e^{-v}\langle x,-y\rangle_\omega,u-v) + \frac{x-y}{\langle x,-y \rangle_\omega} f_\odd (e^{-v}\langle x,-y\rangle_\omega,u-v) \right] \\
		&\qquad\times (1-\omega)\,\dpr(\omega)\\
		&= e^{-\nu v} \tau_{\nu,\rD}^{[e^{-v}x]} f(-e^{-v}y,u-v)\\
		&=\modu_\nu(\bdy) r^\rD_{\bdx}f(\bdy^{-1}).
\end{split}\]

\ref{en:conv_inv_D_adj}. The proof is analogous to that of Lemma \ref{lm:conv_inv}\ref{en:conv_inv} and is omitted.
\end{proof}

The next lemma is the Dunkl counterpart of Lemma \ref{lm:11}.
Notice that the vector fields $\vfX_0 = \partial_u$ and $\vfX_1 = e^u \partial_x$ can be thought of as defined on the whole $G_\nu^\rD$, so it also makes sense to apply them and the gradient $\nabla_\nu = (\vfX_0,\vfX_1)$ to functions on $G_\nu^\rD$. We also introduce the multiplication operator $\tilde U$ on $G_\nu$ given by
\begin{equation}\label{eq:tildeX1}
\tilde U = e^u U = e^u x^{-1},
\end{equation}
which can be thought of as a lifting to $G_\nu$ of the operator $U$ on $X_\nu$, see \eqref{eq:opU}.
We further introduce the notation $\tilde\bdx \defeq (|x|,u)$ for any $\bdx = (x,u) \in G^\rD_\nu$.

\begin{lm}\label{lm:12}
Let $x\in \RR$ and $\bdx\in G^\rD_\nu$.
\begin{enumerate}[label=(\roman*)]	
\item\label{lm:12(1)} Let $g \in C(G_\nu^\rD)$. Then
\[
 \supp (r_\bdx^\rD g) \subseteq \{ \bdz \in G_\nu^\rD \tc \dist(\supp(|g_\even|+|g_\odd|),\tilde\bdz) \leq |\tilde\bdx|_\dist\}.
\]
\item\label{lm:12(2)}
Let $f\in C^1(X^\rD_\nu)$ and $g\in C^1(G^\rD_\nu)$. Then the function $(x,y) \mapsto \tau_{\nu,\rD}^{[x]} f(y)$ is in $C^1(X^\rD_\nu \times X^\rD_\nu)$, the function $(\bdx,\bdy) \mapsto r^\rD_{\bdx} g(\bdy)$ is in $C^1(G^\rD_\nu \times G^\rD_\nu)$ and
\begin{equation}\label{eq:73}
\begin{aligned}
|\partial_x \tau_{\nu,\rD}^{[x]} f|  &\leq 2\tau_\nu^{[|x|]} (|(f_\even)'|+|(f_\odd)'|+2|U f_\odd|),\\
| \nabla_{\nu}^{\bdx} r^\rD_{\bdx} g| &\leq 2r_{\tilde\bdx}(|\nabla_\nu g_\even|+|\nabla_\nu g_\odd|+2| \tilde{U} g_\odd|)
\end{aligned}
\end{equation}
pointwise.
\end{enumerate}
\end{lm}
\begin{proof}
\ref{lm:12(1)}. This follows immediately from Lemma \ref{lm:11}\ref{lm:11(1)} and \eqref{eq:72}.

\ref{lm:12(2)}.
Let us first discuss the result for $f \in C^1(X^\rD_\nu)$. If $\nu = 1$, then \eqref{eq:Dtau} simply reduces to
\[
\tau_{1,\rD}^{[x]} f(y) = f(x+y),
\]
whence the stated $C^1$ regularity and estimate follow.

Suppose instead that $\nu > 1$, and notice that $\langle x,y \rangle_\omega > 0$ for all $(x,y) \in \RR^2 \setminus \{(0,0)\}$ and $\omega \in (-1,1)$. Thus, from the definition of the Dunkl translations \eqref{eq:Dtau} and the bounds \eqref{eq:bd_xycos} and \eqref{eq:est_Dtrk}, by differentiating under the integral sign we obtain that, for all $(x,y) \neq (0,0)$,
\[\begin{split}
\partial_x \tau_{\nu,\rD}^{[x]} f(y) 
&= \int_{-1}^1 (f_\even)'(\langle x,y\rangle_\omega) \frac{x-y\omega}{\langle x,y\rangle_\omega} (1-\omega) \,\dpr(\omega) \\
&\quad+ \int_{-1}^1 \left( (f_\odd)'(\langle x,y\rangle_\omega) - \frac{f_\odd(\langle x,y \rangle_\omega)}{\langle x,y\rangle_\omega} \right)  \frac{x-y\omega}{\langle x,y\rangle_\omega} \frac{x+y}{\langle x,y \rangle_\omega} (1-\omega) \,\dpr(\omega) \\
&\quad+ \int_{-1}^1 \frac{f_\odd(\langle x,y \rangle_\omega)}{\langle x,y\rangle_\omega} (1-\omega) \,\dpr(\omega),
\end{split}\]
as well as an analogous formula (with $x$ and $y$ switched) for $\partial_y \tau_{\nu,\rD}^{[x]} f(y)$. This expression, together with the fact that $f_\even,f_\odd \in C^1(\RR)$, prove the $C^1$ regularity of $(x,y) \mapsto \tau_{\nu,\rD}^{[x]} f(y)$ on $\RR^2 \setminus \{(0,0)\}$. Moreover, by using the fact that $(f_\even)'(0) = 0$ and $\lim_{h \to 0} f_\odd(h)/h = (f_\odd)'(0)$, together with the bounds \eqref{eq:bd_xycos} and \eqref{eq:est_Dtrk}, from the above formulas one sees that 
\[
\lim_{(x,y) \to (0,0)} \partial_x \tau_{\nu,\rD}^{[x]} f(y) = \lim_{(x,y) \to (0,0)} \partial_y \tau_{\nu,\rD}^{[x]} f(y) = (f_\odd)'(0) = f'(0),
\]
thus completing the proof of the $C^1$ regularity of $(x,y) \mapsto \tau_{\nu,\rD}^{[x]} f(y)$ on $\RR^2$. Finally, by the bounds \eqref{eq:bd_xycos} and \eqref{eq:est_Dtrk}, together with \eqref{eq:hankel_transl} and \eqref{eq:opU}, we also see that
\begin{align*}
\left| \int_{-1}^1 (f_\even)'(\langle x,y\rangle_\omega) \frac{x-y\omega}{\langle x,y\rangle_\omega} (1-\omega) \,\dpr(\omega) \right| &\leq 2\tau_\nu^{[|x|]}|(f_\even)'|(y),\\
\left| \int_{-1}^1 (f_\odd)'(\langle x,y\rangle_\omega)  \frac{x-y\omega}{\langle x,y\rangle_\omega} \frac{x+y}{\langle x,y \rangle_\omega} (1-\omega) \,\dpr(\omega) \right| &\leq 2\tau_\nu^{[|x|]}|(f_\odd)'|(y), \\
\left| \int_{-1}^1 \frac{f_\odd(\langle x,y \rangle_\omega)}{\langle x,y\rangle_\omega}   \frac{x-y\omega}{\langle x,y\rangle_\omega} \frac{x+y}{\langle x,y \rangle_\omega} (1-\omega) \,\dpr(\omega) \right| &\leq 2\tau_\nu^{[|x|]} |U f_\odd|(y), \\
\left|\int_{-1}^1 \frac{f_\odd(\langle x,y \rangle_\omega)}{\langle x,y\rangle_\omega} (1-\omega) \,\dpr(\omega)\right| &\leq 2\tau_\nu^{[|x|]} |U f_\odd|(y),
\end{align*}
whence the first estimate in \eqref{eq:73} follows.
	
We now move to the results for $g \in C^1(G_\nu^\rD)$. By the definition \eqref{eq:D-transl} of the right translations on $G_\nu^\rD$, the $C^1$ regularity of $(\bdx,\bdy) \mapsto r^\rD_{\bdx} g(\bdy)$ follows from the analogous one on $X_\nu^\rD$ that we have just proved. It remains to justify the second inequality in \eqref{eq:73}. Notice that, by the definition \eqref{eq:D-transl} of the right translations on $G_\nu^\rD$, arguing as in \eqref{eq:dutransl} and using \eqref{eq:72},
\begin{equation*}
	| \partial_u r^\rD_{(x,u)}g(y,v)| 
	= |r^\rD_{(x,u)} \vfX_0 g(y,v)| 
	\leq 2 r_{(|x|,u)} [|\vfX_0 g_\even| + |\vfX_0 g_\odd|](y,v);
\end{equation*}
similarly, using also the first bound in \eqref{eq:73},
\[\begin{split}
| e^u \partial_x r^\rD_{(x,u)}g(y,v)| 
&= | e^{u+v} (\partial_z \tau_{\nu,\rD}^{[z]}g)(y,u+v)|_{z=e^v x} | \\
&\leq 2 e^{u+v} \tau_\nu^{[e^v |x|]} [|\partial_y (g_\even)| + |\partial_y (g_\odd)| + 2|U g_\odd|](y,u+v) \\
&= 2 r_{(|x|,u)} [| \vfX_1 g_\even| + | \vfX_1 g_\odd| + 2| \tilde U g_\odd| ](y,v).
\end{split}\]
By arguing much as in \eqref{eq:18}, we can combine these two bounds and obtain the second inequality in \eqref{eq:73}.
\end{proof}

The next result is an analogue of Proposition \ref{prop:SIO_diamond}. To state it, we introduce the notation
\[
\vfX_0^\rD = \partial_u, \qquad \vfX_1^\rD = e^u D_\nu.
\]
Notice that $\vfX_1^\rD$ can be thought of as a lifting of the Dunkl operator $D_\nu$ to $G_\nu^\rD$.

\begin{prop}\label{prop:SIO_diamond_D}
Let $K \in L^1_\loc(G_\nu^\rD) \cap C^1(G_\nu^\rD \setminus \{\bdzero\})$. Then, for all $f \in C_c(G_\nu^\rD)$, the function $f \diamond_\nu^\rD K$ is continuosly differentiable on $G_\nu^\rD \setminus \supp (|f_\even|+|f_\odd|)$ and, for $j=0,1$, 
\[
\vfX_j^\rD(f \diamond_\nu^\rD K)(\bdx) = \int_{G_\nu^\rD} \ell_{\bdx}^\rD f(\bdy^-) \, \vfX_j^\rD K(\bdy) \,\dmub^\rD(\bdy) \quad\text{for all } \bdx \notin \supp(|f_\even|+|f_\odd|).
\]
Moreover, for all $f,g \in C_c(G_\nu^\rD)$ with $\supp(|f_\even|+|f_\odd|) \cap \supp(|g_\even|+|g_\odd|) = \emptyset$,
\begin{equation}\label{eq:SIO_der_adj_D}
\begin{split}
\langle \vfX_j^\rD(f \diamond_\nu K), g \rangle_{L^2(G_\nu^\rD)} 
&= \int_{G_\nu^\rD} \overline{g(\bdx)} \int_{G_\nu} \ell_{\bdx} f(\bdy^{-}) \, \vfX_j^\rD K(\bdy) \,\dmub^\rD(\bdy) \,\dmub^\rD(\bdx) \\
&= \int_{G_\nu^\rD} f(\bdx) \overline{\int_{G_\nu} \ell_{\bdx} g(\bdy^{-}) \, (\vfX_j^\rD K)^\bullet(\bdy) \,\dmub^\rD(\bdy)} \,\dmub^\rD(\bdx).
\end{split}
\end{equation}
\end{prop}
\begin{proof}
From \eqref{eq:61} and \eqref{eq:D-transl} one sees that the operators $\vfX_j^\rD$ commute with left translations on $G_\nu^\rD$, i.e.,
\begin{equation}\label{eq:XjD_li}
\vfX_j^\rD \ell^\rD_\bdx g = \ell^\rD_\bdx \vfX_j^\rD g \qquad \forall \bdx \in G_\nu^\rD, \ g \in C^1(G_\nu^\rD).
\end{equation}
Thus, one can essentially repeat the proof of Proposition \ref{prop:SIO_diamond}, using Lemmas \ref{lm:conv_inv_D} and \ref{lm:12} in place of Lemmas \ref{lm:conv_inv} and \ref{lm:11}.

Notice that the operators $D_\nu$ and $\vfX_1^\rD$ are not differential operators and are not local; however, they are local when restricted to even functions and to odd functions. Therefore, one can replace the kernel $K$ with its regularisation $K_\varepsilon$ as in the proof of Proposition \ref{prop:SIO_diamond} by making sure that both even and odd parts of $K_\varepsilon$ coincide with those of $K$ off a suitable neighbourhood of $\bdzero$.
\end{proof}

\section{Calder\'on--Zygmund theory}\label{s:CZ}

We recall the definition of Calder\'on--Zygmund spaces (see \cite{HeSt,Va-PhD} and compare \cite[Definition~3.1]{MOV}).

\begin{defin}\label{def:1}
A metric-measure space $(X,d,\mu)$ is a \emph{Calder\'on--Zygmund space} if there exists a constant $\kappa\in[1,\infty)$ and a family $\admR$ of Borel subsets of $X$ satisfying the following properties.
\begin{enumerate}[label=(\roman*)]
\item For all $R \in \admR$ there exist $x(R)\in X$ and $r(R)>0$ such that
\[
R\subseteq B_X(x(R),\kappa r(R)).
\]
where $B_X(x,r)$ denotes the open ball of centre $x$ and radius $r$ in $X$.
\item $\mu(R^\ast)\leq \kappa \mu(R)$ for all $R\in \admR$, where
\[
	R^\ast \defeq \{x\in X \tc d(x,R)<r(R) \}.
\]
\end{enumerate}
Moreover, for all $f\in L^1(X)$ and $\lambda>\kappa \Vert f\Vert_{L^1(X)}/\mu(X)$ (if $\mu(X)=\infty$, then $\lambda>0$) there exists a decomposition
\begin{equation}\label{eq:CZdec}
f=g+\sum_{i\in\NN} b_i
\end{equation}
satisfying
\begin{enumerate}[label=(\roman*),resume]
\item\label{it:3,def:1} $\Vert g\Vert_{L^\infty(X)}\leq \lambda$;
\item\label{it:4,def:1} $\int_X b_i\,\dd\mu=0$ for all $i\in\NN$, and $\sum_{i\in\NN} \Vert b_i\Vert_{L^1(X)}\leq 2\kappa \Vert f\Vert_{L^1(X)}$;
\item\label{it:5,def:1} there exists a family $\{R_i\}_{i\in\NN}\subseteq \admR$ such that $\supp b_i\subseteq R_i$, and $\sum \mu(R_i)\leq \kappa \Vert f\Vert_{L^1(X)}/\lambda$;
\item\label{it:6,def:1} if $f \in L^p(X)$ for some $p \in [1,\infty)$, then $g,b_i \in L^p(X)$ too and the sum in \eqref{eq:CZdec} converges in $L^p(X)$.
\end{enumerate}
\end{defin}

For convenience, we shall work with a slightly less general notion of Calder\'on--Zygmund spaces, where we replace the conditions \ref{it:3,def:1}--\ref{it:6,def:1} referring to the decomposition of functions by assumptions concerning only the family $\admR$.

\begin{defin}
Let $(X,\mu)$ be a measure space and $E \subseteq X$ be measurable. A \emph{quasi-partition} of $E$ is an at most countable family $\cP$ of pairwise disjoint measurable subsets of $E$ such that $E \setminus \bigcup \cP$ is $\mu$-negligible.
\end{defin}

\begin{defin}\label{def:1'}
A \emph{restrictive Calder\'on--Zygmund space} is a metric-measure space $(X,d,\mu)$, where $\mu$ is a Radon measure on $X$, equipped with a family $\admR$ of Borel subsets of $X$, called \emph{admissible sets}, such that, for some constant $\kappa\in[1,\infty)$, the following properties hold.
\begin{enumerate}[label=(\arabic*)]
\item\label{it:1,def:1'} Every $R\in \admR$ has a \emph{centre} $x(R)\in X$ and a \emph{radius} $r(R)>0$ such that
\[
R\subseteq B_X(x(R),\kappa r(R)).
\]
\item\label{it:2,def:1'} $0 < \mu(R^\ast)\leq \kappa \mu(R) < \infty$ for all $R \in \admR$, where
\begin{equation}\label{eq:augment}
R^\ast \defeq \{x\in X \tc d(x,R) < r(R) \}.
\end{equation}
\item\label{it:3,def:1'} To each $R\in \admR$ is associated a collection $\chil(R)\subseteq \admR$ of \emph{children} of $R$ such that
\begin{itemize}
\item $\chil(R)$ is a quasi-partition of $R$ with $\#\chil(R)\leq \kappa$, and
\item $\mu(R)\leq \kappa\mu(R')$ for all $R'\in\chil(R)$.
\end{itemize}
We define iteratively $\chil^0(R) = \{ R \}$, $\chil^{n+1}(R) = \bigcup_{R' \in \chil^{n}(R)} \chil(R')$ for $n \in \NN$; the elements of $\chil^*(R) \defeq \bigcup_{n \in \NN} \chil^n(R)$ are called \emph{descendants} of $R$.
\item\label{it:4,def:1'} For all $\varepsilon>0$ and $R\in \admR$, and for almost all $x\in R$, there exists a descendant $R'$ of $R$ such that $x\in R'$ and $r(R')\leq \varepsilon$.
\item\label{it:5,def:1'} For all $\lambda\in (0,\mu(X))$ there exists a quasi-partition $\cP \subseteq \admR$ of $X$ such that $\mu(R)>\lambda$ for all $R \in \cP$.
\end{enumerate}
\end{defin}

The proof of the following result is analogous to that of \cite[Proposition 3.19]{MOV} and is omitted; we point out that the property \ref{it:4,def:1'} in Definition \ref{def:1'} plays the role of \cite[Lemma 3.18(v)]{MOV} in the proof of the inequality \eqref{eq:max_ineq}. Notice that we use the notation $\dashint_R |f| \,\dd\mu \defeq \mu^{-1}(R) \int_R |f| \,\dd\mu$ for the average of $|f|$ over $R$.

\begin{lm}
Let $(X,d,\mu)$ be a restrictive Calder\'on--Zygmund space with admissible family $\admR$.
Let $\cP \subseteq \admR$ be a quasi-partition of $X$, and let $\cP_* \defeq \bigcup_{R \in \cP} \chil^*(R)$. Then the maximal operator $M_{\cP}$ associated with $\cP$, defined by
\[
M_{\cP} f(x) = \begin{cases}
\sup\limits_{R \in \cP_* \tc x \in R} \dashint_R |f| \,\dd\mu &\text{if } x \in \bigcup \cP_*, \\
0 &\text{otherwise},
\end{cases}
\]
is of weak type $(1,1)$. Moreover
\begin{equation}\label{eq:max_ineq}
|f| \leq M_{\cP} f \ \text{$\mu$-a.e.}
\end{equation}
for all $f \in L^1_\loc(X)$.
\end{lm}

Following the proof of \cite[Theorem 3.20]{MOV}, we now show that Definition \ref{def:1'} is indeed more restrictive than Definition \ref{def:1}.

\begin{prop}
A restrictive Calder\'on--Zygmund space is a Calder\'on--Zygmund space.
\end{prop}
\begin{proof}
Let $(X,d,\mu)$ be a restrictive Calder\'on--Zygmund space with the admissible family $\admR$ and the constant $\kappa$. We shall prove that $(X,d,\mu)$ with $\admR$ and $\kappa$ satisfies Definition \ref{def:1}. It suffices to verify conditions \ref{it:3,def:1}--\ref{it:6,def:1}.
	
Fix $f\in L^1(X)$ and  $\lambda>\kappa \Vert f\Vert_{L^1(X)}/\mu(X)$. By \ref{it:5,def:1'} we can find a quasi-partition $\cP \subseteq \admR$ of $X$, such that every $R \in \cP$ satisfies $\mu(R) > \kappa \Vert f \Vert_{L^1(X)}/\lambda$. Thus $\dashint_R |f| \,\dd\mu <\lambda/\kappa$ for all $R \in \cP$. For any $R_0 \in \cP$, we say that a descendant $R$ of $R_0$ is \emph{stopping} if $\dashint_{R} |f| \,\dd\mu \geq\lambda/\kappa$, but $\dashint_{R'} |f| \,\dd\mu<\lambda/\kappa$ for any other $R'$ in the line of descendants from $R_0$ to $R$.
Let $\cS \subseteq \admR$ be the family of all stopping descendants of elements of $\cP$.

Let $\Omega=X\setminus\bigcup_{R\in\cS} R$. We set
\begin{equation*}
b_R = \left(f-\dashint_R f \,\dd\mu\right)\ind_R,\quad R\in \cS,
\qquad \text{and} \qquad 
g=f\ind_{\Omega}+ \sum_{R\in \cS} \left(\dashint_R f \,\dd\mu\right)\ind_R,
\end{equation*}
so that
\begin{equation}\label{eq:CZdec_R}
f = g + \sum_{R \in \cS} b_R
\end{equation}
pointwise.
Clearly, $\int b_R=0$ and, by \ref{it:3,def:1'}, $\dashint_R |f| \,\dd\mu \leq \lambda$ for all $R \in \cS$. Thus, $\Vert b_R\Vert_{L^1(X)}\leq 2\lambda\mu(R)$. By the definition of $\cS$ its elements are pairwise disjoint and
\begin{equation*}
	\sum_{R \in \cS} \mu(R) \leq \frac{\kappa}{\lambda} \sum_{R \in \cS} \int_R |f| \,\dd\mu \leq \frac{\kappa}{\lambda} \Vert f\Vert_{L^1(X)},
\end{equation*}
which proves \ref{it:4,def:1} and \ref{it:5,def:1}. Furthermore,
\[
|g|+\sum_{R \in \cS} |b_R| \leq |f| + 2 \lambda \ind_{X \setminus \Omega},
\]
where again we used that the $R \in \cS$ are pairwise disjoint; as $X \setminus \Omega$ has finite measure, this pointwise domination shows that, if $f \in L^p(X)$ for some $p \in [1,\infty)$, then the same is true for $g$ and the $b_R$, and moreover the sum in \eqref{eq:CZdec_R} converges in $L^p(X)$, thus proving \ref{it:6,def:1}.
	
For the remaining property \ref{it:3,def:1} notice that $\left|g|_R\right| \leq \dashint_R |f| \leq \lambda$ for any $R \in \cS$, and therefore $\|g|_{X \setminus \Omega}\|_\infty \leq \lambda$. On the other hand $\Vert g|_\Omega\Vert_{L^\infty}= \Vert f|_\Omega\Vert_{L^\infty} \leq \|M_\cP f|_\Omega\|_\infty$ by \eqref{eq:max_ineq}. By the definition of $\Omega$, for all $x \in \Omega$ and $R \in \cP_*$ with $R \ni x$, we have $\dashint_R |f| < \lambda/\kappa$, whence $M_\cP f(x) \leq \lambda/\kappa$. In conclusion, $\|g\|_{L^\infty} \leq \lambda$, as desired.
\end{proof}

We now recall from \cite{Va-PhD} the definition of Hardy and bounded mean oscillation spaces on a Calder\'on--Zygmund space.

\begin{defin}\label{def:(C)}
We say that a (restrictive) Calder\'on--Zygmund space satisfies condition (C) if there exists a subfamily $\admR' \subseteq \admR$ such that
\begin{itemize}
\item for each $R \in \admR$ there exists $R'\in \admR'$ such that $R \subseteq R'$;
\item if two sets from $\admR'$ are not disjoint, then one is contained in the other.
\end{itemize}
\end{defin}

\begin{defin}\label{defin:HardyBMO}
Let $(X,d,\mu)$ be a restrictive Calder\'on--Zygmund space satisfying condition (C), and let $\admR$ be the family of admissible sets.
\begin{enumerate}[label=(\alph*)]
\item An \emph{atom} is a function $a \in L^1(X)$ supported in some admissible set $R \in \admR$ and such that
\[
\int_X a \,\dd\mu = 0, \qquad \|a\|_2 \leq \mu(R)^{-1/2}.
\]

\item The Hardy space $H^1(X)$ is defined as the set of the functions $f \in L^1(X)$ admitting a decomposition
\begin{equation}\label{eq:atomic_dec}
f = \sum_{j \in \NN} \lambda_j a_j
\end{equation}
for some atoms $a_j$ and coefficients $\lambda_j \in \CC$ with $\sum_j |\lambda_j|<\infty$. The norm $\|f\|_{H^1(X)}$ of $f \in H^1(X)$ is defined as the infimum of the quantities $\sum_j |\lambda_j|$ over all decompositions \eqref{eq:atomic_dec}.

\item The bounded mean oscillation space $\BMO(X)$ is the quotient, modulo constant functions, of the space of the functions $f \in L^2_\loc(X)$ such that
\[
\|f\|_{\BMO(X)} \defeq \sup_{R \in \admR} \left( \frac{1}{\mu(R)} \int_R \left|f-\dashint_R f \,\dd\mu \right|^2 \,\dd\mu\right)^{1/2} < \infty.
\]
\end{enumerate}
\end{defin}

\begin{rem}
By \cite[Theorem 3.9]{Va-PhD}, we know that $\BMO(X)$ is the dual of $H^1(X)$ with respect to the natural pairing.
\end{rem}

We shall now show that $G_\nu$ is a restrictive Calder\'on--Zygmund space satisfying condition (C), so the general theory developed in \cite{HeSt,Va-PhD} can be applied. We start with a few technical results.

\begin{lm}\label{lm:tech}
The following hold.
\begin{enumerate}[label=(\roman*)]
\item\label{en:tech_dpr} Let $(x,u),(y,v)\in G_\nu$. Then
\begin{equation*}
|x-y| \leq e^{\min\{u,v\}} \sinh \dist((x,u),(y,v)),\qquad |u-v|\leq \dist((x,u),(y,v)).
\end{equation*}
\item\label{en:tech_dpr0} For any $(x,u)\in G_\nu$,
\begin{equation*}
|x|\leq \sinh|(x,u)|_\dist \leq \min \{ 1,|(x,u)|_\dist \} \, e^{|(x,u)|_\dist},\qquad |u|\leq |(x,u)|_\dist. 
\end{equation*}
\item\label{en:tech_rect}	Let $(x,u) \in G_\nu$ and $r>0$. Then,
\begin{equation*}
B_{G_\nu}((x,u),r) \subseteq [(x-e^u \sinh r)_+,x+e^u \sinh r)\times [u-r,u+r)\subseteq B_{G_\nu}((x,u),3r).
\end{equation*}
\end{enumerate}
\end{lm}
\begin{proof}
\ref{en:tech_dpr}. Set $r = \dist((x,u),(y,v))$. We may assume $r>0$. By \eqref{eq:dist},
\begin{equation*}
	\cosh(u-v)+\frac{(x-y)^2}{2e^{u+v}} = \cosh r.
\end{equation*}
Clearly, this yields $|u-v| \leq r$. Moreover,
\begin{equation*}
	|x-y|^2 = 2 e^{u+v} ( \cosh r - \cosh(u-v) ) = \frac{e^{v-u} - e^{-r}}{\sinh r}\left(2- \frac{e^{v-u} - e^{-r}}{\sinh r} \right) e^{2u}\sinh^2 r.
\end{equation*}
	Since $A(2-A)\leq 1$ for $A>0$, we obtain $|x-y| \leq e^u \sinh r$, and by symmetry
\begin{equation*}
	|x-y| \leq e^{\min\{u,v\}} \sinh r.
\end{equation*}

\ref{en:tech_dpr0}. Apply part \ref{en:tech_dpr} with $(y,v)=\bdzero$ and notice that $\sinh z\leq \min\{ 1,z\} e^z$ for $z \geq 0$.

\ref{en:tech_rect}. Fix $(x,u)\in G_\nu$ and $r>0$. For the first inclusion let $(y,v)\in B_{G_\nu}((x,u),r)$. Part \ref{en:tech_dpr} immediately gives
\begin{equation*}
|y-x| < e^{u}\sinh r ,\qquad |v-u|<r.
\end{equation*}

For the second inclusion, let $(y,v) \in G_\nu$ be such that
\begin{equation*}
|y-x|\leq e^u \sinh r,\qquad |v-u|\leq r.
\end{equation*}
Then
\begin{multline*}
\dist((x,u),(y,v)) =\cosh(u-v)+\frac{(x-y)^2}{2e^{u+v}} \\
\leq \cosh r + \frac{e^{u-v}\sinh^2 r}{2}\leq \cosh r +\frac{e^r\sinh^2 r}{2}< \cosh(3r),
\end{multline*}
where the last inequality follows from the identity
\[
\cosh(3r) - \cosh r = 4\sinh^2 r \, \cosh r;
\]
so $(y,v) \in B_{G_\nu}((x,u),3r)$.
\end{proof}

\begin{defin}\label{def:2}
	Let $m\in\NN$, $l\in\ZZ$, $u\in\RR$, and $r>0$. We define
	\begin{equation*}
	R_{m,l,u,r}= [m 2^l, (m+1) 2^l) \times [u-r,u+r).
	\end{equation*} 
	We say that a set $R_{m,l,u,r}$ is admissible if it satisfies the admissibility condition
	\begin{equation}\label{eq:adm}
	e^{u} \sinh(2r) \leq 2^{l-1} < e^u \sinh(9r).
	\end{equation}
	We denote the family of all admissible sets $R_{m,l,u,r}$ by $\adm$.
\end{defin}

\begin{thm}\label{thm:2}
	For any $\nu\geq 1$, $G_\nu$ with the family of sets $\adm$ is a (restrictive) Calder\'on--Zygmund space satisfying condition (C). 
\end{thm}
\begin{proof}
We shall justify that $G_\nu$ with $\kappa=\max\{27,2^{\nu+1}\}$ and $\adm$ satisfies \ref{it:1,def:1'}--\ref{it:5,def:1'} of Definition \ref{def:1'} and the condition (C) of Definition \ref{def:(C)}.
	
Fix $R\defeq R_{m,l,u,r}\in\adm$. By the admissibility condition \eqref{eq:adm} and Lemma \ref{lm:tech},
\begin{align*}
	R &\subseteq [(m2^l+2^{l-1}-e^u\sinh(9r))_+,m2^l+2^{l-1}+e^u\sinh(9r)) \times [u-r,u+r)\\
	&\subseteq B_{G_\nu}((m2^l+2^{l-1},u),27r),
\end{align*}
which proves \ref{it:1,def:1'} with
\[
x(R) \defeq (m2^l+2^{l-1},u), \qquad r(R)\defeq r.
\]
	
For \ref{it:2,def:1'} observe that, by \eqref{eq:augment} and Lemma \ref{lm:tech}\ref{en:tech_rect},
\begin{multline*}
	R^\ast 
	= \bigcup_{(y,v)\in R} B_{G_\nu}((y,v),r)
	\subseteq \bigcup_{(y,v)\in R} [(y-e^v\sinh r)_+,y+e^v\sinh r ) \times [v-r,v+r)\\
	\subseteq [(m 2^l - e^{u} \sinh(2r))_+,(m+1)2^l +e^{u}\sinh(2r) ) \times [u-2r,u+2r),
\end{multline*}
where the fact that 
\[
e^r \sinh r \leq 2 \cosh r \, \sinh r = \sinh(2r)
\]
was used. Since $e^u \sinh(2r)\leq 2^{l-1}$ by the admissibility condition \eqref{eq:adm}, we obtain
\begin{equation*}
	\mu_\nu(R^\ast)\leq\frac{4 r 2^{\nu l}}{\nu} ( (m+3/2)^\nu -(m-1/2)_+^\nu ),
\end{equation*}
whereas
\begin{equation}\label{eq:meas_R}
	\mu_\nu(R)=\frac{2 r 2^{\nu l}}{\nu} ( (m+1)^\nu -m^\nu )
	=\frac{2 r 2^{\nu l} }{\nu 2^{\nu}} ( (2m+2)^\nu -(2m)^\nu ).
\end{equation}
Thus, by the monotonicity of $x\mapsto (x+2)^{\nu}-x^\nu$,
\begin{equation*}
	\mu_\nu(R^\ast)\leq 2^{\nu+1} \mu_\nu(R).
\end{equation*}
	
Now we pass to verifying \ref{it:3,def:1'}. Here we need to decompose $R = R_{m,l,u,r} \in \adm$ as a disjoint union of family $\chil(R)$ of admissible sets, which we declare to be the children of $R$. We consider two potential decompositions,
\[
R = R_{1,1} \cup R_{1,2} = R_{2,1} \cup R_{2,2},
\]
corresponding to splitting $R$ with respect to the first or the second component, i.e.,
\begin{align*}
R_{1,1} &\defeq R_{2m,l-1,u,r}, & R_{1,2} &\defeq R_{2m+1,l-1,u,r} & \text{(first type)},\\
R_{2,1} &\defeq R_{m,l,u-r/2,r/2}, & R_{2,2} &\defeq R_{m,l,u+r/2,r/2}  & \text{(second type)}.
\end{align*}
	
Observe that since $R\in\adm$ we have $R_{1,1},R_{1,2}\in\adm$ whenever
\begin{equation}\label{eq:10}
	e^u \sinh 2r \leq 2^{l-2};
\end{equation}
so, when \eqref{eq:10} is satisfied, we declare $R_{1,1}$ and $R_{1,2}$ to be the children of $R$.
In the opposite case, denoting $u'=u\pm r/2$, we have
\begin{equation*}
	2^{l-1}<2 e^u \sinh(2r) \leq 2 e^{u'+r/2}\sinh(2r) \leq e^{u'} \sinh(9r/2),
\end{equation*}
where in the last inequality we used the simple fact that
\begin{equation*}
	2 e^a \sinh(b/2) \leq  e^a \sinh b  \leq \sinh(a+b),\qquad a,b>0.
\end{equation*}
Moreover, by the admissibility of $R$,
\begin{equation*}
	e^{u'}\sinh r \leq e^{u} \sinh(2r) \leq 2^{l-1}.
\end{equation*}
Thus, if \eqref{eq:10} fails, then $R_{2,1},R_{2,2}\in\adm$; in this case, we declare $R_{2,1},R_{2,2}$ to be the children of $R$.
	
To conclude the proof of \ref{it:3,def:1'} notice that
\begin{equation*}
	\mu_\nu(R_{2,j})= \frac{1}{2} \mu_\nu(R),\qquad j=1,2,
\end{equation*}
and
\begin{equation*}
	\mu_\nu(R_{1,1})= \frac{2r 2^{\nu l}}{\nu 2^\nu} ((2m+1)^\nu-(2m)^\nu ),\qquad \mu_\nu(R_{1,2})= \frac{2r 2^{\nu l}}{\nu 2^\nu} ((2m+2)^\nu-(2m+1)^\nu ).
\end{equation*}
Therefore, by \eqref{eq:meas_R} and the monotonicity of $x\mapsto (x+1)^\nu -x^\nu$,
\begin{equation*}
	\mu_\nu(R) \leq 2^{\nu} \mu_\nu(R_{i,j}),\qquad i,j=1,2.
\end{equation*}
	
In order to justify \ref{it:4,def:1'} it suffices to notice that, for any $R\in\adm$, if the children of $R$ are of the second type, then $r(R')=r(R)/2$ for $R'\in \chil(R)$. But the children can be of the first type only for a finite number of generations in a row, since \eqref{eq:10} eventually will be false. Moreover, $R$ is clearly a disjoint union of all descendants of any fixed generation.
	
We now prove \ref{it:5,def:1'}. Fix $r_0>0$ and choose the smallest $l_0\in\ZZ$ such that
\begin{equation*}
	\sinh(2r_0) \leq 2^{l_0-1}<\sinh(9r_0).
\end{equation*}
Thus, $\{R_{m,l_0,0,r_0}\}_{m\in\NN}\subseteq\adm$ and 
\begin{equation}\label{eq:dec0}
\bigcup_{m\in\NN} R_{m,l_0,0,r_0}=\Rnon\times[-r_0,r_0).
\end{equation}
	
Set $u_0 = 0$. We now define inductively, for any $n \geq 1$,
\[
r_{n} = 3r_{n-1}, \quad u_{n} = u_{n-1} + r_{n-1} + r_n
\]
and let $l_n^\pm \in \ZZ$ be the smallest integers such that
\[
	e^{\pm u_n}\sinh(2r_n) \leq 2^{l_n^\pm -1}<e^{\pm u_n}\sinh(9r_n).
\]
Again, $\{ R_{m,l_n^\pm,\pm u_n,r_n} \}_{m \in \NN} \subseteq \adm$ and
\begin{equation}\label{eq:decn}
\bigcup_{m \in \NN} R_{m,l_n^\pm,\pm u_n,r_n} = \Rnon\times [\pm u_n-r_n,\pm u_n+r_n).
\end{equation}
Moreover, since $u_n \geq u_{n-1}$, $r_n \geq r_{n-1}$, clearly $e^{u_n}\sinh (2r_n) \geq e^{u_{n-1}}\sinh(2r_{n-1})$ and $l_n^+ \geq l_{n-1}^+$. On the other hand, $u_n = u_{n-1} + 4r_n/3$, so
\begin{equation*}
	e^{-u_{n-1}}\sinh(2r_{n-1}) = e^{-u_n} e^{4r_n/3} \sinh (2r_n/3) \leq e^{-u_n} \sinh(2r_n),
\end{equation*}
which gives $l^{-}_{n}\geq l^{-}_{n-1}$.

Since $u_n-r_n=u_{n-1}+r_{n-1}$, from \eqref{eq:dec0} and \eqref{eq:decn} we deduce that $G_\nu = \Rnon\times\RR$ is the disjoint union of the family $\{R_{m,l^\pm_n,\pm u_n,r_n} \}_{m,n\in\NN} \subseteq \adm$.
Moreover, notice that
\begin{equation*}
	\mu_\nu(R_{m,l_n^{\pm},\pm u_n,r_n})=\frac{2r_n 2^{\nu l_n^\pm}}{\nu} ( (m+1)^\nu-m^\nu)
\end{equation*}
is increasing in both $n$ and $m$. Thus,
\begin{equation*}
	\mu_\nu(R_{m,l_n^{\pm},\pm u_n,r_n})\geq \mu_\nu(R_{0,l_0,0,r_0})=\frac{2r_02^{\nu l_0}}{\nu}.
\end{equation*}
By an appropriate choice of $r_0$ we can make this measure arbitrarily large. This completes the proof of \ref{it:5,def:1'}.
	
Finally, we discuss the condition (C) of Definition \ref{def:(C)}. Let $\{r_l\}_{l\in\ZZ}$ be a sequence of positive real numbers satisfying $\sinh(2r_l) = 2^{l-1}$. Thus, the role of family $\admR'$ can be played by $\{R_{0,l,0,r_l} \}_{l\in\ZZ}$. One can immediately verify the necessary assumptions.
\end{proof}

Thanks to Theorem \ref{thm:2}, we can apply to $G_\nu$ the singular integral theory  on Calder\'on--Zygmund spaces developed in \cite{HeSt,Va-PhD}. In particular, we have the following result concerning the boundedness properties of a class of singular integral operators, which combines \cite[Theorem~1.2]{HeSt} and \cite[Theorem 3.10]{Va-PhD}.

\begin{thm}\label{thm:3}
Let $T$ be a linear operator bounded on $L^2(G_\nu)$.
Assume that
\[
\langle T f,g \rangle_{L^2(G_\nu)} = \sum_{j \in \ZZ} \langle T_j f,g \rangle_{L^2(G_\nu)}
\]
for all $f,g \in L^2(G_\nu)$ with disjoint compact supports, where the $T_j$ are integral operators with kernels $\iker_j$
satisfying for some positive constants $B,c,\varepsilon$, with $c\neq 1$, the following conditions:
\begin{equation}\label{eq:64}
\begin{aligned}
 \sup_{\bdy\in G_\nu}\int_{G_\nu} |\iker_j(\bdx,\bdy)| \,  (1+c^j \dist(\bdx,\bdy))^\varepsilon \,\dmub(\bdx) &\leq B, \\
 \int_{G_\nu} |\iker_j(\bdx,\bdy)-\iker_j(\bdx,\bdy')|\, \dmub(\bdx) &\leq B \dist(\bdy,\bdy') c^j \quad\forall \bdy,\bdy' \in G_\nu.
\end{aligned}
\end{equation}
Then, $T$ is of weak type $(1,1)$, bounded on $L^p(G_\nu)$, $p\in(1,2]$, and bounded from $H^1(G_\nu)$ to $L^1(G_\nu)$.
\end{thm}

Lemma \ref{lm:10} immediately gives the following corollary.

\begin{cor}\label{cor:2}
If the $\iker_j$ are continuously differentiable with respect the second variable, then the second line in \eqref{eq:64} can be replaced by the stronger condition
\begin{equation*}
	\sup_{\bdy\in G_\nu}\int_{G_\nu} | \nabla_\nu^{\bdy} \iker_j(\bdx,\bdy) | \,\dmub(\bdx) \leq Bc^j.
\end{equation*}
\end{cor}

We can also state a version of Theorem \ref{thm:3} for $\diamond_\nu$-convolution operators, which follows immediately from Lemma \ref{lm:l1est_int_conv} and Corollary \ref{cor:2}.

\begin{cor}\label{cor:4}
Assume that the operators $T_j$ from Theorem \ref{thm:3} are right $\diamond_\nu$-convolution operators, i.e.,
\[
T_j f = f \diamond_\nu K_j,
\]
with convolution kernels $K_j : G_\nu \to \CC$. Then in place of \eqref{eq:64} it suffices to assume that 
\begin{equation}\label{eq:65}
\begin{aligned}
\int_{G_\nu} | K_j(\bdx)| \, (1+c^j |\bdx|_\dist)^\varepsilon \,\dmub(\bdx) &\leq B,\\
\int_{G_\nu} |\nabla_\nu K^\ast_j(\bdx)|\, \dmub(\bdx) &\leq B c^j.
\end{aligned}
\end{equation}
\end{cor}

We can also prove a variant of Corollary \ref{cor:4} for $\diamond_\nu^\rD$-convolution operators. The following statement makes use of the notation of Section \ref{ss:dunkl}, as well as the notation $|\bdx|_\dist \defeq |\tilde\bdx|_\dist$ for all $\bdx \in G^\rD_\nu$. We shall further write $\ext_\even f$ and $\ext_\odd f$ for the even and odd extensions to $G_\nu^\rD$ of a function $f$ on $G_\nu$.

\begin{cor}\label{cor:6}
Assume that the operators $T_j$ from Theorem \ref{thm:3} are given by
\begin{equation}\label{eq:Tj_diamD}
T_j f = ((\ext_{\sigma_j} f) \diamond_\nu^\rD K_j)|_{G_\nu}
\end{equation}
for certain kernels $K_j : G^\rD_{\nu} \to \CC$ and $\sigma_j \in \{\even,\odd\}$. Then in place of \eqref{eq:64} 
it suffices to assume that
\begin{equation}\label{eq:70}
\begin{aligned}
\int_{G_\nu^\rD} | K_j(\bdx)| \, (1+c^j |\bdx|_\dist)^\varepsilon \,\dmub^\rD(\bdx) &\leq B,\\
\int_{G_\nu^\rD} ( |\nabla_\nu K^{\bullet}_j(\bdx) | + |\tilde{U} (K^{\bullet}_j)_\odd|)\, \dmub^\rD(\bdx) &\leq B c^j.
\end{aligned}
\end{equation}
\end{cor}
\begin{proof}
From \eqref{eq:Tj_diamD} and \eqref{eq:Ddiamond2} we see that, if we set $\epsilon_\pm^\even \defeq 1$ and $\epsilon_\pm^\odd \defeq \pm 1$, then the integral kernels $\iker_j : G_\nu \times G_\nu \to \CC$ of the operators $T_j$ are given by
\[\begin{split}
\iker_j(\bdx,\bdy) 
&= \frac{1}{2} \sum_{\pm} \epsilon_\pm^{\sigma_j} \modu_\nu((\pm y,v)) r^\rD_\bdx K_j((\pm y,v)^{-1}) \\
&= \frac{1}{2} \sum_{\pm} \epsilon_\pm^{\sigma_j} \modu_\nu(\bdx) \overline{r^\rD_{(\pm y,v)}K_j^{\bullet}(\bdx^{-1} )},
\end{split}\]
for all $\bdx \in G_\nu$ and $\bdy=(y,v) \in G_\nu$, where the second equality is due to Lemma \ref{lm:conv_inv_D}\ref{en:conv_inv_D_magic}. Thus, the operators $T_j$ satisfy the assumptions of Theorem \ref{thm:3} (strengthened as in Corollary \ref{cor:2}) provided
\begin{equation}\label{eq:71}
\begin{aligned}
\sup_{\bdy \in G^\rD_\nu} \int_{G^\rD_\nu \setminus G_\nu} | r^\rD_{\bdy} K^{\bullet}_j(\bdx)| \, (1+c^j \dist(\bdx^{-1},\tilde\bdy))^\varepsilon \,\dmub^\rD(\bdx) &\leq B,\\
\sup_{\bdy\in G^\rD_\nu} \int_{G^\rD_\nu \setminus G_\nu} | \nabla_\nu^\bdy r^\rD_{\bdy}K^{\bullet}_j(\bdx)|\, \dmub^\rD(\bdx) &\leq B c^j.
\end{aligned}
\end{equation}
We shall justify that \eqref{eq:70} implies \eqref{eq:71}.
	
Firstly, fix $\bdy = (y,v) \in G_\nu^\rD$ and set $g(\bdx)\defeq (1+c^j\dist(\bdx,\tilde\bdy))^\varepsilon$. By \eqref{eq:72} and \eqref{eq:diamond},
\[\begin{split}
\int_{G^\rD_\nu\setminus G_\nu} |r^\rD_{\bdy} K^{\bullet}_j(\bdx)| g(\bdx^{-1})\,\dmub^\rD(\bdx)
&\leq \int_{G_\nu}  r_{\tilde\bdy} ( |(K^{\bullet}_j)_\even(\bdx)|+|(K^{\bullet}_j)_\odd(\bdx)|) g(\bdx^{-}) \,\dmub(\bdx) \\
&= g \diamond_\nu ( |(K^{\bullet}_j)_\even|+|(K^{\bullet}_j)_\odd|) (\tilde\bdy).
\end{split}\]
Thus, by arguing much as in \eqref{eq:dist_conv_est} and using Lemma \ref{lm:11}\ref{lm:11(3)}, we obtain
\[\begin{split}
&\int_{G^\rD_\nu \setminus G_\nu} | r^\rD_{\bdy} K^{\bullet}_j(\bdx)| (1+c^j \dist(\bdx^{-1},\tilde\bdy))^\varepsilon \,\dmub^\rD(\bdx) \\
&\leq \int_{G_\nu} ( |(K^{\bullet}_j)_\even(\bdx) |+|(K^{\bullet}_j)_\odd(\bdx)| )  \, (1+c^j|\bdx|_\dist)^\varepsilon\,\dmub(\bdx) \\
&\leq 2 \int_{G^\rD_\nu}  |K_j(\bdx)|  \, (1+c^j|\bdx|_\dist)^\varepsilon\,\dmub^\rD(\bdx).
\end{split}\]
On the other hand, by Lemma \ref{lm:12},
\[\begin{split}
&\int_{G^\rD_\nu \setminus G_\nu} |\nabla_\nu^\bdy r^\rD_{\bdy}K^{\bullet}_j(\bdx)| \, \dmub^\rD(\bdx) \\
&\leq \int_{G_\nu} r_{\tilde\bdy}(|\nabla_\nu (K^{\bullet}_j)_\even|+|\nabla_\nu (K^{\bullet}_j)_\odd|+2| \tilde{U} (K^{\bullet}_j)_\odd|)(\bdx) \, \dmub(\bdx) \\
&\leq \int_{G_\nu} (|\nabla_\nu (K^{\bullet}_j)_\even|+|\nabla_\nu (K^{\bullet}_j)_\odd|+2| \tilde{U} (K^{\bullet}_j)_\odd|)(\bdx) \, \dmub(\bdx) \\
&\leq 2\int_{G_\nu^\rD} (|\nabla_\nu K^{\bullet}_j|+| \tilde{U} (K^{\bullet}_j)_\odd|)(\bdx) \, \dmub^\rD(\bdx).
\end{split}\]
These estimates show that, up to relabelling $2B$ as $B$, the assumptions \eqref{eq:70} imply \eqref{eq:71}.
\end{proof}

\section{Heat kernel estimates}\label{s:heatestimates}

By using  finite propagation speed and the methods of \cite{CouSi} we derive a Gaussian-type pointwise upper bound for the heat kernel on $G_\nu$.

\begin{lm}\label{lm:2}
Let $\epsilon\in (0,\infty)$. Then
\begin{equation*}
0\leq K_{e^{-t\Delta_\nu}}(x,u) 
\lesssim_{\nu,\epsilon} e^{-\nu u/2} t^{[-\frac{\nu+1}{2},-\frac{3}{2}]} \exp\left( -\frac{|(x,u)|_\dist^2}{(4+\epsilon)t}\right) \quad \forall (x,u)\in G_\nu,\ t>0.
\end{equation*}
Moreover,
\[
\|K_{e^{-t\Delta_\nu}}\|_{L^1(G_\nu)} = 1, \qquad \|K_{e^{-t\Delta_\nu}}\|_{L^2(G_\nu)} \lesssim_\nu t^{[-(\nu+1)/4,-3/4]} \qquad \forall t > 0.
\]
\end{lm}
\begin{proof}
It was checked in \cite[Proposition 4.4]{MaPl} that $K_{e^{-t\Delta_\nu}}$ is nonnegative and $L^1$-normalised. On the other hand, by Proposition \ref{prop:6},
\[
\Vert K_{e^{-t\Delta_\nu}}\Vert_{L^2(G_\nu)}^2 \lesssim_\nu \int_0^\infty e^{-2t\lambda} \lambda^{[\frac{3}{2},\frac{\nu+1}2]}\,\ddlam \simeq_\nu t^{[-\frac{\nu+1}2,-\frac{3}{2}]},
\]
and a similar computation, using Lemma \ref{lm:1}, gives, for all $u \in \RR$, that
\begin{equation}\label{eq:2}
\Vert K_{e^{-t\Delta_\nu}}(\cdot,u)\Vert_{L^\infty(X_\nu)} 
\lesssim_\nu e^{-\nu u/2} \int_0^\infty e^{-t\lambda} \lambda^{[\frac{3}{2},\frac{\nu+1}2]}\,\ddlam 
\simeq_\nu e^{-\nu u/2} t^{[-\frac{\nu+1}2,-\frac32]}.
\end{equation}
	
Let $\iheat_t$ denote the integral kernel of $e^{-t\Delta_\nu}$.
By \eqref{eq:26}, \eqref{eq:transl} and \eqref{eq:hankel_transl} we have the following relation between convolution and integral kernels:
\begin{equation*}
	\iheat_{t}((x,u),(y,v)) = e^{-\nu v}\tau_\nu^{[e^{-v} y]} K_{e^{-t\Delta_\nu}}(e^{-v}x,u-v)
\end{equation*}
for all $(x,u),(y,v) \in G_\nu$. In particular,
\begin{equation}\label{eq:12}
	\iheat_{t}(\bdx,\bdzero) = K_{e^{-t\Delta_\nu}}(\bdx) \qquad \forall \bdx \in G_\nu.
\end{equation}
	
Since the Hankel translations are contractions on $L^\infty(X_\nu)$ and preserve positivity, by applying \eqref{eq:2} we obtain
\begin{equation}\label{eq:4}
	0\leq \iheat_{t}((x,u),(y,v)) \leq e^{-\nu v} \Vert K_{e^{-t\Delta_\nu}}(\cdot,u-v)\Vert_{L^\infty(X_\nu)}\lesssim_\nu e^{-\nu(u+v)/2} t^{[-\frac{\nu+1}{2},-\frac{3}{2}]}
\end{equation}
for all $(x,u),(y,v) \in G_\nu$. In particular, we have the on-diagonal estimate
\begin{equation*}
0\leq \iheat_{t}(\bdx,\bdx) \lesssim_\nu e^{-\nu u} t^{[-\frac{\nu+1}{2},-\frac{3}{2}]} \qquad \forall \bdx \in G_\nu;
\end{equation*}
in other words, $\iheat_t$ satisfies \cite[eq.\ (4.12)]{CouSi} with $V((x,u),s) \defeq e^{\nu u}s^{[\nu+1,3]}$. Moreover, the function $V$ satisfies the doubling condition \cite[eq.\ (4.11)]{CouSi} with $\delta=\min\{\nu+1,3\}$; this follows from the bound
\begin{equation*}
\frac{s^{[\nu+1,3]}}{r^{[\nu+1,3]}}\leq \left( \frac{s}{r}\right)^{\min\{\nu+1,3\}},\qquad 0<s\leq r.
\end{equation*}
Furthermore, since $\Delta_\nu$ has finite propagation speed (see \cite[Lemma 4.17]{MaPl}), by \cite[Theorem~3.4]{CouSi} it also satisfies the Davies--Gaffney estimate \cite[eq.\ (3.2)]{CouSi}. Thus, \cite[Theorem 4.4]{CouSi} implies an improved version of \eqref{eq:4}: for any $\epsilon > 0$,
\begin{equation*}
	0\leq \iheat_{t}((x,u),(y,v))\lesssim_{\nu,\epsilon} e^{-(u+v)\nu/2} t^{[-\frac{\nu+1}{2},-\frac{3}{2}]}\exp\left(-\frac{\dist((x,u),(y,v))^2}{(4+\epsilon)t} \right)
\end{equation*}
for all $(x,u),(y,v) \in G_\nu$. 
Finally, in light of \eqref{eq:12}, taking $(y,v)=0$ proves the lemma.
\end{proof}

We also record here an easy weigthed $L^1$-estimate for derivatives of the heat kernel on $X_\nu$, which follows immediately by homogeneity.

\begin{lm}\label{lm:3}
Let $\alpha_0 > 0$. Then, for all $n \in \NN$ and $\alpha \in [0,\alpha_0]$,
\begin{equation}\label{eq:5}
\int_0^\infty x^\alpha |\partial_x^n K_{e^{-tL_\nu}}(x)|\,\dmu(x) \simeq_{\nu,n,\alpha_0} t^{(\alpha-n)/2} \qquad \forall t>0.
\end{equation} 
\end{lm}
\begin{proof}
From \eqref{eq:13} we see that
$K_{e^{-tL_\nu}}(x) =t^{-\nu/2} K_{e^{-L_\nu}}(x/t^{1/2})$
and
\[
\int_0^\infty x^\alpha |\partial_x^n K_{e^{-tL_\nu}}(x)|\,\dmu(x) = t^{(\alpha-n)/2} \int_0^\infty x^\alpha |\partial_x^n K_{e^{-L_\nu}}(x)|\,\dmu(x)
\]
by homogeneity considerations; moreover, the latter integral remains bounded as long as $\alpha$ does, because $K_{e^{-L_{\nu}}} \in \Sz_\even(\Rnon)$.
\end{proof}

We now obtain some formulas for derivatives of the heat kernel on $G_\nu$ with respect to the vector fields $\vfX_0$ and $\vfX_1$ of \eqref{eq:vfs}.

\begin{lm}\label{lm:4}
For all $\nu \geq 1$ and $\bdx = (x,u) \in G_\nu$,
\begin{align*}
\vfX_0 K_{e^{-t\Delta_\nu}}(\bdx) &= \frac{\nu}{2} (x^2-2e^u\sinh u ) K_{e^{-t\Delta_{\nu+2}}}(\bdx) -\frac{\nu}{2} K_{e^{-t\Delta_{\nu}}}(\bdx), \\
\vfX_1 K_{e^{-t\Delta_\nu}}(\bdx) &=- \nu x e^u K_{e^{-t\Delta_{\nu+2}}}(\bdx).
\end{align*}
\end{lm}
\begin{proof}
From 
Proposition \ref{prop:heat_dim}\ref{en:heat_dim_der} we see that, if $r = |\bdx|_\dist$, then
\[
\vfX_j (\modu_\nu^{1/2} H_{e^{-t\Delta_\nu}}(r))
= (\vfX_j \modu_\nu^{1/2}) H_{e^{-t\Delta_\nu}}(r) - \modu_\nu^{1/2} H_{e^{-t\Delta_{\nu+2}}}(r) (\sinh r) \vfX_j r,
\]
that is, by \eqref{eq:H_notat},
\[
\vfX_j K_{e^{-t\Delta_\nu}} = \frac{\vfX_j \modu_\nu^{1/2}}{\modu_\nu^{1/2}} K_{e^{-t\Delta_\nu}} - \nu e^u (\vfX_j \cosh r)  K_{e^{-t\Delta_{\nu+2}}} .
\]
On the other hand, as $\modu_\nu^{1/2}(x,u) = e^{-\nu u/2}$ and $\cosh r = \cosh u + e^{-u} x^2/2$ by \eqref{eq:modular} and \eqref{eq:dist}, from \eqref{eq:vfs} we see that
\begin{equation}\label{eq:der_modu}
\begin{aligned}
\frac{\vfX_0 \modu_\nu^{1/2}}{\modu_\nu^{1/2}} &= -\nu/2, & \frac{\vfX_1 \modu_\nu^{1/2}}{\modu_\nu^{1/2}} &= 0, \\
\vfX_0 \cosh r &= \sinh u - e^{-u} x^2/2, & \vfX_1 \cosh r &= x,
\end{aligned}
\end{equation}
and the desired formulas follow.
\end{proof}

We can now prove the main result of this section, namely, weighted $L^1$-estimates for the heat kernel on $G_\nu$ and its gradient, which crucially are valid both for small and large times.

\begin{prop}\label{prop:2}
Let $\varepsilon\in[0,\infty)$. For all $t>0$,
\begin{align*}
\left\Vert e^{\varepsilon |\cdot|_\dist/\sqrt{t}} K_{e^{-t\Delta_\nu}} \right\Vert_{L^1(G_\nu)} &\lesssim_{\varepsilon,\nu} 1,\\
\left\Vert e^{\varepsilon |\cdot|_\dist/\sqrt{t}} | \nabla_{\nu}K_{e^{-t\Delta_\nu}} |\right\Vert_{L^1(G_\nu)} &\lesssim_{\varepsilon,\nu} t^{-1/2}.
\end{align*}
\end{prop}
\begin{proof}
Clearly the second bound above is equivalent to the two bounds
\begin{equation}\label{eq:8}
\left\Vert e^{\varepsilon |\cdot|_\dist/\sqrt{t}} \vfX_j K_{e^{-t\Delta_\nu}}\right\Vert_{L^1(G_\nu)} 
\lesssim_{\varepsilon,\nu} t^{-1/2},\qquad j=0,1.
\end{equation}
	
Fix $\varepsilon\in[0,\infty)$. We first discuss the small-time bounds. By applying Lemma \ref{lm:2} with $\epsilon=1$,
\begin{equation}\label{eq:9}
\begin{split}
&\left\Vert e^{\varepsilon |\cdot|_\dist/\sqrt{t}} K_{e^{-t\Delta_\nu}} \right\Vert_{L^1(G_\nu)} \\
&\lesssim_\nu t^{[-\frac{\nu+1}{2},-\frac{3}{2}]} \int_{G_\nu} e^{-\nu u/2} \exp\left( -\frac{|(x,u)|_\dist^2}{5t} + \frac{\varepsilon |(x,u)|_\dist}{\sqrt{t}} \right) \,\dmu(x)\,\du \\
&\lesssim_\nu   t^{[-\frac{\nu+1}{2},-\frac{3}{2}]} \int_0^\infty \exp\left(-\frac{r^2}{5t} + \frac{\varepsilon r}{\sqrt{t}} +\frac{\nu}{2} r \right) r^\nu \,\dd r;
\end{split}
\end{equation}
in the last estimate, we used the 
integration formula from Lemma \ref{lm:13}.
In particular, for all $t_0>0$, from \eqref{eq:9} we deduce that
\[
\left\Vert e^{\varepsilon |\cdot|_\dist/\sqrt{t}} K_{e^{-t\Delta_\nu}} \right\Vert_{L^1(G_\nu)} \lesssim_{\nu,\varepsilon,t_0} t^{-(\nu+1)/2} \int_0^\infty \exp\left(-\frac{r^2}{10t}\right) r^\nu \,\dd r \simeq_\nu 1 \quad\forall t \in (0,t_0].
\]
As for the gradient bound, observe that, by Lemma \ref{lm:tech}\ref{en:tech_dpr0},
\[
xe^u, x^2, e^u \left|\sinh u\right| \leq |(x,u)|_\dist\, e^{2|(x,u)|_\dist} \qquad \forall (x,u)\in G_\nu.
\]
Thus, by Lemma \ref{lm:4},
\begin{multline*}
\left\Vert e^{\varepsilon |\cdot|_\dist/\sqrt{t}} |\nabla_\nu K_{e^{-t\Delta_\nu}} |\right\Vert_{L^1(G_\nu)} \\
\lesssim_{\nu} \left\Vert e^{\varepsilon |\cdot|_\dist/\sqrt{t}} K_{e^{-t\Delta_\nu}}\right\Vert_{L^1(G_\nu)} + \left\Vert |\cdot|_\dist\, e^{(2+\varepsilon/\sqrt{t}) |\cdot|_\dist} K_{e^{-t\Delta_{\nu+2}}}\right\Vert_{L^1(G_\nu)}.
\end{multline*}
We already know that the first summand is bounded for $t \in (0,t_0]$. As for the second summand, much as in \eqref{eq:9} we estimate
\[\begin{split}
&\left\Vert |\cdot|_\dist\, e^{(2+\varepsilon/\sqrt{t}) |\cdot|_\dist} K_{e^{-t\Delta_{\nu+2}}}\right\Vert_{L^1(G_\nu)}\\
&\lesssim_{\nu} t^{[-(\nu+3)/2,-3/2]} \int_0^\infty \exp\left(-\frac{r^2}{5t} +\frac{\varepsilon r}{\sqrt{t}} +\frac{\nu+6}{2} r  \right) r^{\nu+1} \,\dd r\\
&\lesssim_{\nu,\varepsilon,t_0} t^{-(\nu+3)/2} \int_0^\infty \exp\left(-\frac{r^2}{10} \right) r^{\nu+1} \,\dd r \simeq_\nu t^{-1/2}
\end{split}\]
for all $t \in (0,t_0]$. About the first inequality, notice that the estimate for $K_{e^{-t\Delta_{\nu+2}}}$ from Lemma \ref{lm:2} involves a factor $\modu_{\nu+2}^{1/2}$, while the integration formulas in Lemma \ref{lm:13} involve $\modu_{\nu}^{1/2}$; nevertheless, the ratio $e^{-u}$ of the two factors is simply bounded by $e^{|(x,u)|_{\dist}}$ by Lemma \ref{lm:tech}\ref{en:tech_dpr0}.

This confirms the validity of the desired bounds \eqref{eq:8} for small time $t \in (0,t_0]$, for any choice of $t_0 > 0$, with a constant depending on $t_0$.

We now set $t_0=4\varepsilon^2+1$, and discuss the bounds for large time $t \geq t_0$.
If $\gamma_t \defeq \varepsilon/\sqrt{t}$, then $\gamma_t \in [0,1/2)$ due to our choice of $t_0$. Consequently,
\begin{equation*}
	e^{\gamma_t |(x,u)|_\dist}
	\leq (2\cosh|(x,u)|_\dist)^{\gamma_t}
	\simeq (\cosh u)^{\gamma_t} + (e^{-u} x^2)^{\gamma_t} \qquad \forall (x,u)\in G_\nu,
\end{equation*}
with implicit constants independent of $t \geq t_0$ and $\varepsilon$.

Let $\ell\in\{0,1\}$. By \eqref{eq:7} and \eqref{eq:5},
\[\begin{split}
&\left\Vert e^{\varepsilon |\cdot|_\dist/\sqrt{t}} \vfX_1^\ell K_{e^{-t\Delta_\nu}} \right\Vert_{L^1(G_\nu)} \\
&\lesssim \int_0^\infty |\Psi_t(\xi)| \int_\RR e^{\ell u}  \exp\left(-\frac{\cosh u}\xi \right) \\
&\quad\times \int_0^\infty ( (\cosh u)^{\gamma_t} + (e^{-u} x^2)^{\gamma_t}) \, |\partial^\ell_x K_{e^{-\xi e^u L_\nu/2}}(x)|\,\dmu(x)\,\du\,\dd\xi\\
&\lesssim_\nu \int_0^\infty |\Psi_t(\xi)| \int_\RR \frac{e^{\ell u/2}}{\xi^{\ell/2}}  \exp\left(-\frac{\cosh u}\xi \right)
( \cosh^{\gamma_t} u + \xi^{\gamma_t})\,\du\,\dd\xi.
\end{split}\]
The latter integral can be bounded much as in \cite[p.~185--186]{MaVa}, and we obtain
\begin{equation*}
	\big\Vert e^{\varepsilon |\cdot|_\dist/\sqrt{t}} \vfX_1^\ell K_{e^{-t\Delta_\nu}}\big\Vert_{L^1(G_\nu)}\lesssim_{\nu,\varepsilon} t^{-\ell/2}
\end{equation*}
for all $t \geq t_0$.
In a similar way, one can prove that
\begin{equation*}
	\left\Vert e^{\varepsilon |\cdot|_\dist/\sqrt{t}} \vfX_0 K_{e^{-t\Delta_\nu}}\right\Vert_{L^1(G_\nu)}
	\lesssim_{\nu, \varepsilon} t^{-1/2}
\end{equation*}
for all $t \geq t_0$ by following almost verbatim the calculations in \cite[p.~186--187]{MaVa}; in our case Lemma \ref{lm:3} plays the role of \cite[Proposition 3.2]{MaVa}.
\end{proof}

From Proposition \ref{prop:2} and Lemma \ref{lm:l1est_int_conv} we immediately deduce the following estimates for the integral kernel $\iheat_t$ of $e^{-t\Delta_\nu}$.

\begin{cor}\label{cor:l1hkb}
Let $\varepsilon\in[0,\infty)$. Then, for all $t > 0$,
\begin{align}
\label{eq:30}
\sup_{\bdy\in G_\nu} \int_{G_\nu} \iheat_t(\bdx,\bdy) \, e^{\varepsilon \dist(\bdx,\bdy)/\sqrt{t} } \,\dmub(\bdx)
&\lesssim_{\varepsilon,\nu} 1, \\
\label{eq:31}
\sup_{\bdy\in G_\nu} \int_{G_\nu} | \nabla_\nu^{\bdy} \iheat_t(\bdx,\bdy)| \, e^{\varepsilon \dist(\bdx,\bdy)/\sqrt{t} } \,\dmub(\bdx)
&\lesssim_{\varepsilon,\nu} t^{-1/2}, \\
\label{eq:32}
\sup_{\bdy\in G_\nu} \int_{G_\nu} | \nabla_\nu^{\bdx} \iheat_t(\bdx,\bdy)| \, e^{\varepsilon \dist(\bdx,\bdy)/\sqrt{t} } \,\dmub(\bdx)
&\lesssim_{\varepsilon,\nu} t^{-1/2},\\
\label{eq:33}
\sup_{\bdy\in G_\nu} \int_{G_\nu} |\nabla_\nu^{\bdy} \nabla_\nu^{\bdx} \iheat_t(\bdx,\bdy)| \, e^{\varepsilon \dist(\bdx,\bdy)/\sqrt{t} } \,\dmub(\bdx)
&\lesssim_{\varepsilon,\nu} t^{-1}.
\end{align}
\end{cor}
\begin{proof}
As $e^{-t\Delta_\nu}$ is self-adjoint, from \eqref{eq:inv_conv} it follows that $K_{e^{-t\Delta_\nu}}^* =K_{e^{-t\Delta_\nu}}$. Thus, the bounds \eqref{eq:30}-\eqref{eq:32} for $\iheat_t$ are immediate consequences, via Lemma \ref{lm:l1est_int_conv}, of the bounds for $K_{e^{-t\Delta_\nu}}$ in Proposition \ref{prop:2}. It only remains to prove \eqref{eq:33}.

On the other hand, by the semigroup property and the triangle inequality,
\begin{multline*}
\int_{G_\nu} |\nabla_\nu^{\bdy} \nabla_\nu^{\bdx} \iheat_t(\bdx,\bdy) | \, e^{\varepsilon\dist(\bdx,\bdy)/\sqrt{t} }  \,\dmub(\bdx) \\ 
\leq \int_{G_\nu} |\nabla_\nu^{\bdy} \iheat_{t/2}(\bdz,\bdy) | \, e^{\varepsilon\dist(\bdz,\bdy)/\sqrt{t} }
\int_{G_\nu} | \nabla_\nu^{\bdx} \iheat_{t/2} (\bdx,\bdz) | \, e^{\varepsilon\dist(\bdx,\bdz)/\sqrt{t} }  \,\dmub(\bdx) \,\dmub(\bdz) ,
\end{multline*}
and therefore the bounds \eqref{eq:32} and \eqref{eq:31} imply \eqref{eq:33}.
\end{proof}

\section{Multiplier theorem and Riesz transforms for \texorpdfstring{$p\leq 2$}{p <= 2}}\label{s:multriesz}

\subsection{Proof of the multiplier theorem}

The following proposition is a counterpart of \cite[Proposition 5.1]{MOV}.
For any $r > 0$, we write $\fctE_r$ for the set of even Schwartz functions on $\RR$ whose Fourier support is contained in $[-r,r]$.

\begin{prop}\label{prop:3}
	Let $F \in \fctE_r$ for some $r>0$. 
	Then,
	\begin{equation*}
\Vert K_{F(\sqrt{\Delta_\nu})}\Vert_{L^1(G_\nu)} 
\lesssim_\nu r^{[(\nu+1)/2,3/2]} \Vert K_{F(\sqrt{\Delta_\nu})} \Vert_{L^2(G_\nu)} .
	\end{equation*}
\end{prop}
\begin{proof}
Observe that since $F\in\fctE_r$ the finite propagation speed property for $\Delta_\nu$ \cite[Lemma~4.17]{MaPl} implies that $\supp K_{F(\sqrt{\Delta_\nu})}\subseteq \overline{B}_{G_\nu}(\bdzero,r)$.
Thus, for $r \leq 1$ the Cauchy--Schwarz inequality gives
\[
\Vert K_{F(\sqrt{\Delta_\nu})} \Vert_{L^1(G_\nu)}
\leq |\overline{B}_{G_\nu}(\bdzero,r)|^{1/2} \, \Vert K_{F(\sqrt{\Delta_\nu})} \Vert_{L^2(G_\nu)}
\simeq_\nu r^{(\nu+1)/2} \Vert K_{F(\sqrt{\Delta_\nu})} \Vert_{L^2(G_\nu)},
\]
where the last estimate follows from \cite[eq.\ (4.32)]{MaPl}.
	
On the other hand, for $r\geq 1$ we again use the Cauchy--Schwarz inequality to get
\begin{multline*}
	\Vert K_{F(\sqrt{\Delta_\nu})} \Vert_{L^1(G_\nu)}
	\leq \left(\int_{\overline{B}_{G_\nu}(\bdzero,r)} (1+x^\nu)^{-1}\,\dmu(x)\,\du \right)^{1/2}\\
	\times \left[ \Vert K_{F(\sqrt{\Delta_\nu})}\Vert_{L^2(G_\nu)} + \left(\int_{\overline{B}_{G_\nu}(\bdzero,r)} |K_{F(\sqrt{\Delta_\nu})}(x,u)|^2 \, x^\nu\,\dmu(x)\,\du \right)^{1/2} \right].
\end{multline*}
By \cite[Corollary 4.16]{MaPl},
\begin{equation*}
	\left(\int_{\overline{B}_{G_\nu}(\bdzero,r)} |K_{F(\sqrt{\Delta_\nu})}(x,u)|^2 \,x^\nu\,\dmu(x)\,\du \right)^{1/2}
	\lesssim_\nu r^{1/2} \, \Vert K_{F(\sqrt{\Delta_\nu})}\Vert_{L^2(G_\nu)}.
\end{equation*}
Finally,
\begin{multline*}
\int_{\overline{B}_{G_\nu}(\bdzero,r)} (1+x^\nu)^{-1} \,\dmu(x)\,\du \\
= \int_0^\infty \int_\RR \ind_{[0,r]}( \arccosh(\cosh u +x^2e^{-u}/2 )) \frac{x^{\nu-1}}{1+x^\nu}\,\dd x\, \du,
\end{multline*}
which is bounded by a constant multiple of $r^2$ (cf. \cite[p.~371]{MOV}). 

Combining the above gives the claim.	
\end{proof}

The following result is an analogue of \cite[Propositions~5.3 and 5.5]{MOV}, and gives us a weighted $L^1$-estimate for convolution kernels of operators $F(\Delta_\nu)$ corresponding to compactly supported multipliers $F$ which are sufficiently smooth.

\begin{prop}\label{prop:4}
Let $F\in L^2(\RR)$ be even and supported in $[-2,2]$. Then
\begin{equation}\label{eq:27}
\| K_{F(t\sqrt{\Delta_\nu})}  \, (1+|\cdot|_\dist/t)^\varepsilon \|_{L^1(G_\nu)}
\lesssim_{\nu,s_0,s_\infty,\varepsilon} \begin{cases}
\Vert F \Vert_{\sobolev{s_0}{2}} &\text{if } t \geq 1,\\
\Vert F \Vert_{\sobolev{s_\infty}{2}} + t^{(2-\nu)_+/2} \Vert F\Vert_{\sobolev{s_0}{2}} &\text{if } t \leq 1,
\end{cases}
\end{equation}
for all $\varepsilon \geq 0$, $t>0$, $s_0 > 3/2+\varepsilon$ and $s_\infty > (\nu+1)/2+\varepsilon$. 
Moreover,
\begin{equation}\label{eq:28}
\Vert |\nabla_\nu K_{F(t\sqrt{\Delta_\nu})}| \Vert_{L^1(G_\nu)} 
\lesssim_{\nu,s_0,s_\infty} t^{-1} \begin{cases}
\Vert F \Vert_{\sobolev{s_0}{2}} &\text{if } t \geq 1,\\
\Vert F \Vert_{\sobolev{s_\infty}{2}} + t^{(2-\nu)_+/2} \Vert F\Vert_{\sobolev{s_0}{2}} &\text{if } t \leq 1,
\end{cases}
\end{equation}
for all $t>0$, $s_0 > 3/2$ and $s_\infty > (\nu+1)/2$. 
\end{prop}

\begin{proof}
In order to justify \eqref{eq:27} we follow, with some modifications, the steps of the proof of \cite[Proposition~5.3]{MOV}, using Propositions \ref{prop:6} and \ref{prop:3} above in place of \cite[Corollary 4.6 and Proposition 5.1]{MOV}. Specifically, as in the proof of \cite[Proposition~5.3]{MOV}, one decomposes $F = \sum_{\ell \in \NN} F_\ell$ and obtains the estimates
\[
\|K_{F_\ell(t\sqrt{\Delta_\nu})} (1+|\cdot|_\dist/t)^\varepsilon \|_{L^1(G_\nu)}
 \lesssim_{\nu,\varepsilon,s} 2^{\ell(\varepsilon-s)} (2^\ell t)^{[(\nu+1)/2,3/2]} t^{[-(\nu+1)/2,-3/2]} \|F\|_{\sobolev{s}{2}}
\]
for all $\ell \in \NN$, $t>0$ and $\varepsilon,s \geq 0$. Differently from \cite{MOV}, here we sum the above estimates using different Sobolev exponents according to the value of $2^\ell t$, thus obtaining
\begin{multline*}
\|K_{F(t\sqrt{\Delta_\nu})} (1+|\cdot|_\dist/t)^\varepsilon\|_{L^1(G_\nu)} 
\lesssim_{\nu,\varepsilon,s,s'} t^{[(2-\nu)/2,0]} \|F\|_{\sobolev{s}{2}} \, \sum_{\ell \in \NN \tc 2^\ell t \geq 1} 2^{\ell(3/2+\varepsilon-s)} \\
+ t^{[0,(\nu-2)/2]} \|F\|_{\sobolev{s'}{2}}  \sum_{\ell \in \NN \tc 2^\ell t \leq 1} 2^{\ell((\nu+1)/2+\varepsilon-s')} .
\end{multline*}
When $t\geq 1$, the second sum is empty, and taking $s=s_0>3/2+\varepsilon$ leads to the first bound in \eqref{eq:27}. When $t \leq 1$, in the case $\nu \geq 2$ we take $s=s'=s_\infty > (\nu+1)/2+ \varepsilon$, so both sums converge and the first one compensates the factor $t^{(2-\nu)/2}$; 
instead, in the case $\nu < 2$, we take $s=s_0$ and $s'=s_\infty$. In both cases, we obtain the second bound in \eqref{eq:27}.

As for \eqref{eq:28}, let $H(\lambda)=F(\lambda)e^{\lambda^2}$. Thus
$K_{F(t\sqrt{\Delta_\nu})} = K_{H(t\sqrt{\Delta_\nu})} \diamond_\nu K_{e^{-t^2 \Delta_\nu}}$ and,
by \eqref{eq:diamond} and \eqref{eq:17},
\begin{equation*}
	|\nabla_\nu K_{F(t\sqrt{\Delta_\nu})}| \leq  |K_{H(t\sqrt{\Delta_\nu})}|  \diamond_\nu |\nabla_\nu K_{e^{-t^2\Delta_\nu}}|.
\end{equation*}
Hence, Young's inequality \eqref{eq:young} and Proposition \ref{prop:2} give
\begin{equation*}
\|| \nabla_\nu K_{F(t\sqrt\Delta_\nu)}| \|_{L^1(G_\nu)} \leq  t^{-1} \|K_{H(t\sqrt{\Delta_\nu})}\|_{L^1(G_\nu)}  
\end{equation*}
and \eqref{eq:28} follows by applying \eqref{eq:27} with $\varepsilon=0$ and $H$ in place of $F$.
\end{proof}

\begin{proof}[Proof of Theorem \ref{thm:main_mult}]
Much as in the proof of \cite[Theorem~1.1]{MOV}, we may assume that $F$ is real, and decompose $F = \sum_{j \in \ZZ} F_j(2^{-j} \cdot)$, where $F_j = F(2^j \cdot) \psi$ for a suitable cutoff $\psi$. Choose $\varepsilon>0$ such that $s_0 > 3/2+\varepsilon$ and $s_\infty > (\nu+1)/2 + \varepsilon$.
By applying Proposition \ref{prop:4} to the functions $F_j$ with $t = 2^{-j}$, we see that the convolution kernels $K_j$ of the operators $F_j(2^{-j} \sqrt{\Delta_\nu})$ satisfy the assumptions \eqref{eq:65} of Corollary \ref{cor:4} with $c = 2$, provided
\[
\sup_{j \leq 0} \|F_j\|_{\sobolev{s_0}{2}} < \infty, \qquad \sup_{j > 0} {(\|F_j\|_{\sobolev{s_\infty}{2}} + 2^{-j(2-\nu)_+/2} \|F_j\|_{\sobolev{s_0}{2}})} < \infty.
\]
This condition is clearly satisfied under the assumption \eqref{eq:ass_mult}, thus Corollary \ref{cor:4} and duality yield the required boundedness properties of $F(\sqrt{\Delta})$.
\end{proof}

\subsection{Sharpness of the smoothness requirement}\label{ss:sharpness}

Theorem \ref{thm:main_mult} implies that, under the scale-invariant smoothness assumption
\[
\|F\|_{L^\infty_{s,\sloc}} \defeq \sup_{t > 0} \|F(t\cdot) \chi \|_{\sobolev{s}{\infty}} < \infty
\]
for some $s > \max\{3,\nu+1\}/2$, the operator $F(\Delta_\nu)$ is of weak type $(1,1)$, bounded from $H^1(G_\nu)$ to $L^1(G_\nu)$, and bounded on $L^p(G_\nu)$ for all $p \in (1,\infty)$. We shall now show that the threshold $\max\{3,\nu+1\}/2$ cannot be replaced by a smaller quantity.

In the case $\nu \geq 2$, we have $\max\{3,\nu+1\}/2 = (\nu+1)/2$. So, in this case, the claimed sharpness is a consequence of the following result, which also gives a $p$-dependent lower bound on the smoothness requirement for $L^p(G_\nu)$-boundedness.

\begin{prop}\label{prop:sharp1}
If $\nu\geq 1$, then, for all $p \in (1,\infty)$,
\begin{align*}
\inf \left\{s>0 \tc \Vert F(\Delta_\nu) \Vert_{L^1 \to L^{1,\infty} } \lesssim_s \|F\|_{L^\infty_{s,\sloc}} \ \forall F \in C^\infty_c(\Rpos) \right\}
&\geq \frac{\nu+1}{2},\\
\inf \left\{s>0 \tc \Vert F(\Delta_\nu) \Vert_{L^p \to L^{p}} \lesssim_s \|F\|_{L^\infty_{s,\sloc}} \ \forall F \in C^\infty_c(\Rpos) \right\}
&\geq (\nu+1)\left|\frac{1}{2}-\frac{1}{p}\right|.
\end{align*} 
\end{prop}

The proof of the above result is based on a transplantation argument, analogous to that used, e.g., in \cite{KST,Ma17}. The idea is that, by an appropriate ``contraction'' procedure via dilations, the semidirect product $G_\nu = X_\nu \rtimes \RR$ can be approximated by the direct product $X_\nu \times \RR$, and the natural Laplacian $-\partial_u^2 + L_\nu$ on the direct product behaves like $L_{\nu+1}$ when restricted to a suitable class of ``radial functions''.

\begin{lm}\label{lm:transplantation}
For all $F \in C_0(\Rnon)$ and $p \in [1,\infty]$,
\begin{align*}
\|F(L_{\nu+1})\|_{L^1(X_{\nu+1}) \to L^{1,\infty}(X_{\nu+1})} &\leq \liminf_{r \to 0^+} \|F(r\Delta_\nu) \|_{L^1(G_\nu) \to L^{1,\infty}(G_\nu)},\\
\|F(L_{\nu+1})\|_{L^p(X_{\nu+1}) \to L^p(X_{\nu+1})} &\leq \liminf_{r \to 0^+} \|F(r\Delta_\nu) \|_{L^p(G_\nu) \to L^p(G_\nu)}.
\end{align*}
\end{lm}
\begin{proof}
Let $\delta_\lambda : G_\nu \to G_\nu$ denote the isotropic dilation of parameter $\lambda>0$, given by $\delta_\lambda(x,u) = (\lambda x,\lambda u)$. Then from \eqref{eq:Deltanu} we see that
\begin{equation*}
(\Delta_\nu (f\circ\delta_\lambda))\circ\delta_\lambda^{-1} = \lambda^2 \Delta_\nu^{(\lambda)}f,
\end{equation*}
where 
\begin{equation*}
\Delta_\nu^{(\lambda)}= -\partial_u^2 +e^{2u/\lambda} L_\nu.
\end{equation*}
As a consequence, if $D_\lambda f = f \circ \delta_\lambda$, then $D_\lambda$ is a multiple of an $L^2(G_\nu)$-isometry and, for any bounded Borel function $F$,
\begin{equation*}
F(\Delta_\nu^{(\lambda)}) = D_\lambda^{-1} F(\lambda^{-2}\Delta_\nu) D_\lambda;
\end{equation*}
as $D_\lambda$ also preserves $L^p$ and $L^{1,\infty}$ norms (up to multiplicative constants), we conclude that
\begin{align*}
\|F(\Delta_\nu^{(\lambda)})\|_{L^1 \to L^{1,\infty}} &= \|F(\lambda^{-2}\Delta_\nu) \|_{L^1 \to L^{1,\infty}}, \\
\quad \|F(\Delta_\nu^{(\lambda)})\|_{L^p \to L^{p}} &= \|F(\lambda^{-2}\Delta_\nu) \|_{L^p \to L^{p}}.
\end{align*}
Moreover, if we set
\begin{equation*}
\Delta_\nu^{(\infty)}\defeq -\partial_u^2+L_\nu,
\end{equation*}
then we see that, for all $f \in \Sz_\even(\Rnon)\otimes C_c^\infty(\RR)$,
\[
\Delta_\nu^{(\lambda)} f \to \Delta_\nu^{(\infty)} f \quad\text{ as } \lambda \to \infty
\]
in $L^2(G_\nu)$. As $\Sz_\even(\Rnon)\otimes C_c^\infty(\RR)$ is a core for $\Delta_\nu^{(\infty)}$, following the proof of \cite[Theorem~5.2]{Ma17} we conclude that, for any $F \in C_0(\Rnon)$, 
\[
F(\Delta_\nu^{(\lambda)}) \to F(\Delta_\nu^{(\infty)}) \quad \text{as } \lambda \to \infty
\]
in the strong operator topology on $L^2(G_\nu)$, and therefore
\begin{align*}
\|F(\Delta_\nu^{(\infty)})\|_{L^1 \to L^{1,\infty}} &\leq \liminf_{\lambda \to \infty} \|F(\lambda^{-2} \Delta_\nu) \|_{L^1 \to L^{1,\infty}},\\
\|F(\Delta_\nu^{(\infty)})\|_{L^p \to L^p} &\leq \liminf_{\lambda \to \infty} \|F(\lambda^{-2} \Delta_\nu) \|_{L^p \to L^p}.
\end{align*}

Notice now that $\Phi : L^2(X_{\nu+1}) \to L^2(G_\nu)$ given by $\Phi g(x,u) = g(\sqrt{x^2+u^2})$ is (up to a multiplicative constant) a linear isometric embedding, and an elementary computation shows that
\[
\Delta_\nu^{(\infty)} \Phi g = \Phi L_{\nu+1} g
\]
for all $g \in \Sz_\even(\Rnon)$. We then conclude that $\Phi F(L_{\nu+1}) = F(\Delta_\nu^{(\infty)}) \Phi$
for all bounded Borel functions $F$, whence also
\begin{align*}
\|F(L_{\nu+1})\|_{L^1 \to L^{1,\infty}} &\leq \|F(\Delta_\nu^{(\infty)})\|_{L^1 \to L^{1,\infty}} \\
\|F(L_{\nu+1})\|_{L^p \to L^p} &\leq \|F(\Delta_\nu^{(\infty)})\|_{L^p \to L^p}.
\end{align*}
The conclusion follows by combining the above bounds.
\end{proof}

\begin{proof}[Proof of Proposition \ref{prop:sharp1}]
We only discuss weak type $(1,1)$ bounds; a similar argument applies to $L^p$-bounds.

Let us assume a contrario that a bound of the form
\begin{equation*}
\Vert F(\Delta_\nu) \Vert_{L^1(G_\nu)\to L^{1,\infty}(G_\nu)} \lesssim \sup_{r>0} \Vert F(r\cdot)\chi\Vert_{L^\infty_s(\RR)}
\end{equation*}
holds for some $s < (\nu+1)/2$. Combining this with Lemma \ref{lm:transplantation} would then give the bound
\begin{equation*}
\Vert F(L_{\nu+1}) \Vert_{L^1(X_{\nu+1}) \to L^{1,\infty}(X_{\nu+1})} \lesssim \sup_{r>0} \Vert F(r\cdot)\chi\Vert_{L^\infty_s(\RR)}.
\end{equation*}
By the discussion in \cite[Section 1.3]{KanPre}, however, we know that the latter bound cannot hold when $s < (\nu+1)/2$.
\end{proof}

We now turn to the discussion of the case where $\nu < 2$. In this case, $\max\{3,\nu+1\}/2 = 3/2$, and the sharpness of this threshold is given by the following result. 

\begin{prop}\label{prop:sharp2}
If $\nu\geq 1$, then 
\begin{equation*}
	\inf \left\{s > 0 \tc \Vert F(\Delta_\nu)\Vert_{H^1(G_\nu)\to L^{1}(G_\nu)}
	\lesssim_s \|F\|_{L^\infty_{s,\sloc}} \ \forall F \in C^\infty_c(\Rpos) \right\}\geq \frac{3}{2}.
\end{equation*} 
\end{prop}

A discrete analogue of this result in the setting of flow trees is contained in \cite[Proposition 6.13 and Remark 6.14]{MSTV}; as in \cite{MSTV}, the proof here exploits ideas from \cite{MM16} in the construction of appropriate functions $F$ to test the above bounds.

We shall need a couple of auxiliary lemmas.

\begin{lm}\label{lm:19}
For all $n \in \Npos$, we can write
\begin{equation}\label{eq:68}
	\left( \frac{-1}{\sinh r} \partial_r\right)^n e^{ir\xi} 
	= \frac{-i\xi e^{ir\xi}}{\sinh^{n} r} P_{n-1}(\xi,\coth r) \qquad \forall \xi\in\RR,\ r>0,
\end{equation}
where $P_k(x,y)$ is a polynomial of degree $k$ such that $P_k(0,1)=k!$ for any $k \in \NN$.
\end{lm}
\begin{proof}
An easy induction argument shows that \eqref{eq:68} holds 
with the polynomials $P_k$
given by the recursive formula
\[
P_k(x,y) = \begin{cases}
1 &\text{if } n=0,\\
(ky-ix)P_{k-1}(x,y)+(y^2-1)\partial_y P_{k-1}(x,y) &\text{if } k>0.
\end{cases}
\]
In particular,
\[
P_k(0,1) = \begin{cases} 
1 &\text{if } k=0,\\
k P_{k-1}(0,1) &\text{if } k>0,
\end{cases}
\]
i.e., $P_k(0,1) = k!$, as claimed.
\end{proof}

The next lemma shows how to produce elements of the Hardy space $H^1(G_\nu)$ by means of functions of $\Delta_\nu$. Recall that, for any $r > 0$, we write $\fctE_r$ for the set of even Schwartz functions on $\RR$ whose Fourier support is contained in $[-r,r]$.

\begin{lm}\label{lm:hardy}
For any $r>0$ and $\psi \in \fctE_r$,
\[
K_{\Delta_\nu \psi(\sqrt{\Delta_\nu})} \in H^1(G_\nu).
\]
\end{lm}
\begin{proof}
We shall show that $K_{\Delta_\nu \psi(\sqrt{\Delta_\nu})}$ is a multiple of an $H^1(G_\nu)$-atom (see Definition \ref{defin:HardyBMO}). Notice that, by Definition \ref{def:2}, any compact subset of $G_\nu$ is contained in an admissible set. Thus, it is enough to show that $K_{\Delta_\nu \psi(\sqrt{\Delta_\nu})}$ is a compactly supported $L^2(G_\nu)$-function with vanishing integral.

Notice that $\Delta_\nu \psi(\sqrt{\Delta_\nu}) = \tilde\psi(\sqrt{\Delta_\nu})$, where $\tilde\psi(\lambda) = \lambda^2 \psi(\lambda)$ is in $\fctE_r$ too, thus $K_{\Delta_\nu \psi(\sqrt{\Delta_\nu})} \in L^2(G_\nu)$ by Proposition \ref{prop:6} and $\supp K_{\Delta_\nu \psi(\sqrt{\Delta_\nu})} \subseteq \overline{B}_{G_\nu}(\bdzero,r)$ by finite propagation speed \cite[Lemma~4.17]{MaPl}. On the other hand,
\[
K_{\Delta_\nu \psi(\sqrt{\Delta_\nu})} = \Delta_\nu K_{\psi(\sqrt{\Delta_\nu})},
\]
where $K_{\psi(\sqrt{\Delta_\nu})}$ is also an $L^2(G_\nu)$-function supported in $\overline{B}_{G_\nu}(\bdzero,r)$. If we now take a compactly supported, real-valued cutoff $\chi \in \Sz_\even(\Rnon)\otimes C_c^\infty(\RR)$ such that $\chi\equiv 1$ on a neighbourhood of $\overline{B}_{G_\nu}(\bdzero,r)$, then integration by parts gives
\[
\int_{G_\nu} K_{\Delta_\nu \psi(\sqrt{\Delta_\nu})} \,\dmub = \int_{G_\nu} \chi \, \Delta_\nu  K_{\psi(\sqrt{\Delta_\nu})} \,\dmub = \int_{G_\nu} \Delta_\nu \chi \, K_{\psi(\sqrt{\Delta_\nu})} \,\dmub = 0,
\]
as $\Delta_\nu \chi$ vanishes on $\overline{B}_{G_\nu}(\bdzero,r)$.
\end{proof}

\begin{proof}[Proof of Proposition \ref{prop:sharp2}]
Let $\delta>0$ be a small parameter to be fixed later, and set $F_t(\lambda)=e^{i t\lambda}\chi_\delta(\lambda)$, where $\chi_\delta$ is a smooth cutoff function supported in $[\delta^2/8,2\delta^2]$ and equal to $1$ in $[\delta^2/4,\delta^2]$. Let moreover $\psi \in \fctE_1$ be such that $\psi(0) = 1$.
	
By Corollary \ref{cor:5} and Lemma \ref{lm:19}, for all $r>0$,
\[\begin{split}
&2\pi \Gamma(1-\{\nu/2 \}) H_{F_t(\Delta_\nu) \Delta_\nu \psi(\sqrt{\Delta_\nu})}(r)  \\
&= \int_r^\infty \frac{\sinh x}{(\cosh x-\cosh r)^{\{\frac{\nu}{2} \} }} \int_\RR \xi^2 \psi(\xi) F_t(\xi^2) \left(\frac{-1}{\sinh x} \partial_x \right)^{\lfloor \frac{\nu}{2}\rfloor+1 } e^{ix\xi}\, \dd x\\
&= -i \int_\RR \xi^3 \psi(\xi) \chi_\delta(\xi^2) e^{it\xi^2} \int_r^\infty \frac{  P_{\lfloor \frac{\nu}{2}\rfloor}(\xi,\coth x) \, e^{i x\xi}}{(\cosh x-\cosh r)^{\{\frac{\nu}{2} \} } \sinh^{\lfloor \frac{\nu}{2}\rfloor} x } \,\dd x\\
&=  -i\int_\RR \xi^3 \psi(\xi) \chi_\delta(\xi^2) e^{i(t\xi^2+r\xi)} \int_0^\infty \frac{  P_{\lfloor \frac{\nu}{2}\rfloor}(\xi,\coth (x+r)) \, e^{i x\xi}}{(\cosh(x+r)-\cosh r)^{\{\frac{\nu}{2} \} } \sinh^{\lfloor \frac{\nu}{2}\rfloor} (x+r) } \,\dd x.
\end{split}\]
Expressing hyperbolic functions in terms of exponentials then shows that
\begin{equation}\label{eq:oscillatory}
H_{F_t(\Delta_\nu) \Delta_\nu \psi(\sqrt{\Delta_\nu})}(r) = e^{-{\frac{\nu}{2}r}} I_\nu(r/t,e^{-2r},t),
\end{equation}
where
\begin{gather*}
I_\nu(h,\sigma,t) = \int_\RR e^{it\phi(h,\xi)} \xi^3 \chi_\delta(\xi^2)  A_\nu(\xi,\sigma)\,\dd\xi, \qquad
\phi(h,\xi) = \xi^2+h\xi,\\
A_\nu(\xi,\sigma)=\frac{-2^{\frac{\nu}{2}}i \psi(\xi)}{2\pi \Gamma(1-\{\frac{\nu}{2}\})} \int_0^\infty \frac{ P_{\lfloor \frac{\nu}{2} \rfloor} \left(\xi,1+\frac{2\sigma e^{-2x} }{1-\sigma e^{-2x}} \right) e^{i x\xi}}{(1-\sigma e^{-x})^{\{\frac{\nu}{2} \} } (1-\sigma e^{-2x})^{\lfloor \frac{\nu}{2} \rfloor}  (e^{x}-1)^{\{\frac{\nu}{2}\}} e^{\lfloor \frac{\nu}{2} \rfloor x}  } \,\dd x.
\end{gather*}
Notice that $A_\nu$ is smooth in a neighbourhood of $(0,0)$ and, by Lemma \ref{lm:19},
\begin{equation*}
	|A_\nu(0,0)|=\frac{2^{ \frac{\nu}{2}}\lfloor\frac{\nu}{2}\big\rfloor! }{2\pi \Gamma(1-\{\frac{\nu}{2} \})} \int_0^\infty (e^{x}-1)^{-\{\frac{\nu}{2} \} } e^{-{\lfloor \frac{\nu}{2}\rfloor}x}   \,\dd x >0.
\end{equation*}
Hence, there is $\varepsilon > 0$ such that $|A_\nu(\xi,\sigma)| \gtrsim 1$ for any $\xi,\sigma \in [-\varepsilon,\varepsilon]$.

The phase function $\phi(h,\cdot)$ in the oscillatory integral $I_\nu(h,\sigma,t)$ has a single critical point $\xi_c(h) = -h/2$, and moreover $\partial_\xi^2 \phi(h,\xi) = 2$.
By the method of stationary phase (see \cite[Theorem~7.7.6]{Ho'90}), up to taking a smaller $\varepsilon >0$, we obtain that, for any $\delta \in (0,\varepsilon]$,
\begin{equation}\label{eq:stphase}
I_\nu(h,\sigma,t) = -2^{-3} \sqrt{\pi} e^{i\pi/4} h^3 \chi_\delta(h^2/4) A_\nu(-h/2,\sigma) e^{-i th^2/4} t^{-1/2} + O(t^{-3/2})
\end{equation}
as $t \to \infty$, uniformly in $h,\sigma \in [-\varepsilon,\varepsilon]$.

Let us now fix $\delta = \varepsilon/2$. From \eqref{eq:oscillatory} and \eqref{eq:stphase} we deduce that
\[
|H_{F_t(\Delta_\nu) \Delta_\nu \psi(\sqrt{\Delta_\nu})}(r)| \gtrsim e^{-\nu r/2} t^{-1/2}
\]
for any sufficiently large $t$ and any $r \in [\delta t,2\delta t]$. Thus, from \eqref{eq:H_notat} and Lemma \ref{lm:13}, we obtain that, for any sufficiently large $t$,
\begin{multline*}
	\Vert K_{F_t(\Delta_\nu) \Delta_\nu \psi(\sqrt{\Delta_\nu})}\Vert_{L^1(G_\nu)}
	\gtrsim \int_1^\infty | H_{F_t(\Delta_\nu) \Delta_\nu \psi(\sqrt{\Delta_\nu})}(r)| \, r e^{\nu r/2}\,\dd r \\
	\gtrsim t^{-1/2}\int_{\delta t}^{2\delta t} r\,\dd r 
	\simeq t^{3/2}.
\end{multline*}
Thus, for any sufficiently large $t$,
\begin{multline}\label{eq:hardy_lowerbd}
t^{3/2} 
\lesssim \Vert K_{F_t(\Delta_\nu) \Delta_\nu \psi(\sqrt{\Delta_\nu})} \Vert_{L^1(G_\nu)} 
= \Vert F_t(\Delta_\nu) K_{\Delta_\nu \psi(\sqrt{\Delta_\nu}) }\Vert_{L^1(G_\nu)}\\
 \lesssim \Vert F_t(\Delta_\nu)\Vert_{H^1(G_\nu)\to L^1(G_\nu)} ,
\end{multline}
where we used that $\Vert K_{\tilde{\psi}(\sqrt{\Delta_\nu})}\Vert_{H^1(G_\nu)} < \infty$ by Lemma \ref{lm:hardy}.
Combining this with the easily verified fact 
\begin{equation}\label{eq:sobolev_upperbd}
\sup_{r> 0} \|\chi F_t(r\cdot) \|_{\sobolev{s}{\infty}(\RR)} \simeq_s \| F_{t} \|_{\sobolev{s}{\infty}(\RR)} \simeq_s t^s, \qquad t \geq 1,
\end{equation}
finishes the proof.
\end{proof}

\begin{rem}
From Proposition \ref{prop:4} we deduce the bound
\[
\sup_{r \geq 1} \|K_{F(r \Delta_\nu)}\|_{L^1(G_\nu)} \lesssim_s \|F\|_{\sobolev{s}{\infty}}
\]
for any $F$ supported in $[-4,4]$ and any $s>3/2$. The threshold $3/2$ here is sharp too, as can be seen from the proof of Proposition \ref{prop:sharp2}: indeed, this follows by comparing \eqref{eq:sobolev_upperbd} with the fact that, by \eqref{eq:hardy_lowerbd}, $\|K_{F_t(\Delta_\nu)}\|_{L^1(G_\nu)} \gtrsim t^{3/2}$ for sufficiently large $t$.
\end{rem}

\subsection{Boundedness of the Riesz transforms for \texorpdfstring{$p\leq 2$}{p <= 2}}
We shall now verify that the Riesz transforms $\Riesz_0$ and $\Riesz_1$ satisfy the assumptions of Theorem \ref{thm:3} and thus are of weak type $(1,1)$, bounded from $H^1(G_\nu)$ to $L^1(G_\nu)$, and bounded on $L^p(G_\nu)$, $p\in(1,2]$. Much as in other works in the literature (see, e.g., \cite{HeSt,MSTV,MaVa}), this verification is an immediate consequence of the gradient heat kernel bounds stated in Corollary \ref{cor:l1hkb}.

\begin{proof}[Proof of Theorem \ref{thm:main_riesz} for {$p\in(1,2]$}]
Let $k\in\{0,1\}$, and notice that
\[
\Riesz_k = \vfX_k \Delta_\nu^{-1/2} = \frac{1}{\sqrt{\pi}} \int_0^\infty \vfX_k e^{-t\Delta_\nu} \frac{\dd t}{\sqrt{t}}.
\]
Thus, we can decompose $\Riesz_k = \sum_{j\in\ZZ} T^{(k)}_j$,  where the $T^{(k)}_j$ are integral operators associated with the kernels
\begin{equation*}
	\iker^{(k)}_j(\bdx,\bdy) = \frac{1}{\sqrt{\pi}} \int_{2^{j-1}}^{2^j} \vfX_k^{\bdx} \iheat_t(\bdx,\bdy)\frac{\dd t}{\sqrt{t}}.
\end{equation*}
By using the heat kernel bounds of Corollary \ref{cor:l1hkb}, we can now check that the operators $T^{(k)}_j$ satisfy the assumptions of Theorem \ref{thm:3} and Corollary \ref{cor:2} with $c=1/\sqrt{2}$ and $\varepsilon=1$. Indeed, for any $\bdy\in G_\nu$, by  \eqref{eq:32},
\begin{multline*}
\int_{G_\nu} | \iker^{(k)}_j(\bdx,\bdy)| \, (1+2^{-{j/2}}\dist(\bdx,\bdy)) \,\dmub(\bdx) \\
	\lesssim \int_{2^{j-1}}^{2^j} \int_{G_\nu} |\nabla_\nu^{\bdx} \iheat_t(\bdx,\bdy)| \, e^{\dist(\bdx,\bdy)/\sqrt{t}}\,\dmub(\bdx)\,\frac{\dd t}{\sqrt{t}} \lesssim 1.
\end{multline*}
Moreover, by \eqref{eq:33},
\begin{equation*}
	\sup_{\bdy\in G_\nu} \Vert |\nabla_\nu^\bdy\iker^{(k)}_j(\cdot,\bdy) | \Vert_{L^1(G_\nu)}
	\lesssim \sup_{\bdy\in G_\nu} \int_{2^{j-1}}^{2^j} \int_{G_\nu} \big| \nabla^{\bdy}_\nu \nabla_\nu^{\bdx}\iheat_t(\bdx,\bdy)\big|\, \dmub(\bdx)\,\frac{\dd t}{\sqrt{t}} \lesssim 2^{-j/2}.
\end{equation*}
By Theorem \ref{thm:3} we then conclude that $\Riesz_k$ is of weak type $(1,1)$, bounded from $H^1(G_\nu)$ to $L^1(G_\nu)$, and bounded on $L^p(G_\nu)$ for all $p\in(1,2]$.
\end{proof}

\section{Convolution kernels of the Riesz transforms}\label{s:convRiesz}
	
\subsection{The kernels \texorpdfstring{$K_{\Riesz_0}$}{KR0} and \texorpdfstring{$K_{\Riesz_1}$}{KR1}}\label{ss:ConvKer}
	
Recall the definition \eqref{eq:def_Riesz}
of the Riesz transforms. By the well-known subordination formula, we can write
\[
	\Delta_\nu^{-1/2} = \frac{1}{\sqrt{\pi}} \int_0^\infty e^{-t\Delta_\nu}\frac{\dd t}{\sqrt{t}}.
\]
We now define
\begin{equation}\label{eq:21}
	K_{\Delta_\nu^{-1/2}} \defeq \frac1{\sqrt{\pi}} \int_0^\infty K_{e^{-t\Delta_\nu}} \frac{\dd t}{\sqrt{t}}.
\end{equation}
Notice that $K_{e^{-t\Delta_\nu}}$ is nonnegative on $G_\nu$, so $K_{\Delta_\nu^{-1/2}}$ is too. Moreover, by applying the $L^1$- and $L^2$-bounds from Lemma \ref{lm:2} for $t$ small and $t$ large respectively, one sees that $K_{\Delta_\nu^{-1/2}} \in L^1(G_\nu) + L^2(G_\nu)$ and that, at least for $f\in C_c(G_\nu)$,
\begin{equation}\label{eq:KDeltanusqrt}
\Delta_\nu^{-1/2} f(\bdx) =f\diamond_\nu  K_{\Delta_\nu^{-1/2}}(\bdx) = \int_{G_\nu} f(\bdy^-) \, \ell_{\bdy} K_{\Delta_\nu^{-1/2}}(\bdx) \,\dmub(\bdy)
\end{equation}
by \eqref{eq:transl} and \eqref{eq:diamond}.
As we shall see in Proposition \ref{prop:7} below, $K_{\Delta_\nu^{-1/2}}$ is actually a smooth function away from $\bdzero$.

As the differential operator $\vfX_0$ commutes with left translations on $G_\nu$
(see \eqref{eq:X0_li}), from \eqref{eq:KDeltanusqrt} we see that, at least formally,
\[
\Riesz_0 f = f\diamond_\nu K_{\Riesz_0}, \qquad K_{\Riesz_0} \defeq \vfX_0 K_{\Delta_\nu^{-1/2}},
\]
and moreover, by \eqref{eq:inv_conv},
\[
\Riesz_0^* f = f\diamond_\nu K_{\Riesz_0^*}, \qquad K_{\Riesz_0^*} \defeq K_{\Riesz_0}^*.
\]
The previous expressions for $\Riesz_0$ and $\Riesz_0^*$ must actually be interpreted with some care, because, as we shall see, the above-defined kernels $K_{\Riesz_0},K_{\Riesz_0^*} \in C^\infty_c(G_\nu \setminus \{\bdzero\})$ are not locally integrable in a neighbourhood of $\bdzero$; this is to be expected, due to the nature of $\Riesz_0$ and $\Riesz_0^*$ as a singular integral operators. Nevertheless, the kernels $K_{\Riesz_0},K_{\Riesz_0^*}$ may still be used as ``off-diagonal kernels'' for the corresponding operators, in the sense of Proposition \ref{prop:SIO_diamond}, as we shall discuss in greater detail in Section \ref{ss:decomposition} below.

The same approach cannot be directly applied to $\Riesz_1$, as $\vfX_1$ does not commute with left translations on $G_\nu$. A similar representation for $\Riesz_1$ can however be obtained by exploiting the relations between $G_\nu$ and $G_\nu^\rD$ discussed in Section \ref{ss:dunkl}.
Indeed, by \eqref{eq:KDeltanusqrt} and \eqref{eq:DHconv_even_lift} we can also write
\begin{equation}\label{eq:KDeltanusqrt_D}
\begin{split}
(\ext_\even \Delta_\nu^{-1/2} f)(\bdx) &= (\ext_\even f) \diamond_\nu^\rD (\ext_\even K_{\Delta_\nu^{-1/2}})(\bdx) \\
&= \int_{G_\nu^\rD} (\ext_\even f)(\bdy^-) \, \ell_{\bdy}^\rD (\ext_\even K_{\Delta_\nu^{-1/2}})(\bdx) \,\dmub^\rD(\bdy).
\end{split}
\end{equation}
Moreover, if $\vfX_1^\rD \defeq e^u D_\nu$ is the lifting to $G_\nu^\rD$ of the Dunkl operator $D_\nu$ on $X_\nu^\rD$, then we know from \eqref{eq:XjD_li} that $\vfX_1^\rD$ commutes with left translations on $G_\nu^\rD$,
and moreover $\vfX_1^\rD$ coincides with $\vfX_1$ when applied to even functions, in the sense that
\begin{equation}\label{eq:X1DX1}
\ext_\odd \vfX_1 f = \vfX_1^\rD \ext_\even f.
\end{equation}
Thus, if we write $\Riesz_1 = \ext_\odd^* \ext_\odd \Riesz_1 = \ext_\odd^* \vfX_1^\rD \ext_\even \Delta_\nu^{-1/2}$, where $\ext_\odd^*$ is the adjoint operator to $\ext_\odd$, then from \eqref{eq:KDeltanusqrt_D} we deduce that, at least formally,
\[
\Riesz_1 f = \ext_\odd^*((\ext_\even f)\diamond_\nu^\rD (\ext_\odd K_{\Riesz_1})) = ((\ext_\even f)\diamond_\nu^\rD (\ext_\odd K_{\Riesz_1}))|_{G_\nu}, \quad K_{\Riesz_1} \defeq \vfX_1 K_{\Delta_\nu^{-1/2}},
\]
and moreover, by \eqref{eq:inv_conv_D} and the fact that $(\ext_\odd K_{\Riesz_1})^\bullet = - \ext_\odd K_{\Riesz_1}^*$,
\[
\Riesz_1^* f = \ext_\even^*((\ext_\odd f)\diamond_\nu^\rD (\ext_\odd K_{\Riesz_1^*})) = ((\ext_\odd f)\diamond_\nu^\rD (\ext_\odd K_{\Riesz_1^*}))|_{G_\nu}, \qquad K_{\Riesz_1^*} \defeq -K_{\Riesz_1}^*.
\]
Again, the above identities must be interpreted with some care, as the kernels $\ext_\odd K_{\Riesz_1},\ext_\odd K_{\Riesz_1^*} \in C^\infty_c(G_\nu^\rD \setminus \{\bdzero\})$ are singular at $\bdzero$, but we postpone to Section \ref{ss:decomposition} a more detailed discussion, based on Proposition \ref{prop:SIO_diamond_D}.

We shall now obtain precise asymptotics at the origin and at infinity for the kernel $K_{\Delta_\nu^{-1/2}}$ and its derivatives $K_{\Riesz_j} = \vfX_j K_{\Delta_\nu^{-1/2}}$ on $G_\nu$, which play a crucial role in deriving $L^p$-boundedness properties for the Riesz transforms for $p \geq 2$.

\subsection{Estimates for \texorpdfstring{$K_{\Delta_\nu^{-1/2}}$}{KDeltanu-12}}
	
By \eqref{eq:21} and \eqref{eq:H_notat}, we can also write
\begin{equation}\label{eq:rieszkernelsradial}
K_{\Delta_\nu^{-1/2}}(\bdx) = \frac{2^{(2-\nu)/2}}{\Gamma(\nu/2)} \modu_\nu^{1/2}(\bdx) H_{\Delta_\nu^{-1/2}}(|\bdx|_\dist),
\end{equation}
where
\begin{equation}\label{eq:HsqrtDelta}
H_{\Delta_\nu^{-1/2}} = \frac1{\sqrt{\pi}} \int_0^\infty H_{e^{-t\Delta_\nu}} \frac{\dd t}{\sqrt{t}}.
\end{equation}

\begin{prop}\label{prop:7}
For all $\nu \geq 1$, the function $H_{\Delta_\nu^{-1/2}}$ is smooth on $\Rpos$. Moreover, the following identities hold for all $r \in \Rpos$.
\begin{enumerate}[label=(\roman*)]
\item\label{en:HsqrtDelta_integer} For all $k \in \Npos$,
\[
H_{\Delta_{2k}^{-1/2}}(r) = \frac{1}{\pi} \left(\frac{-1}{\sinh r}\partial_r\right)^{k-1} \frac{1}{r\sinh r}.
\]
\item\label{en:HsqrtDelta_frac} For all $\nu \geq 1$ and $\ell \in \Npos$,
\begin{equation}\label{eq:23}
\left(\frac{-1}{\sinh r}\partial_r\right)^{\ell} H_{\Delta_\nu^{-1/2}}(r) 
=\frac{1}{\Gamma(1-\{\frac{\nu}{2}\})} 
 \int_r^\infty \frac{\sinh x}{(\cosh x-\cosh r)^{\{\frac\nu2 \}}} 
 H_{\Delta_{2 + 2 \lfloor\frac{\nu}{2}\rfloor + 2\ell}^{-1/2}}(x) \,\dd x.
\end{equation}
\end{enumerate}
In particular, $K_{\Delta_\nu^{-1/2}}$ is smooth on $G_\nu \setminus \{\bdzero\}$.
\end{prop}
\begin{proof}
Notice first that, by \eqref{eq:3}, for all $x \in \Rpos$,
\begin{equation*}
\int_0^\infty H_{e^{-t\Delta_2}}(x) \frac{\dd t}{\sqrt{t}} = \frac{1}{\sqrt{\pi}} \frac{1}{x\sinh x},
\end{equation*}
and clearly differentiation in $x$ and integration in $t$ can be exchanged in the above integral.
Thus, by \eqref{eq:HsqrtDelta} and Proposition \ref{prop:heat_dim}\ref{en:heat_dim_der}, for all $k \in \Npos$,
\[
H_{\Delta_{2k}^{-1/2}}(x) = \frac{1}{\sqrt{\pi}} \int_0^\infty \left(\frac{-1}{\sinh x}\partial_x\right)^{k-1} H_{e^{-t\Delta_2}}(x) \frac{\dd t}{\sqrt{t}}
= \frac{1}{\pi} \left(\frac{-1}{\sinh x}\partial_x\right)^{k-1} \frac{1}{x\sinh x},
\]
as claimed in part \ref{en:HsqrtDelta_integer}.

Similarly, by Proposition \ref{prop:heat_dim}\ref{en:heat_dim_der}-\ref{en:heat_dim_intder}, for all $\nu \geq 1$, $\ell \in \NN$ and $r \in \Rpos$,
\[\begin{split}
&\left(\frac{-1}{\sinh r}\partial_r\right)^{\ell} H_{e^{-t\Delta_\nu}}(r) \\
&= \frac{1}{\Gamma(1-\{\frac{\nu}{2}\})} \int_r^\infty \frac{\sinh x}{(\cosh x-\cosh r)^{\{\frac\nu2 \}}} \left(\frac{-1}{\sinh x}\partial_x\right)^{\lfloor \frac{\nu}{2}\rfloor+\ell} H_{e^{-t\Delta_2}}(x) \,\dd x\\
&= \frac{1}{\Gamma(1-\{\frac{\nu}{2}\})} \int_r^\infty \frac{\sinh x}{(\cosh x-\cosh r)^{\{\frac\nu2 \}}} H_{e^{-t\Delta_{2+2\lfloor\frac{\nu}{2}\rfloor+2\ell}}}(x) \,\dd x;
\end{split}\]
by \eqref{eq:HsqrtDelta}, integration of the previous identity with respect to $\dd t/\sqrt{t}$ leads to the formula of part \ref{en:HsqrtDelta_frac}.
\end{proof}

\begin{rem}
Proposition \ref{prop:7} agrees with \cite[Proposition~2.2]{Ma23} for integer values of $\nu$, up to the constant $\tilde S_\nu$ in \eqref{eq:tS_nu}.
\end{rem}
	
Now we are ready to study the asymptotic behaviour of the derivatives of $H_{\Delta_\nu^{-1/2}}$. An analogue of the following result for integer $\nu$ is \cite[Proposition~2.3]{Ma23}.

\begin{prop}\label{prop:8}
For any $\ell\in\NN$ and $r\in \Rpos$,
\begin{equation}\label{eq:25}
\left( \frac{-1}{\sinh r} \partial_r \right)^\ell  H_{\Delta_\nu^{-1/2}}(r) = \begin{cases}
		\frac{ \Gamma(\ell+\frac{\nu}{2})}{\pi} \frac{(\cosh r)^{-\left(\ell+\frac\nu2\right)}}{ \log \cosh r} \left( 1+O\left(\frac{1}{\log \cosh r}\right)\right), &r\to \infty,\\
		\frac{2^{\ell+\frac\nu2}\Gamma(\ell+\frac{\nu}{2})}{2\pi}  r^{-(2\ell+\nu)}\left( 1+O(r^{2}+r^{\nu+2\ell})\right), & r\to 0^+.
\end{cases}
\end{equation}
\end{prop}
\begin{proof}
By \cite[Proposition~2.3]{Ma23} the claim holds for integer values of $\nu$. We shall now justify it for any noninteger $\nu\geq 1$.

Let us first prove the asymptotics for $r \to \infty$.	
We recall the generalised binomial formula:
\begin{equation}\label{eq:Newton}
	(x-y)^{-\alpha} = \sum_{k\in\NN} \binom{\alpha +k-1}{k} x^{-\alpha-k} y^k
\end{equation}
for all $\alpha>0$ and $x>y>0$, where
\begin{equation*}
	\binom{\beta}{k} = \frac{\beta(\beta-1)\ldots(\beta-k+1)}{k!},\qquad \beta\in\RR,\ k\in\NN.
\end{equation*}
Applying this to \eqref{eq:23} gives
\[\begin{split}
	&\left( \frac{-1}{\sinh r} \partial_r \right)^\ell H_{\Delta_\nu^{-1/2}}(r)\\
	&=  \frac{1}{\Gamma(1-\{\frac{\nu}{2}\})} \sum_{k\in\NN} \binom{\{\frac{\nu}{2}\} +k-1}{k}  (\cosh r)^k  \int_r^\infty \frac{\sinh x}{(\cosh x)^{\{\frac{\nu}{2}\}+k}} H_{\Delta_{2\lfloor \frac{\nu}{2}\rfloor+2\ell+2}^{-1/2}}(x) \,\dd x.
\end{split}\]
By using the known asymptotics for $H_{\Delta_{2\lfloor \frac{\nu}{2}\rfloor+2\ell+2}^{-1/2}}$ we arrive at
\begin{multline}\label{eq:24}
	\left( \frac{-1}{\sinh r} \partial_r \right)^\ell H_{\Delta_\nu^{-1/2}}(r)
	= \frac{\Gamma(\lfloor\frac{\nu}{2}\rfloor+\ell+1)}{\pi \Gamma(1-\{\frac{\nu}{2}\})} \sum_{k\in\NN} \binom{\{\frac{\nu}{2}\}+k-1}{k}  (\cosh r)^k \\
	\times \left(1+O\left(\frac{1}{\log\cosh r}\right)\right) \int_r^\infty \frac{\sinh x}{(\cosh x)^{\frac{\nu}{2}+k+\ell+1} \log \cosh x}  \,\dd x
\end{multline}
as $r \to \infty$,
where the implicit constant in the Big-O term is independent of $k$.

Observe that, for any $A>0$, by integration by parts,
\begin{equation*}
\frac{1}{A e^{A}}-\int_A^{\infty} \frac{\dd x}{x e^x} = \int_A^\infty \frac{\dd x}{x^2 e^x} \leq \int_A^\infty \frac{x+2}{x^3 e^x}\, \dd x = \frac1{A^2 e^A}.
\end{equation*}
Thus, by the substitution $\cosh x =\exp(y(k+\ell+\nu/2)^{-1})$,
\begin{multline*}
	\int_r^\infty \frac{\sinh x}{(\cosh x)^{\frac{\nu}{2}+k+\ell+1} \log \cosh x}  \,\dd x
=	\int_{(k+\ell+\nu/2)\log \cosh r}^\infty \frac{\dd y}{y \, e^y}\\
	= \frac{(\cosh r)^{-(k+\ell+\nu/2)}}{ (k+\ell+\nu/2) \log \cosh r}\left( 1 + O\left( \frac{1}{\log \cosh r}\right)\right)
\end{multline*}
as $r\to\infty$.
By plugging this into \eqref{eq:24} we obtain
\begin{multline}\label{eq:almostinftyasymp}
\left( \frac{-1}{\sinh r} \partial_r \right)^\ell H_{\Delta_\nu^{-1/2}}(r)\\
	= \frac{\Gamma(\lfloor\frac{\nu}{2}\rfloor+\ell+1)}{\pi\Gamma(1-\{\frac{\nu}{2}\})} \frac{(\cosh r)^{-\ell-\nu/2}}{\log \cosh r}\sum_{k\in\NN} \frac{\binom{\{\frac{\nu}{2}\} +k-1}{k}}{k+\ell+\nu/2}  \left( 1 + O\left( \frac{1}{\log \cosh r}\right)\right) .
\end{multline}
Now, notice that, by \eqref{eq:Newton}, for $t\in(0,1)$,
\[\begin{split}
	\sum_{k\in\NN} \binom{\{\frac{\nu}{2}\} +k-1}{k} \frac{t^{k+\ell+\nu/2}}{k+\ell+\nu/2} &= \int_0^t \sum_{k\in\NN} \binom{\{\frac{\nu}{2}\} +k-1}{k} s^{k+\ell-1+\nu/2}\, \dd s\\
	 &= \int_0^t s^{\ell-1+\nu/2} (1-s)^{-\{\frac{\nu}{2}\} }\,\dd s,
\end{split}\]
and taking the limit as $t\to 1^-$ yields
\begin{equation*}
	\sum_{k\in\NN} \binom{\{\frac{\nu}{2}\} +k-1}{k} \frac{ 1}{k+\ell+\nu/2}= \frac{\Gamma(\frac{\nu}{2}+\ell)\Gamma(1-\{\frac{\nu}{2}\})}{\Gamma(\lfloor\frac{\nu}{2}\rfloor + \ell+1)}.
\end{equation*}
Plugging this expression into \eqref{eq:almostinftyasymp} yields the desired asymptotics \eqref{eq:25} for $r \to \infty$.
	
Now we consider $r\in(0,1)$. In \eqref{eq:23} we split the interval of integration into $(r,1)$ and $ (1,\infty)$. In the latter case, by using \eqref{eq:25} for large arguments,
\begin{multline*}
	0 \leq \frac{1}{\Gamma(1-\{\frac{\nu}{2}\})} \int_1^\infty \frac{\sinh x}{(\cosh x-\cosh r)^{\{\frac\nu2 \}}}  H_{\Delta_{2 + 2 \lfloor\frac{\nu}{2}\rfloor + 2\ell}^{-1/2}}(x) \,\dd x \\
	\lesssim_{\nu,\ell} \int_1^\infty \frac{(\cosh x)^{-(\lfloor \frac{\nu}{2}\rfloor+\ell+1)} \sinh x}{(\cosh x-\cosh 1)^{\{\frac\nu2 \}}} \,\dd x \lesssim_{\nu,\ell} 1.
\end{multline*}

For the integration over the interval $(r,1)$, much as in the proof of \cite[Proposition~2.3]{Ma23}, we first notice that
\begin{equation*}
	\frac{\sinh x}{(\cosh x-\cosh r)^{\{\frac\nu2 \}}}= \frac{2^{\{\frac\nu2 \}}x}{(x^2-r^2)^{\{\frac\nu2 \}}}\big(1+O(x^2)\big),\qquad 0<r<x<1.
\end{equation*}
Thus, by applying also the known asymptotics $\eqref{eq:25}$ for $H_{\Delta_{2 + 2 \lfloor\frac{\nu}{2}\rfloor + 2\ell}}$, we obtain
\[\begin{split}
&\frac{1}{\Gamma(1-\{\frac\nu2\})} \int_r^1 \frac{\sinh x}{(\cosh x-\cosh r)^{\{\frac\nu2 \}}} H_{\Delta_{2 + 2 \lfloor\frac{\nu}{2}\rfloor + 2\ell}^{-1/2}}(x) \,\dd x \\
&= \frac{2^{1+\ell+\frac\nu2}\Gamma(\ell+\lfloor\frac\nu2\rfloor+1) }{2\pi\Gamma(1-\{\frac\nu2\})} \int_r^1 \frac{x^{-(2\lfloor \frac{\nu}{2}\rfloor+2\ell+1 )}}{(x^2- r^2)^{\{\frac\nu2 \}}} (1+O(x^2)) \,\dd x\\
&= \frac{2^{\ell+\frac\nu2} \Gamma(\ell+\lfloor\frac\nu2\rfloor+1)}{\pi \Gamma(1-\{\frac\nu2\}) r^{\nu+2l}} \int_{1}^{1/r}  \frac{x^{-(2\lfloor \frac{\nu}{2}\rfloor+2\ell+1)}}{(x^2-1)^{\{\frac\nu2 \}}} (1+O((rx)^2))\,\dd x\\
&=\frac{2^{\ell+\frac\nu2} \Gamma(\ell+\lfloor\frac\nu2\rfloor+1)}{\pi\Gamma(1-\{\frac\nu2\}) r^{\nu+2l}} \left(\int_{1}^{\infty}  \frac{x^{-(2\lfloor \frac{\nu}{2}\rfloor+2\ell+1)}}{(x^2-1)^{\{\frac\nu2 \}}} \,\dd x +E(r) \right),
\end{split}\]
where $E(r)$ is the error term given by
\begin{equation*}
	E(r)= \int_{1}^{1/r}  \frac{x^{-(2\lfloor \frac{\nu}{2}\rfloor+2\ell+1)}}{(x^2-1)^{\{\frac\nu2 \}}}O((xr)^2)\,\dd x -\int_{1/r}^{\infty}  \frac{x^{-(2\lfloor \frac{\nu}{2}\rfloor+2\ell+1)}}{(x^2-1)^{\{\frac\nu2 \}}} \,\dd x.
\end{equation*}
	
For the main term we have
\begin{equation*}
\int_{1}^{\infty} \frac{x^{-(2\lfloor \frac{\nu}{2}\rfloor+2\ell+1)}}{(x^2-1)^{\{\frac\nu2 \}}} \,\dd x 
= \frac{1}{2} \int_{0}^{1} (1-u)^{-\{\frac\nu2 \}} u^{\nu/2+\ell-1} \,\dd u
= \frac{\Gamma(1-\{\frac\nu2\})\Gamma(\frac\nu2+\ell)}{2\Gamma(\lfloor\frac\nu2\rfloor +\ell+1)},
\end{equation*}
whereas, for the error term,
\[
|E(r)|
\lesssim_{\nu,\ell} r^2 \int_{1}^{1/r} \frac{x^{-(2\lfloor \frac{\nu}{2}\rfloor+2\ell-1)}}{(x^2-1)^{\{\frac\nu2 \}}}\,\dd x 
+ \int_{1/r}^{\infty} \frac{x^{-(2\lfloor \frac{\nu}{2}\rfloor+2\ell+1)}}{(x^2-1)^{\{\frac\nu2 \}}}\,\dd x  
\simeq_{\nu,\ell} r^2 + r^{\nu+2\ell};
\]
notice that the first integral has a different asymptotic behaviour for $r \to 0^+$ according to whether $\nu+2\ell$ is less than or greater than $2$ (the case $\nu+2\ell = 2$ does not occur, as $\nu$ is not an integer).

Combining the above yields the asymptotics \eqref{eq:25} for $r \to 0^+$.
\end{proof}

\subsection{Estimates for \texorpdfstring{$K_{\Riesz_0}$}{KR0} and \texorpdfstring{$K_{\Riesz_1}$}{KR1}}
Recall from Section \ref{ss:ConvKer} that
\begin{equation}\label{eq:rieszkernels}
	K_{\Riesz_0} = \vfX_0 K_{\Delta_\nu^{-1/2}} ,\qquad K_{\Riesz_1} = \vfX_1 K_{\Delta_\nu^{-1/2}} 
\end{equation}
and
\begin{equation}\label{eq:adjrieszkernels}
	K_{\Riesz_0^*} = K_{\Riesz_0}^* ,\qquad K_{\Riesz_1^*} = - K_{\Riesz_1}^* .
\end{equation}

\begin{prop}\label{prop:9}
For all $\bdx = (x,u)\in G_\nu\setminus\{\bdzero\}$,
\begin{align*}
\frac{\Gamma(\nu/2)}{2^{(2-\nu)/2}} \modu_\nu^{-1/2}(\bdx) K_{\Riesz_1}(\bdx) 
&=\frac{x}{\sinh r} \partial_r H_{\Delta_\nu^{-1/2}}(r),\\ 
\frac{\Gamma(\nu/2)}{2^{(2-\nu)/2}} \modu_\nu^{-1/2}(\bdx) K_{\Riesz^\ast_1}(\bdx)
&=-\frac{e^{-u} x}{\sinh r} \partial_r H_{\Delta_\nu^{-1/2}}(r),
\end{align*}
and 
\begin{align*}
\frac{\Gamma(\nu/2)}{2^{(2-\nu)/2}} \modu_\nu^{-1/2}(\bdx) K_{\Riesz_0}(\bdx) &= -\frac{\nu}{2} H_{\Delta_\nu^{-1/2}}(r) +\frac{\sinh u-e^{-u} x^2/2}{\sinh r} \partial_r H_{\Delta_\nu^{-1/2}}(r),\\
\frac{\Gamma(\nu/2)}{2^{(2-\nu)/2}} \modu_\nu^{-1/2}(\bdx) K_{\Riesz^\ast_0}(\bdx) &= -\frac{\nu}{2} H_{\Delta_\nu^{-1/2}}(r) -\frac{\sinh u+e^{-u}x^2/2}{\sinh r} \partial_r H_{\Delta_\nu^{-1/2}}(r),
\end{align*}
where
\begin{equation}\label{eq:distance_origin}
	r=|(x,u)|_\dist = \arccosh (\cosh u + e^{-u}x^2/2).
\end{equation}
\end{prop}
\begin{proof}
The expression \eqref{eq:distance_origin} for $r$ follows from \eqref{eq:dist}.

Now, much as in \eqref{eq:der_modu}, by \eqref{eq:vfs}, \eqref{eq:modular} and \eqref{eq:distance_origin} we see that
\begin{equation}\label{eq:44}
\begin{aligned}
\vfX_0 \modu_\nu^{1/2}(\bdx) &= -\frac{\nu}{2} \modu_\nu^{1/2}(\bdx), & \vfX_0 r(\bdx) &= \frac{\sinh u - e^{-u} x^2/2}{\sinh r} ,\\
\vfX_1 \modu_\nu^{1/2}(\bdx) &= 0, & \vfX_1 r(\bdx) &= \frac{x}{\sinh r}.
\end{aligned}
\end{equation}
Thus, the claimed formulas for $K_{\Riesz_0}$ and $K_{\Riesz_1}$ follow by \eqref{eq:rieszkernels} and \eqref{eq:rieszkernelsradial}, the Chain Rule and the Leibniz Rule.

Now, by \eqref{eq:involution}, \eqref{eq:modular} and \eqref{eq:distance_origin},
\begin{equation*}
(\modu_\nu^{1/2} f)^\ast(\bdx) = \modu_\nu^{1/2}(\bdx) f(\bdx^-), \qquad r(\bdx^-)=r(\bdx),
\end{equation*}
and the formulas for the kernels $K_{\Riesz_j^*}$ follow by \eqref{eq:adjrieszkernels} from those proved for the kernels $K_{\Riesz_j}$.
\end{proof}

\begin{prop}\label{prop:10}
For all $\bdx = (x,u) \in G_\nu \setminus \{\bdzero\}$,
\begin{align*}
K_{\Riesz_0}(\bdx) &= -\frac{\nu}{\pi} \frac{u}{r^{\nu+2}} + E_0(\bdx),
& K_{\Riesz^\ast_0}(\bdx) &= \frac{\nu}{\pi} \frac{u}{r^{\nu+2}} + \widetilde{E}_0(\bdx),\\
K_{\Riesz_1}(\bdx) &= -\frac{\nu}{\pi} \frac{x}{r^{\nu+2}} + E_1(\bdx), 
& K_{\Riesz^\ast_1}(\bdx) &= \frac{\nu}{\pi} \frac{x}{r^{\nu+2}} + \widetilde{E}_1(\bdx),
\end{align*}	
where $r$ is as in \eqref{eq:distance_origin}, while $E_j,\widetilde{E}_j \in L^1_\loc(G_\nu)$ for $j=0,1$.
Analogous expressions hold where $r^{\nu+2}$ is replaced by $(x^2+u^2)^{\frac{\nu}{2}+1}$.
\end{prop}
\begin{proof}
Propositions \ref{prop:8} and \ref{prop:9} yield that all the kernels are locally integrable off the origin. Thus, we only focus on the case $r<1$. 

Much as in \cite[eq.\ (3.1)]{Ma23}, notice that
\begin{equation}\label{eq:41}
r^2=(x^2+u^2) (1+O(r)).
\end{equation}
In particular,
\begin{gather*}
x = O(r), \quad u = O(r),\\
e^u = 1+O(r), \quad \sinh u = u(1+O(r^2)) \qquad \modu_\nu^{1/2}(\bdx)=1+O(r).
\end{gather*}
Thus, by combining the formulas of Proposition \ref{prop:9} and the asymptotics of Proposition \ref{prop:8},
\begin{align*}
K_{\Riesz_0}(\bdx)
&= -\frac{\nu}{\pi} \frac{u}{r^{\nu+2}} + O(r^{-\nu}), & K_{\Riesz_0^*}(\bdx) &= \frac{\nu}{\pi} \frac{u}{r^{\nu+2}} + O(r^{-\nu}), \\
K_{\Riesz_1}(\bdx) 
&= -\frac{\nu}{\pi} \frac{x}{r^{\nu+2}} + O(r^{-\nu}), & K_{\Riesz_1^*}(\bdx) &= \frac{\nu}{\pi} \frac{x}{r^{\nu+2}} + O(r^{-\nu}).
\end{align*}
Notice that we can replace $r^{\nu+2}$ by $(x^2+u^2)^{\nu/2+1}$ in the above formulas, because
\[
r^{-(\nu+2)} = (x^2+u^2)^{-\frac{\nu}{2}-1}(1+O(r))
\]
by \eqref{eq:41}. By \cite[Proposition~4.15]{MaPl}, $r^{-\nu}$ is integrable at the origin, thus we obtain the claimed expressions.
\end{proof}

It will be convenient to consider the difference kernel $K_{\Riesz_0-\Riesz_0^*} \defeq K_{\Riesz_0} -K_{\Riesz_0^*}$, which serves as off-diagonal kernel for the skewsymmetric part of $\Riesz_0-\Riesz_0^*$.

\begin{prop}\label{prop:11}
For all $\bdx = (x,u) \in G_\nu \setminus \{\bdzero\}$,
\begin{align*}
K_{\Riesz^\ast_1}(\bdx)
&= -\frac{2\nu}{\pi} \frac{\ind_{\{u\leq -1\}}}{u} x\aux^\nu(x) + E^\infty_1(\bdx), \\
K_{\Riesz_0-\Riesz_0^\ast}(\bdx) 
&= -\frac{2\nu}{\pi} \left[ \frac{\ind_{\{|u|\geq 1\}}}{u} \aux^\nu(x)+ \frac{\ind_{\{u\geq 1\}}}{u}\big(\aux^\nu_{(e^u)}(x)- \aux^\nu(x)\big) \right]+ E^\infty_0(\bdx),
	\end{align*}
where 
the functions $E^\infty_j$, $j=0,1$, are integrable at infinity, while
\begin{equation}\label{eq:aux}
\aux^\nu(x) \defeq (1+x^2)^{-(1+\nu/2)},
\end{equation}
and $\aux^\nu_{(\lambda)}$ is defined as in \eqref{eq:dilation}.
\end{prop}
\begin{proof}
The proof is similar to the proof of \cite[Proposition~3.4]{Ma23}. We provide a sketch.

By Propositions \ref{prop:9} and \ref{prop:8},
\begin{equation*}
K_{\Riesz_1}(\bdx)= -\frac{\nu}{\pi 2^{\nu/2}} \frac{x \modu_\nu^{1/2}(\bdx)}{(\cosh r)^{1+\nu/2} \log \cosh r} (1+\BigO(r^{-1})),
\end{equation*} 
where the Big-O notation refers to the decay at infinity. By Lemma \ref{lm:13} we see that the Big-O term and the main term for $u\leq 1$ are integrable at infinity; on the other hand, as in \cite[eq.\ (3.7)]{Ma23}, for $u\geq 1$ we have
\begin{equation}\label{eq:14}
\frac{\modu_\nu^{1/2}(\bdx)}{(\cosh r)^{1+\nu/2} \log \cosh r} = \frac{2^{1+\nu/2} e^{-(1+\nu)u} }{u(1+e^{-2u}x^2)^{1+\nu/2}} \left[ 1+\BigO\left( \frac{1+\log(1+e^{-2u}x^2)}{u}\right)\right].
\end{equation}
Hence,
\begin{equation*}
K_{\Riesz_1}(\bdx)=  -\frac{2\nu}{\pi} \frac{e^{-(\nu+1) u}x  \ind_{\{u\geq 1\} }}{u(1+e^{-2u}x^2)^{1+\nu/2}} \left[ 1+\BigO\left( \frac{1+\log(1+e^{-2u}x^2)}{u}\right)\right]+ \IT,
\end{equation*}
where $\IT$ stands for an integrable term at infinity. Again, the Big-O summand is integrable, since
\begin{multline*}
\int_1^\infty  \int_0^\infty \frac{x (1+\log(1+e^{-2u}x^2))}{u^2 e^{(\nu+1) u}  (1+e^{-2u}x^2)^{1+\nu/2}} \,\dmu(x)\,\du \\
= \int_1^\infty \frac{\dd u}{u^2} \int_0^\infty \frac{x^\nu (1+\log(1+x^2))}{(1+x^2)^{1+\nu/2}} \,\dd x < \infty.
\end{multline*}
Thus,
\begin{equation*}
K_{\Riesz_1}(\bdx)=  -\frac{2\nu}{\pi} \frac{e^{-(\nu+1) u}x  \ind_{\{u\geq 1\} }}{u(1+e^{-2u}x^2)^{1+\nu/2}} + \IT,
\end{equation*}
and by taking adjoints, as in \eqref{eq:adjrieszkernels}, one obtains the claimed expression for $K_{\Riesz_1^\ast}$.

Now we study $K_{\Riesz_0-\Riesz_0^\ast}$. By Propositions \ref{prop:9} and \ref{prop:8},
\begin{equation*}
K_{\Riesz_0-\Riesz_0^\ast}(\bdx)
= -\frac{\nu}{\pi 2^{\nu/2}} \frac{ \modu_\nu^{1/2}(\bdx) \, (2\sinh u)}{(\cosh r)^{1+\nu/2} \log \cosh r} (1+\BigO(r^{-1})).
\end{equation*}
Much as above, by Lemma \ref{lm:13},
\begin{equation*}
K_{\Riesz_0-\Riesz_0^\ast}(\bdx)
= -\frac{\nu}{\pi 2^{\nu/2}} \frac{ \modu_\nu^{1/2}(\bdx) (2\sinh u) \ind_{\{|u|\geq 1 \} } }{(\cosh r)^{1+\nu/2} \log \cosh r} + \IT.
\end{equation*}
Firstly, we focus on the case $u\geq 1$. By \eqref{eq:14} and the fact that $2\sinh u=e^u (1+\BigO(e^{-2u}))$ we obtain
\begin{multline*}
K_{\Riesz_0-\Riesz_0^\ast}(\bdx) \ind_{\{u\geq 1\} }\\
=- \frac{2\nu}{\pi} \frac{ e^{-\nu u}  \ind_{\{u\geq 1 \} } }{u (1+e^{-2u}x^2)^{1+\nu/2}}\left[1+\BigO\left(\frac{1+\log(1+e^{-2u}x^2)}{u} \right)\right] + \IT.
\end{multline*}
As previously, the Big-O term is integrable. This gives 
\begin{equation*}
K_{\Riesz_0-\Riesz_0^\ast}(\bdx) \ind_{\{u\geq 1\} }= -\frac{2\nu}{\pi} \frac{ e^{-\nu u} \ind_{\{u\geq 1 \} } }{u (1+e^{-2u}x^2)^{1+\nu/2}} + \IT.
\end{equation*}
Now clearly $K_{\Riesz_0-\Riesz_0^\ast}^* = - K_{\Riesz_0-\Riesz_0^\ast}$ by \eqref{eq:adjrieszkernels}. Thus, by taking the adjoint,
\begin{equation*}
K_{\Riesz_0-\Riesz_0^\ast}(\bdx) \ind_{\{u\leq -1\} }= -\frac{2\nu}{\pi} \frac{  \ind_{\{u\leq -1 \} } }{u (1+x^2)^{1+\nu/2}} + \IT.
\end{equation*}
By summing the above we get
\begin{equation*}
K_{\Riesz_0-\Riesz_0^\ast}(\bdx)= -\frac{2\nu}{\pi} \left(\frac{  \ind_{\{u\leq -1 \} } }{u (1+x^2)^{1+\nu/2}}+ \frac{ e^{-\nu u} \ind_{\{u\geq 1 \} } }{u (1+e^{-2u}x^2)^{1+\nu/2}}\right) + \IT.
\end{equation*}
One can immediately see that the last expression is of the desired form.
\end{proof}

\subsection{Decomposition into local and global parts}\label{ss:decomposition}

Since we already know that $\Riesz_0$ is bounded on $L^p(G_\nu)$, $p\in(1,2]$, in order to prove the same for $\Riesz_0^\ast$ it suffices to justify the boundedness of $\Riesz_0-\Riesz_0^\ast$.

\begin{prop}
We can decompose
\begin{equation}\label{eq:R_0_decomp_kernel}
K_{\Riesz_0-\Riesz_0^\ast} = K_{0}^{(1)}+K_{0}^{(2)}+K_{0}^{(3)},
\end{equation}
where, with the notation $\bdx = (x,u)$ and $r=|\bdx|_\dist$ as in \eqref{eq:distance_origin},
\begin{equation}\label{eq:20}
\begin{aligned}
K_0^{(1)}(\bdx)&= -\frac{2\nu}{\pi} \frac{u}{r^{\nu+2}} \ind_{(0,1)}(r)+E_0(\bdx),\\
K_0^{(2)}(\bdx)&= -\frac{2\nu}{\pi} \frac{\ind_{\{u\geq1 \} }}{u} (\aux^\nu_{(e^u)}(x) -\aux^\nu(x)),\\
K_0^{(3)}(\bdx)&= -\frac{2\nu}{\pi} \frac{\ind_{\{|u|\geq1 \} }}{u}\aux^\nu(x) ,
\end{aligned}
\end{equation}
and moreover $E_0 \in L^1(G_\nu)$. 
Correspondingly, we can split
\begin{equation}\label{eq:R_0_decomp}
\Riesz_0-\Riesz_0^\ast = \tilde\Riesz_0^{(1)} + \tilde\Riesz_0^{(2)} + \tilde\Riesz_0^{(3)},
\end{equation}
where, for $j=2,3$ and all $f \in C_c(G_\nu)$,
\begin{equation}\label{eq:Rz0j}
\tilde\Riesz_0^{(j)} f = f \diamond_\nu K_{0}^{(j)}
\end{equation}
while the remaining operator $\tilde\Riesz_0^{(1)}$ satisfies an analogous relation in the sense of off-diagonal kernels:
\[
\tilde\Riesz_0^{(1)} f(\bdx) = \int_{G_\nu} \ell_{\bdx} f(\bdy^-) \, K_{0}^{(1)}(\bdy) \,\dmub(\bdy) \quad\text{for a.a.\ } \bdx \notin \supp f.
\]
\end{prop}
\begin{proof}
The decomposition \eqref{eq:R_0_decomp_kernel} readily follows by Propositions \ref{prop:10} and \ref{prop:11}.

Moreover, from \eqref{eq:KDeltanusqrt}, \eqref{eq:rieszkernels} and Proposition \ref{prop:SIO_diamond}, it follows that
\[
\Riesz_0 f(\bdx) = \vfX_0 (f \diamond_\nu K_{\Delta_\nu^{-1/2}})(\bdx) = \int_{G_\nu} \ell_{\bdx} f(\bdy^-) \, K_{\Riesz_0}(\bdy) \,\dmub(\bdy) \quad\forall \bdx \notin \supp f.
\]
Since $\langle g, \Riesz_0^* f\rangle_{L^2(G_\nu)} = \langle \Riesz_0 g,f\rangle_{L^2(G_\nu)}$,
from \eqref{eq:SIO_der_adj} and \eqref{eq:adjrieszkernels} we also deduce that
\[
\Riesz_0^* f(\bdx) = \int_{G_\nu} \ell_{\bdx} f(\bdy^-) \, K_{\Riesz_0^*}(\bdy) \,\dmub(\bdy) \quad\forall \bdx \notin \supp f
\]
and
\[
(\Riesz_0-\Riesz_0^*) f(\bdx) = \int_{G_\nu} \ell_{\bdx} f(\bdy^-) \, K_{\Riesz_0-\Riesz_0^*}(\bdy) \,\dmub(\bdy) \quad\forall \bdx \notin \supp f.
\]
Thus, if we define the operators $\tilde\Riesz_0^{(j)}$ by \eqref{eq:Rz0j} for $j=2,3$, and the remaining operator $\tilde\Riesz_0^{(1)}$ by difference, so that \eqref{eq:R_0_decomp} holds, then we deduce the claimed relation between $\tilde\Riesz_0^{(1)}$ and $K_{0}^{(1)}$.
\end{proof}

Thus, the boundedness of $\Riesz_0^\ast$ on $L^p(G_\nu)$, $p\in(1,2]$, effectively boils down to the $L^p(G_\nu)$-boundedness of $\tilde\Riesz_0^{(j)}$, $j=1,2,3$.

A similar analysis can be performed for $\Riesz_1^*$.

\begin{prop}
We can decompose
\begin{equation}\label{eq:R_1_decomp_kernel}
K_{\Riesz_1^\ast} = K_1^{(1)}+K_1^{(2)}, 
\end{equation}
where
\begin{equation}\label{eq:62}
\begin{aligned}
K_1^{(1)}(\bdx)&= \frac{\nu}{\pi} \frac{x}{r^{\nu+2}} \ind_{(0,1)}(r) +E_1(\bdx),\\
K_1^{(2)}(\bdx)&= -\frac{2\nu}{\pi} \frac{\ind_{\{u\leq -1 \} }}{u} x\aux^\nu(x),
\end{aligned}
\end{equation}
and $E_1 \in L^1(G_\nu)$. Correspondingly, we can split
\begin{equation}\label{eq:R_1_decomp}
\Riesz_1^\ast = \tilde\Riesz_1^{(1)} + \tilde\Riesz_1^{(2)},
\end{equation}
where, for all $f \in C_c(G_\nu)$,
\begin{equation}\label{eq:Rz1j}
\tilde\Riesz_1^{(2)} f = ((\ext_\odd f) \diamond_\nu^\rD (\ext_\odd K_1^{(2)}))|_{G_\nu},
\end{equation}
while the remaining operator $\tilde\Riesz_1^{(1)}$ satisfies an analogous relation in the sense of off-diagonal kernels:
\[
\tilde\Riesz_1^{(1)} f(\bdx) 
= \int_{G_\nu^\rD} \ell_{\bdx}^\rD (\ext_\odd f)(\bdy^{-1}) \, (\ext_\odd K_1^{(1)})(\bdy) \,\dmub^\rD(\bdy) 
\quad\text{for a.a.\ } \bdx \notin \supp f.
\]
\end{prop}
\begin{proof}
The decomposition \eqref{eq:R_1_decomp_kernel} readily follows by Propositions \ref{prop:10} and \ref{prop:11}.

Moreover, by \eqref{eq:KDeltanusqrt_D}, \eqref{eq:X1DX1}, \eqref{eq:rieszkernels} and Proposition \ref{prop:SIO_diamond_D},
\begin{multline*}
\ext_\odd \Riesz_1 f(\bdx) = \vfX_1^\rD((\ext_\even f) \diamond_\nu^\rD (\ext_\even K_{\Delta_\nu^{-1/2}}))(\bdx) \\
=  \int_{G_\nu^\rD} \ell_{\bdx}^\rD (\ext_\even f)(\bdy^-) \, K_{\Riesz_1}^\rD(\bdy) \,\dmub^\rD(\bdy) \quad\forall \bdx \in G_\nu^\rD \setminus \supp (\ext_\even f),
\end{multline*}
where $K_{\Riesz_1}^\rD = \vfX_1^\rD \ext_\even K_{\Delta_\nu^{-1/2}} = \ext_\odd K_{\Riesz_1}$; in particular,
\[
\Riesz_1 f(\bdx) = \int_{G_\nu^\rD} \ell_{\bdx}^\rD (\ext_\even f)(\bdy^-) \, K_{\Riesz_1}^\rD(\bdy) \,\dmub^\rD(\bdy) \quad\forall \bdx \in G_\nu \setminus \supp f.
\]

Using the fact that $\langle g, \Riesz_1^* f \rangle_{L^2(G_\nu)} = \langle \Riesz_1 g, f \rangle_{L^2(G_\nu)} = \langle \ext_\odd \Riesz_1 g, \ext_\odd f \rangle_{L^2(G_\nu^\rD)}$, from \eqref{eq:SIO_der_adj_D} and \eqref{eq:adjrieszkernels} we also deduce that
\[
\Riesz_1^* f(\bdx) = \int_{G_\nu^\rD} \ell_{\bdx}^\rD (\ext_\odd f)(\bdy^-) \, K_{\Riesz_1^*}^\rD(\bdy) \,\dmub^\rD(\bdy) \quad\forall \bdx \in G_\nu \setminus \supp f,
\]
where $K_{\Riesz_1^*}^\rD = (K_{\Riesz_1}^\rD)^\bullet = \ext_\odd K_{\Riesz_1^*}$.

Thus, if we define the operator $\tilde\Riesz_1^{(2)}$ by \eqref{eq:Rz1j} and the remaining operator $\tilde\Riesz_1^{(1)}$ by difference, so that \eqref{eq:R_1_decomp} holds, then we deduce the claimed relation between $\tilde\Riesz_1^{(1)}$ and $K_{1}^{(1)}$.
\end{proof}

Thus, the $L^p(G_\nu)$-boundedness of $\Riesz_1^\ast$ for $p \in (1,2]$ boils down to the analogous boundedness of $\tilde\Riesz_1^{(j)}$ for $j=1,2$.

\section{An operator-valued spectral multiplier theorem}\label{s:opval}

As discussed in the introduction, this section is devoted to the proof of a conditional result, which states, roughly speaking, that if a self-adjoint operator satisfies an $L^p$ spectral multiplier theorem of Mihlin--H\"ormander type, then under certain assumptions it also satisfies an operator-valued spectral multiplier theorem.

First, let us recall some definitions. We begin with the R-boundedness for families of operators (see for instance \cite[Definition~8.1.1(1) and Remark 8.1.2]{HNVW2}).

\begin{defin}\label{def:Rbound}
Let $B$ be a Banach space equipped with a norm $\Vert \cdot\Vert_B$. Consider a family of bounded operators $\cT \subseteq \LinBnd(B)$. We say that $\cT$ is \emph{R-bounded} if 
there exists $C \in \Rpos$ such that,
for all finite sequences $(T_n)_{n=1}^N$ in $\cT$ and $(x_n)_{n=1}^N$ in $B$, $N\in\NN$, there holds
\begin{equation*}
	\int_0^1 \left\Vert \sum_{n=1}^N \rade_n(t) T_n x_n \right\Vert_{B}\,\dd t \leq C \int_0^1\left\Vert \sum_{n=1}^N \rade_n(t) x_n \right\Vert_{B}\,\dd t.
\end{equation*}
Here $\rade_n$ are the Rademacher functions. The smallest possible constant $C$ in the above inequality shall be denoted by $\RBd_{B}(\cT)$.
\end{defin}
We remark that the R-boundedness of a family of operators implies their uniform boundedness (see \cite[Theorem~8.1.3]{HNVW2}).

We also recall the definition of $0$-bisectorial operators (see, e.g., \cite[Definition~10.1.1]{HNVW2} and \cite{DelKri23}).
\begin{defin}
Let $B$ be a Banach space.	A closed, densely defined operator $L$ on $B$ is \emph{$0$-bisectorial} if its spectrum lies in $\RR$ and for every $\omega\in(0,\pi/2)$ there holds
\begin{equation*}
\sup_{\zeta\in\CC \tc \left|\arg (\pm\zeta)\right|>\omega} \Vert (I-\zeta L)^{-1}\Vert_{\LinBnd(B)}<\infty,
\end{equation*}
where $\arg\zeta\in(-\pi,\pi]$ denotes the argument of $\zeta \in \CC \setminus \{0\}$.
\end{defin}

Let $\breve\chi \in C_c^\infty(\RR)$ be a nontrivial nonnegative even cutoff supported in $\Rnoz\defeq\RR\setminus\{0 \}$. For $s>1/2$ we define the {\it H\"ormander class} $\hormander{s}$ by
\begin{equation*}
\hormander{s} = \left\{F\in L^2_\loc(\RR) \tc \Vert F\Vert_{\hormander{s}} \defeq \sup_{t>0} \Vert F(t\cdot) \breve\chi \Vert_{\sobolev{s}{2}}<\infty  \right\}.
\end{equation*}
Notice that different choices of the cutoff $\breve\chi$ give rise to equivalent norms on $\hormander{s}$. Moreover, by the Sobolev embedding theorem, the elements of $\hormander{s}$ are continuous and bounded on $\Rnoz$. Notice further that the norm $\|\cdot\|_{\hormander{s}}$ coincides with the norm $\|\cdot\|_{L^2_{s,\sloc}}$ of \eqref{eq:ass_mult_MH} when applied to functions supported in $\Rnon$.

In this section we assume that $(X,\mu)$ is a $\sigma$-finite measure space, and $L$ is a self-adjoint operator on $L^2(X,\mu)$, which we abbreviate as $L^2(X)$. Moreover, we assume that, for some $p\in(1,\infty)$ and $s_L\geq 1/2$, the operator $L$ has a \emph{bounded $\hormander{s}$-calculus on $L^p(X)$} for all $s>s_L$: this means that, for all bounded Borel functions $F : \RR \to \CC$ such that $F\in\hormander{s}$, the operator $F(L)$, initially defined on $L^2(X)$ by the spectral theorem, extends to a bounded operator on $L^p(X)$ and 
\[
\Vert F(L)\Vert_{\LinBnd(L^p(X))}\lesssim_s \Vert F\Vert_{\hormander{s}}.
\]
This condition implies that $\ind_{\{0\}}(L) = 0$, i.e., $L$ is injective on $L^2(X)$, thus $F(L)$ does not depend on $F(0)$. Hence, it suffices to consider $F$ defined on $\Rnoz$.

In the following lemma we gather some properties of such an operator $L$.

\begin{lm}\label{lm:14}
Fix $p\in(1,\infty)$ and $s_L\geq 1/2$. Let $L$ be a self-adjoint operator on $L^2(X)$, which has a bounded $\hormander{s}$-calculus on $L^p(X)$ for all $s>s_L$. Then the following assertions hold.
\begin{enumerate}[label=(\roman*)]
\item\label{lm:14(1)} The operator $L$ is $0$-bisectorial on $L^p(X)$.
\item\label{lm:14(tensor)} If $(Y,\nu)$ is a $\sigma$-finite measure space, then $L \otimes \id$ is $0$-bisectorial on $L^p(X \times Y)$ and has a bounded $\hormander{s}$-calculus on $L^p(X \times Y)$ for all $s>s_L$.
\item\label{lm:14(2)} If $\breve\chi \in C^\infty_c(\Rnoz)$ satisfies $\supp \breve\chi \subseteq [-2,-1/2] \cup [1/2,2]$ and
\begin{equation}\label{eq:dyadicpart}
\sum_{m \in \ZZ} \breve\chi(2^m \xi) = 1 \qquad \forall \xi \in \Rnoz,
\end{equation}
then there holds
\[\begin{split}
\Vert f\Vert_{L^p(X)} 
&\simeq_{L,\breve\chi,p} \int_0^1 \left\Vert \sum_{m\in\ZZ} \rade_m(t) \breve\chi(2^m L) f\right\Vert_{L^p(X)}\,\dd t
\end{split}\]
for all $f\in L^p(X)$.
\item\label{lm:14(3)} If $s> s_L+|1/p-1/2|$, then
\begin{equation}\label{eq:RbdFC}
\RBd_{L^p(X)} \{ F(L) \tc F \in \hormander{s},\ \Vert F \Vert_{\hormander{s}}\leq 1 \} < \infty.
\end{equation}
We then say that \emph{$L$ has an R-bounded $\hormander{s}$-calculus} for $s>s_L+|1/2-1/p|$.
\end{enumerate}  
\end{lm}
\begin{proof}
\ref{lm:14(1)}. To justify that $L$ is $0$-bisectorial, let us fix $\omega\in (0,\pi/2)$ and $\zeta\in\CC$ such that $\left|\arg(\pm\zeta)\right|>\omega$. Consider the function $F_\zeta : \Rnoz\to\CC$ given by
\begin{equation*}
F_\zeta(\lambda) =(1-\zeta \lambda)^{-1}.
\end{equation*} 
As $L$ has a $\hormander{s}$-bounded functional calculus on $L^p(X)$ for any $s$ sufficiently large, to prove the uniform $L^p$-boundedness of $F_\zeta(L) = (1-\zeta L)^{-1}$ it suffices to show that $\Vert F_\zeta\Vert_{\hormander{s}}\lesssim_\omega 1$ for any $s \geq 1/2$; on the other hand, this follows from the easily observed bound
\begin{equation*}
\sup_{\lambda\in\Rnoz} |\lambda^n\partial_\lambda^n F_\zeta(\lambda) | \lesssim_{n,\omega} 1,\qquad n\in\NN. 
\end{equation*}

\ref{lm:14(tensor)}. 
As $L$ is $0$-bisectorial and has a $\hormander{s}$-bounded functional calculus on $L^p(X)$,
from the definitions and Fubini's theorem it follows immediately that the operator $L \otimes \id$ is $0$-bisectorial and has a $\hormander{s}$-bounded functional calculus on $L^p(X \times Y)$ too; indeed, notice that 
\[
F(L \otimes \id) = F(L) \otimes \id, \qquad \|F(L \otimes \id)\|_{\LinBnd(L^p(X \times Y))} = \|F(L)\|_{\LinBnd(L^p(X))}.
\]
for any bounded Borel function $F$.

\ref{lm:14(2)}.
As $L$ is $0$-bisectorial and has a bounded $\hormander{s}$-functional calculus on $L^p(X)$ for some $s > 1/2$, the claimed estimate follows from \cite[Theorem~4.1 and Section 7]{KriWe16} (see also \cite[Theorem~3.1]{GalMia}).
	
\ref{lm:14(3)}. We shall apply \cite[Proposition~2.14]{DelKri23} to the operator $L$ (observe that it is also valid for $0$-bisectorial operators, cf.\ \cite[Section~7]{KriWe16}). For that purpose recall that $L^p(X)$ has type $\min\{2,p\}$ and cotype $\max\{2,p\}$ (see \cite[pp.~54--58]{HNVW2}), and it also has Pisier's property $(\alpha)$ (see \cite[Sections~4.9--4.10]{KunWei}). As $1/\min\{2,p\}-1/\max\{2,p\} = |1/2-1/p|$ and $L$ has a bounded $\hormander{s}$-functional calculus on $L^p(X)$ for any $s > s_L$, from \cite[Proposition~2.14]{DelKri23} we deduce the bound \eqref{eq:RbdFC} for any $s > s_L + |1/2-1/p|$.
\end{proof}

We now briefly recall some basic definitions and results for operator-valued functions.

\begin{defin}\label{def:opval_smoothness}
Let $(Y,\nu)$ be a separable $\sigma$-finite measure space, and let $M : \Rnoz\to\LinBnd(L^2(Y))$.
\begin{enumerate}[label=(\alph*)]
\item We say that $M$ is a \emph{(weakly) measurable} function if the ``matrix coefficients''
\begin{equation}\label{eq:Mdiagcoeff}
\xi\mapsto\langle M(\xi)f,g\rangle
\end{equation}
are measurable for all $f,g\in L^2(Y)$. Since $L^2(Y)$ is separable, the function
\[
\xi\mapsto \Vert M(\xi)\Vert_{\LinBnd(L^2(Y))}
\]
is also measurable, as the supremum of a countable family of measurable functions.

\item We say that $M$ is \emph{(weakly) continuous} if for any $f,g\in L^2(Y)$ the matrix coefficients \eqref{eq:Mdiagcoeff}
are continuous on $\Rnoz$. In that case, by the Banach--Steinhaus theorem, $M$ is locally bounded, that is,
\[
\sup_{\xi\in K}\Vert M(\xi)\Vert_{\LinBnd(L^2(Y))} < \infty
\]
for any compact set $K\subseteq \Rnoz$. 

\item Let $k \in \NN$. We say that the function $M$ is \emph{(weakly) of class $C^k$} if, for all $f,g\in L^2(Y)$, the matrix coefficients \eqref{eq:Mdiagcoeff} are of class $C^k$.
In this case, another application of the Banach--Steinhaus theorem shows that, for any $\xi\in\Rnoz$, there exist operators in $\LinBnd(L^2(Y))$, which we denote by $\partial^j_\xi M(\xi)$, $j=1,\ldots,k$, such that
\begin{equation}\label{eq:52}
\langle\partial^j_\xi M(\xi)f,g\rangle = \partial^j_\xi \langle M(\xi)f,g\rangle, \qquad f,g\in L^2(Y),\ j=1,\ldots,k;
\end{equation}
moreover each $\xi \mapsto \partial^j_\xi M(\xi)$ is continuous.

\item Assume that $M$ is measurable. We say that $M$ is \emph{(weakly) integrable over a Borel set $A\subseteq \Rnoz$} if the function $\xi\mapsto\Vert M(\xi)\Vert_{\LinBnd(L^2(Y))}$ is $L^1$-integrable over $A$. Then the formula
\begin{equation}\label{eq:54}
\left\langle \int_A M(\xi)\,\dd\xi\, f,g\right\rangle = \int_A \langle M(\xi)f,g\rangle\,\dd\xi, \qquad f,g\in L^2(Y)
\end{equation}
defines an operator $\int_A M(\xi)\,\dd\xi\in\LinBnd(L^2(Y))$. Moreover, there holds
\begin{equation}\label{eq:53}
\left\Vert \int_A M(\xi)\,\dd\xi \right\Vert_{\LinBnd(L^2(Y))} \leq \int_A \Vert M(\xi)\Vert_{\LinBnd(L^2(Y))}\, \dd\xi.
\end{equation}
\end{enumerate}
\end{defin}

For the rest of this section, $(Y,\nu)$ shall be a separable $\sigma$-finite measure space. Consider an operator-valued function $M : \Rnoz\to\LinBnd(L^2(Y))$ and let $p \in (1,\infty)$. If for a fixed $\xi\in\Rnoz$ the operator $M(\xi)$ has a (necessarily unique) bounded extension to $L^p(X)$, then we simply say that $M(\xi)$ is bounded on $L^p(X)$, and use the same notation $M(\xi)$ for the operator on $L^p(X)$ as well. In the next lemma we describe how the properties of $M(\xi)$ on $L^2(X)$ discussed in Definition \ref{def:opval_smoothness} translate to similar properties of $M(\xi)$ on $L^p(X)$.
 
\begin{lm}
Let $M : \Rnoz\to\LinBnd(L^2(Y))$ be measurable. Fix $p\in(1,\infty)$ and let $p'\in(1,\infty)$ be such that $1/p+1/p'=1$. We assume that each $M(\xi)$, $\xi\in\Rnoz$, is also bounded on $L^p(Y)$. Then the following hold.
\begin{enumerate}[label=(\roman*)]
\item\label{it:1,lm:15} The function 
\begin{equation}\label{eq:MLpnorm}
\xi\mapsto\Vert M(\xi)\Vert_{\LinBnd(L^p(Y))}
\end{equation}
is measurable, and the matrix coefficients \eqref{eq:Mdiagcoeff} are measurable for all $f\in L^p(Y)$ and $g\in L^{p'}(Y)$.
\item\label{it:2,lm:15} If $M$ is integrable over a Borel set $A\subseteq \Rnoz$ and the function \eqref{eq:MLpnorm}
is also integrable over $A$, then $\int_A M(\xi)\,\dd\xi$ extends to a bounded operator on $L^p(Y)$ such that
\begin{equation*}
\left\Vert \int_A M(\xi)\,\dd\xi \right\Vert_{\LinBnd(L^p(Y))} \leq \int_A \Vert M(\xi)\Vert_{\LinBnd(L^p(Y))}\, \dd\xi
\end{equation*}
and \eqref{eq:54} is valid for all $f\in L^p(Y)$, $g\in L^{p'}(Y)$.
\item\label{it:3,lm:15} If $M$ is continuous and the function \eqref{eq:MLpnorm} is locally bounded, then the matrix coefficients \eqref{eq:Mdiagcoeff} are continuous for all $f\in L^p(Y)$ and $g\in L^{p'}(Y)$.
\item\label{it:4,lm:15} Let $k\in\NN$. If $M$ is of class $C^k$ and the functions
\begin{equation*}
\xi\mapsto \Vert \partial^j_\xi M(\xi)\Vert_{\LinBnd(L^p(Y))},\qquad j=1,\ldots,k,
\end{equation*}
are locally bounded, then the matrix coefficients \eqref{eq:Mdiagcoeff} are of class $C^k$
and \eqref{eq:52} holds for all $f\in L^p(Y)$ and $g\in L^{p'}(Y)$.
\end{enumerate}
\end{lm}
\begin{proof}
\ref{it:1,lm:15}.
Since $Y$ is $\sigma$-finite we can find a countable family $\dnsA_q \subseteq L^q(Y)\cap L^2(Y) \setminus \{0\}$ dense in $L^q(Y)$ for any $q \in (1,\infty)$. Thus,
\begin{equation*}
\Vert M(\xi)\Vert_{\LinBnd(L^p(Y))}=\sup_{f\in \dnsA_p} \sup_{g\in \dnsA_{p'}} \frac{\vert \langle M(\xi)f,g\rangle\vert}{\Vert f\Vert_{L^p(Y)} \Vert g\Vert_{L^{p'}(Y)}},
\end{equation*}
and therefore \eqref{eq:MLpnorm} is measurable 
as the supremum of countable family of measurable functions.
	
Fix $f\in L^p(Y)$ and $g\in L^{p'}(Y)$, and let $\{f_n \}_{n\in\NN}\subseteq L^p(Y)\cap L^2(Y)$ be a sequence converging to $f$ in $L^p(Y)$; analogously define an approximating sequence $\{g_n\}_{n\in\NN} \subseteq L^{p'}(Y)\cap L^2(Y)$ for $g$. Observe that, for all $\xi \in \Rnoz$,
\begin{equation}\label{eq:49}
\begin{split}
&|\langle M(\xi)f,g\rangle -\langle M(\xi)f_n,g_n\rangle | \\
&\leq |\langle M(\xi)(f-f_n),g\rangle| + |\langle M(\xi)f_n,g-g_n\rangle|\\
&\leq \Vert M(\xi)\Vert_{\LinBnd(L^p(Y))} \left( \Vert f-f_n\Vert_{L^p(Y)} \Vert g\Vert_{L^{p'}(Y)} +  \sup_n \Vert f_n\Vert_{L^p(Y)} \Vert g-g_n\Vert_{L^{p'}(Y)} \right).
\end{split}
\end{equation}
This means that $\langle M(\cdot)f_n,g_n\rangle$ tends to $\langle M(\cdot)f,g\rangle$ pointwise, hence the limit is measurable.
	
\ref{it:2,lm:15}.
For all $f \in L^p(Y)$ and $g \in L^{p'}(Y)$,
\[
\left| \int_A \langle M(\xi)f,g\rangle\, \dd\xi \right|
\leq \int_A \Vert M(\xi)\Vert_{\LinBnd(L^p(Y))}\,\dd\xi\ \Vert f_n\Vert_{L^p(Y)} \Vert g_n\Vert_{L^{p'}(Y)} .
\]
Thus, $(f,g) \mapsto \int_A \langle M(\xi)f,g\rangle \,d\xi$ is a bounded sesquilinear form, and there exists a unique bounded operator $S$ on $L^p(X)$ with
\[
\langle S f, g \rangle = \int_A \langle M(\xi)f,g\rangle \,d\xi
\]
for all $f \in L^p(Y)$ and $g \in L^{p'}(Y)$. In particular, $\langle S f, g \rangle = \langle \int M(\xi) \,d\xi f,g \rangle$ for all $f \in L^2 \cap L^p(Y)$ and $g \in L^2 \cap L^{p'}(Y)$. By density, we deduce that $S= \int M(\xi) \,d\xi$ is bounded on $L^p(Y)$ and the identity \eqref{eq:54} holds for all $f \in L^p(Y)$ and $g \in L^{p'}(Y)$.

\ref{it:3,lm:15} From \eqref{eq:49} we see that $\langle M(\cdot)f,g\rangle$ is a uniform limit on compact sets of $\langle M(\cdot)f_n,g_n\rangle$. Since the latter functions are continuous, we proved \ref{it:3,lm:15}.
	
\ref{it:4,lm:15} As above, $\xi \mapsto \partial_\xi\langle  M(\xi)f_n,g_n\rangle = \langle  \partial_\xi M(\xi)f_n,g_n\rangle$ tends uniformly on compact sets to $\xi \mapsto \langle \partial_\xi M(\xi)f,g\rangle$, hence the latter is a continuous function. Since $\xi \mapsto \langle M(\xi)f_n,g_n\rangle$ tends uniformly on compact sets to $\xi \mapsto \langle M(\xi)f,g\rangle$ too, we deduce that
\begin{equation*}
\partial_\xi \langle M(\xi)f,g\rangle=\langle \partial_\xi M(\xi)f,g\rangle.
\end{equation*}
For higher-order derivatives we proceed similarly.
\end{proof}

Now we want to consider an operator-valued functional calculus for $L$. For that purpose, we are now going to give a meaning to $M(L)$, where $M : \Rnoz \to \LinBnd(L^2(Y))$ is a suitable operator-valued function. 

As $L$ is self-adjoint and injective on $L^2(X)$, recall that by the spectral theorem (see, e.g., \cite[Theorem~VIII.4]{RS1}) there exist a measure space $\Omega$, a measurable function $\ell : \Omega\to\Rnoz$ and a unitary operator $\Upsilon : L^2(X)\to L^2(\Omega)$, such that 
\[
\Upsilon(Lf)(\omega) = \ell(\omega) \Upsilon f(\omega),\qquad f\in \Dom(L),
\]
i.e., $\Upsilon$ intertwines $L$ with the operator of multiplication by $\ell$.
Moreover, for any bounded Borel function $F : \Rnoz\to\CC$,
\[
\Upsilon(F(L)f)(\omega) = F(\ell(\omega)) \Upsilon f(\omega),\qquad f\in L^2(X),
\]
i.e., $\Upsilon$ intertwines the (scalar) Borel functional calculus for $L$ with that for the multiplication operator.

Let $M : \Rnoz\to\LinBnd(L^2(Y))$ be measurable and bounded. Much in this spirit we can define $M(L)$ on $L^2(X\times Y)$. Let $\widetilde{\Upsilon} : L^2(X\times Y)\to L^2(\Omega\times Y)$ be the unitary operator given by $\widetilde{\Upsilon}=\Upsilon\otimes \id$. Then, we define $M(L)$ so that
\begin{equation}\label{eq:55}
\widetilde{\Upsilon}(M(L)f)(\omega,y) = \left(M(\ell(\omega)) \, \widetilde{\Upsilon} f(\omega,\cdot)\right)(y),\qquad f\in L^2(X\times Y).
\end{equation}

\begin{rem}
Some remarks are in order.
\begin{enumerate}[label=(\alph*)]
\item The above definition of $M(L)$ can be equivalently expressed in the language of direct integrals (see, e.g., \cite[Section 7.4]{Fo}): indeed, 
the function $\omega \mapsto M(\ell(\omega))$ is a measurable field of operators with respect to the constant field $\omega \mapsto L^2(Y)$ of Hilbert spaces on $\Omega$, and
under the identification $L^2(\Omega \times Y) \cong \int_\Omega^{\oplus} L^2(Y) \,d\omega$
we can write
\[
\widetilde{\Upsilon} M(L) \widetilde{\Upsilon}^{-1} = \int_{\Omega}^\oplus M(\ell(\omega)) \,d\omega.
\]

\item From the definition, it is clear that
\[
\|M(L)\|_{\LinBnd(L^2(X\times Y))} = \esssup_{\omega \in \Omega} \|M(\ell(\omega))\|_{\LinBnd(L^2(Y))} \leq \sup_{\xi \in \Rnoz} \|M(\xi)\|_{\LinBnd(L^2(Y))}.
\]
\item The definition of $M(L)$ is independent of the choice of the unitary operator $\Upsilon$ given by the spectral theorem for $L$. Indeed, observe that, if we take $f=f_1\otimes f_2$, $g=g_1\otimes g_2\in L^2(X\times Y)$, then
\begin{equation}\label{eq:matcoeff}
\begin{aligned}
\langle M(L)f,g\rangle_{L^2(X\times Y)} 
&= \int_\Omega \Upsilon f_1(\omega) \, \overline{\Upsilon g_1(\omega)} \int_Y M(\ell(\omega)) (f_2)(y) \, \overline{g_2(y)}\,\dd\nu(y)\,\dd\omega\\
&= \int_\Omega \Upsilon f_1(\omega) \, \overline{\Upsilon g_1(\omega)} F(\ell(\omega)) \,\dd\omega\\
&= \langle F(L)f_1,g_1\rangle_{L^2(X)},
\end{aligned}
\end{equation}
where $F(\xi)= \langle M(\xi)f_2,g_2\rangle_{Y}$; by the density of the span of tensor products in $L^2(X\times Y)$ we obtain the required independence.
\item The expression \eqref{eq:matcoeff} for matrix coefficients also shows that, if $M_n \to M$ pointwise boundedly, i.e.,
\begin{gather*}
\sup_{n \in \NN} \sup_{\xi \in \Rnoz} \|M_n\|_{\LinBnd(L^2(Y))} < \infty, \\
\lim_{n \to \infty} \langle M_n(\xi) \phi,\psi \rangle = \langle M(\xi) \phi,\psi \rangle \qquad\forall \xi \in \Rnoz, \ \phi,\psi \in L^2(Y),
\end{gather*}
then $M_n(L) \to M(L)$ in the weak operator topology of $\LinBnd(L^2(X \times Y))$.
\item If $F : \Rnoz \to \CC$ is a bounded Borel function and $\tilde F : \Rnoz \to \LinBnd(L^2(Y))$ is defined by $\tilde F(\xi) \defeq F(\xi) \id$, then clearly
\[
\tilde F(L) = F(L) \otimes \id;
\]
in other words, the operator-valued functional calculus is an extension of the scalar-valued one. Thus, in what follows, when there is no risk of confusion, we shall just write $F(L)$ in place of $F(L)\otimes \id$.
\end{enumerate}
\end{rem}

We are now interested in $L^p$-boundedness properties of operators of the form $M(L)$, i.e., operator-valued spectral multipliers of $L$.

\begin{thm}\label{thm:5}
Fix $p\in(1,\infty)$ and $s_L \geq 1/2$. Let $L$ be a self-adjoint operator on $L^2(X)$ having a bounded $\hormander{s}$-calculus on $L^p(X)$ for all $s>s_L$. Let $M : \Rnoz\to\LinBnd(L^2(Y))$ be measurable and bounded. Assume that $M$ is of class $C^{N}$ for some integer $N>s_L+3/2$ and that
\begin{equation}\label{eq:46}
C_{N,p}(M) \defeq \RBd_{L^p(Y)} \{\xi^j\partial_\xi^j M(\xi) \tc j=0,1,\ldots,N,\ \xi\in\Rnoz\} <\infty.
\end{equation}
Then, $M(L)$ is a bounded operator on $L^p(X\times Y)$ with
\[
\|M(L)\|_{\LinBnd(L^p(X \times Y))} \lesssim_{N,p,L} C_{N,p}(M).
\]
\end{thm}
\begin{proof}
Fix an even smooth cutoff $\breve\chi$ supported in $[-2,-1/2]\cup[1/2,2]$ and satisfying \eqref{eq:dyadicpart}.
Thus, if we set
\begin{equation}\label{eq:dyadicM}
M_m(\xi)\defeq \breve\chi(\xi)M(2^{-m}\xi),
\end{equation}
then we can decompose
\begin{equation*}
	M(\xi)=\sum_{m\in\ZZ} M_m (2^m\xi),\qquad \xi\in\Rnoz,
\end{equation*}
with convergence in the weak operator topology of $\LinBnd(L^2(Y))$. Since the $M(\xi)$ are uniformly bounded on $L^2(Y)$ we also have
\begin{equation*}
	M(L) = \sum_{m\in\ZZ} M_m(2^m L)
\end{equation*}
in the weak operator topology of $\LinBnd(L^2(X\times Y))$.

Thus, by Lemma \ref{lm:14}, for all $f \in L^2 \cap L^p(X\times Y)$,
\[\begin{split}
\Vert M(L) f\Vert_{L^p(X\times Y)} 
&\simeq_{p} \int_0^1 \left\Vert \sum_{m\in\ZZ} \rade_m(t) \breve\chi(2^m L) M(L) f\right\Vert_{L^p(X\times Y)}\,\dd t\\
&\leq \sum_{j\in\{-1,0,1\}} \int_0^1 \left\Vert \sum_{m\in\ZZ} \rade_m(t) \breve\chi(2^m L) M_{m+j}(2^{m+j}L) f\right\Vert_{L^p(X\times Y)}\,\dd t,
\end{split}\]
since the supports of $\breve\chi(2^{m_1}\cdot)$ and $\breve\chi(2^{m_2}\cdot)$ overlap only if $|m_1-m_2|\leq 1$. Further, by Lemma \ref{lm:14} and Definition \ref{def:Rbound},
\[\begin{split}
&\Vert M(L) f \Vert_{L^p(X\times Y)} \\
&\lesssim_{p} \sum_{j\in\{-1,0,1\}} \int_0^1 \left\Vert \sum_{m\in\ZZ} \rade_m(t) M_{m+j}(2^{m+j} L) \breve\chi(2^m L)  f\right\Vert_{L^p(X\times Y)}\,\dd t\\
&\lesssim \RBd_{L^p(X\times Y)} \{ M_m(2^m L)\}_{m\in\ZZ} \int_0^1 \left\Vert \sum_{m\in\ZZ} \rade_m(t) \breve\chi(2^m L)  f\right\Vert_{L^p(X\times Y)}\,\dd t\\
&\simeq_{p} \RBd_{L^p(X\times Y)} \{ M_m(2^m L)\}_{m\in\ZZ} \ \Vert f\Vert_{L^p(X\times Y)}.
\end{split}\]
This means that
\[
\| M(L)\|_{\LinBnd(L^p(X \times Y))} \lesssim_p \RBd_{L^p(X\times Y)} \{ M_m(2^m L) \}_{m\in\ZZ},
\]
so we are reduced to proving the R-boundendess of $\{M_m(2^m L)\}_{m\in \ZZ}$ on $L^p(X\times Y)$.
	
For that purpose, we introduce Fourier coefficients of the functions $\xi \mapsto M_m(\xi)$, $m \in \ZZ$, which are supported in $[-2,-1/2] \cup [1/2,2] \subseteq (-\pi,\pi)$ by construction. Namely, let $\hat{M}_m(k)\in\LinBnd(L^2(Y))$, $k\in\ZZ$, be given by
\[
\hat{M}_m(k)\defeq \frac{1}{2\pi} \int_{-\pi}^{\pi} M_m(\xi) \, e^{-i k \xi}\,\dd\xi.
\]
As $M$ is measurable and bounded, the above integral is well defined.
Actually, as $M$ is of class $C^N$, for all $n\leq N$ and $m,k \in \ZZ$ there holds
\begin{equation}\label{eq:38}
	(-i k)^n \hat{M}_m(k) = \frac{1}{2\pi}  \int_{-\pi}^\pi\partial_\xi^n M_m(\xi) \,e^{-i k\xi}\,\dd\xi,
\end{equation}
as one can readily see via repeated integration by parts. As each $M_m$ is of class $C^N$, from \eqref{eq:38} and 
\eqref{eq:53} we deduce that
\[
\Vert \hat{M}_m(k)\Vert_{\LinBnd(L^2(Y))}\lesssim_{m} (1+|k|)^{-N}.
\]

As $N \geq 2$, the decay of the Fourier coefficients justifies the application of the Fourier inversion formula; namely, if we fix $\bar{\chi}\in C_c^\infty(\RR)$ such that $\ind_{[-2,-1/2]\cup[1/2,2]}\leq \bar{\chi}\leq\ind_{[-9/4,-1/4]\cup[1/4,9/4]}$, then
\begin{equation*}
M_m(\xi)=\sum_{k\in\ZZ} \bar{\chi}(\xi) e^{i k \xi} \hat{M}_{m}(k) ,\qquad \xi\in\Rnoz,
\end{equation*}
with convergence in the weak operator topology of $\LinBnd(L^2(Y))$.
Consequently, if we set
\begin{equation}\label{eq:48}
E_{k,m}(\xi)\defeq \bar{\chi}(2^m\xi)e^{i k2^m\xi},
\end{equation}
then
\begin{equation}\label{eq:47}
	M_m(2^m L) = \sum_{k\in\ZZ} E_{k,m}(L) \otimes \hat{M}_m(k) 
\end{equation}
in the weak operator topology of $\LinBnd(L^2(X\times Y))$.

The decomposition \eqref{eq:47} is the key tool that allows us to ``decouple'' the problem of $L^p(X \times Y)$-boundedness for operator-valued multipliers $M_m(2^m L)$ of $L$ into $L^p(X)$-bounds for scalar-valued multipliers $E_{k,m}(L)$ and $L^p(Y)$-bounds for the Fourier coefficients $\hat{M}_m(k)$, which in turn are related to the smoothness of $M$.

Indeed, let $\varepsilon>0$ be such that $N=s_L+3/2+2\varepsilon$. Then, from the decomposition \eqref{eq:47}, by \cite[Propositions~8.1.24~and~8.1.19(3)]{HNVW2} we deduce that
\[\begin{split}
&\RBd_{L^p(X\times Y)} \{ M_m(2^m L) \}_{m\in\ZZ} \\ 
&\leq \sum_{k\in\ZZ} (1+|k|)^{-(1+\varepsilon)} \, \RBd_{L^p(X\times Y)} \{(1+|k|)^{-(s_L+1/2 +\varepsilon)} E_{k,m}(L) \otimes \id\}_{m\in\ZZ}\\
&\qquad\times \RBd_{L^p(X\times Y)} \{ (1+|k|)^{N}\id\otimes \hat{M}_m(k)\}_{m\in\ZZ} \\
&\lesssim_{\varepsilon}  \RBd_{L^p(X\times Y)} \{(1+|k|)^{-(s_L+1/2 +\varepsilon)} E_{k,m}(L) \otimes \id\}_{m,k\in\ZZ }\\
&\qquad\times \RBd_{L^p(X\times Y)} \{ (1+|k|)^{N} \id \otimes \hat{M}_m(k)\}_{m,k\in\ZZ} \\
&= \RBd_{L^p(X)} \{(1+|k|)^{-(s_L+1/2 +\varepsilon)} E_{k,m}(L)\}_{m,k\in\ZZ}\
 \RBd_{L^p(Y)} \{ (1+|k|)^{N} \hat{M}_m(k)\}_{m,k\in\ZZ},
\end{split}\]
where the last equality is due to Fubini's theorem.
Thus, the R-boundedness on $L^p(X\times Y)$ of $\{M_m(2^mL)\}_m$ boils down to that of $\{(1+|k|)^{-(s_L+1/2 +\varepsilon)} E_{k,m}(L) \}_{m,k}$ on $L^p(X)$ and that of $\{ (1+|k|)^{N} \hat{M}_m(k)\}_{m,k}$ on $L^p(Y)$.
	
For the R-boundedness of the former family, notice that, by \eqref{eq:48},
\begin{equation*}
\sup_{\xi\in\Rnoz} |\xi^j \partial_\xi^j E_{k,m}(\xi)| \lesssim_j (1+|k|)^j, \qquad j\in\NN,
\end{equation*}
whence it follows that
\begin{equation}\label{eq:EkmHormander}
\Vert (1+|k|)^{-s} E_{k,m} \Vert_{\mathcal{H}_2^{s}} \lesssim_s 1
\end{equation}
for all $s \geq 1/2$.
Now, recall that $L$ has a bounded $\hormander{s}$-calculus on $L^p(X)$ for all $s>s_L$. Hence, by Lemma \ref{lm:14}\ref{lm:14(3)}, $L$ has an R-bounded $\hormander{s}$-calculus on $L^p(X)$ for all $s > s_L+|1/2-1/p|$, and therefore from \eqref{eq:EkmHormander} we deduce that
\begin{equation*}
	\RBd_{L^p(X)} \{(1+|k|)^{-(s_L+1/2 +\varepsilon)} E_{k,m}(L) \tc m,k\in\ZZ \}<\infty.
\end{equation*}
	
For the R-boundedness of the other family, observe that by approximating the integrals in \eqref{eq:38} by Riemann sums we obtain, for any $n \leq N$,
\begin{equation*}
\{|k|^n \hat{M}_m(k) \tc m,k\in\ZZ \} \subseteq \overline{\absconv \{ \partial_\xi^n M_m(\xi) \tc \xi\in\Rnoz,\ m\in\ZZ \}},
\end{equation*}
where $\absconv$ denotes the absolute convex hull \cite[Definition~3.2.12(2)]{HNVW1} and the closure is in the weak operator topology of $\LinBnd(L^p(Y))$.
Hence, by \cite[Propositions~8.1.19(1),~8.1.21~and~8.1.22]{HNVW2},
\begin{equation*}
\RBd_{L^p(Y)} \{ |k|^n \hat{M}_m(k) \tc m,k\in\ZZ \} \leq \RBd_{L^p(Y)} \{ \partial_\xi^n M_m(\xi) \tc \xi\in\Rnoz,\ m\in\ZZ \},
\end{equation*}
thus also, by \cite[Proposition~8.1.19(2)]{HNVW2},
\begin{multline*}
\RBd_{L^p(Y)} \{ (1+|k|)^{N}\hat{M}_m(k) \tc m,k\in\ZZ \} \\
\lesssim_N \RBd_{L^p(Y)} \{ \partial_\xi^{n} M_m(\xi) \tc \xi\in\Rnoz,\ m\in\ZZ, \ n \leq N \}.
\end{multline*}
Now, by \eqref{eq:dyadicM},
\begin{equation*}
\partial_\xi^n M_m(\xi) = \sum_{j=0}^n \binom{n}{j} \frac{\breve\chi^{(n-j)}(\xi)}{\xi^j} ( \xi^j\partial_\xi^j M )(2^{-m}\xi);
\end{equation*}
thus, by the Kahane contraction principle \cite[Proposition~3.2.10]{HNVW1},
\begin{multline*}
\RBd_{L^p(Y)} \{ \partial_\xi^{n} M_m(\xi) \tc \xi\in\Rnoz,\ m\in\ZZ, \ n \leq N \}\\
\lesssim_{\chi,N} \RBd_{L^p(Y)} \{\xi^j \partial_\xi^j M(\xi) \tc j\leq N,\ \xi\in\Rnoz \} .
\end{multline*}
By combining the above estimate one finally gets that
\begin{multline*}
\RBd_{L^p(Y)} \{ (1+|k|)^{N}\hat{M}_m(k) \tc m,k\in\ZZ \}  \\
\lesssim_{\chi,N} \RBd_{L^p(Y)} \{\xi^j \partial_\xi^j M(\xi) \tc j\leq N,\ \xi\in\Rnoz \} = C_{N,p}(M),
\end{multline*}
and the latter quantity is finite by our assumption.
\end{proof}

Let $A_2(\RR^d)$, $d\geq 1$, be the Muckenhoupt class; for a weight $w\in A_2(\RR^d)$ we denote its $A_2$-characteristic by $[w]_{A_2}$. In the case where $Y=\RR^d$ with the Lebesgue measure, the smoothness condition on the multiplier $M$ in Theorem \ref{thm:5} admits a more concrete reformulation in terms of uniform weighted $L^2$-boundedness.

\begin{cor}\label{cor:3}
Let $s_L \geq 1/2$ and $p \in (1,\infty)$. Let $L$ be a self-adjoint operator on $L^2(X)$ having a bounded $\hormander{s}$-calculus on $L^p(X)$ for all $s>s_L$. Let $N > s_L + 3/2$ be an integer. Let $M : \Rnoz\to\LinBnd(L^2(\RR^d))$ be of class $C^{N}$. Assume that there exists a nondecreasing function $\psi : [1,\infty) \to [0,\infty)$ such that, for all $w\in A_2(\RR^d)$,
\begin{equation}\label{eq:A2bound}
\max_{j=0,\dots,N} \sup_{\xi\in\Rnoz} \Vert \xi^j\partial_\xi^j M(\xi) \Vert_{\LinBnd(L^2(w))} \leq \psi([w]_{A_2}).
\end{equation}
Then, $M(L)$ is a bounded operator on $L^p(X\times \RR^d)$.
\end{cor}
\begin{proof}
By \cite[Theorem~8.2.6]{HNVW2}, from \eqref{eq:A2bound} it follows that the family
\[
\cT = \{ \xi^j \partial_\xi^j M(\xi) \tc \xi \in \Rnoz,\ j=0,\dots,N \}
\]
is $\ell^2$-bounded on $L^p(\RR^d)$ (in the sense of \cite[Definition~8.1.1]{HNVW2}). In turn, by \cite[Theorem 8.1.1(3)]{HNVW2}, this implies the R-boundedness of $\cT$ on $L^p(\RR^d)$, so Theorem \ref{thm:5} can be applied.
\end{proof}

\begin{rem}
If the self-adjoint operator $L$ is non-negative, then $M(L)$ only depends on the restriction of $M$ to the positive half-line $\Rpos$. As a consequence, in this case one can replace $\Rnoz$ with $\Rpos$ in the smoothness conditions \eqref{eq:46} and \eqref{eq:A2bound} of Theorem \ref{thm:5} and Corollary \ref{cor:3}.
\end{rem}

\section{\texorpdfstring{$L^p$}{Lp}-boundedness of the Riesz transforms for \texorpdfstring{$p\geq 2$}{p>=2}}\label{s:adjriesz}

In this section, we complete the proof of Theorem \ref{thm:main_riesz}, by proving the $L^p$-boundedness of the Riesz transforms $\Riesz_0$ and $\Riesz_1^*$ for $p \geq 2$. As discussed in Section \ref{ss:decomposition}, we are reduced to proving the boundedness of the operators $\tilde\Riesz_0^{(j)}$, $j=1,2,3$, and $\tilde\Riesz_1^{(j)}$, $j=1,2$, for $p \in (1,2]$.

\subsection{An auxiliary estimate}
We shall make use of the following estimate for the function $\aux^\nu$ defined in \eqref{eq:aux}.

\begin{lm}\label{lm:16}
For any $n\in\NN$ and $\varepsilon \in (0,1/2]$ the following estimates hold:
\begin{align*}
|\xi^n (\Hank_\nu \aux^\nu)^{(n)}(\xi)| &\lesssim_{n,\nu,\varepsilon} e^{-(1-\varepsilon)\xi},\\
|\partial_\xi(\xi^n (\Hank_\nu \aux^\nu)^{(n)}(\xi))| &\lesssim_{n,\nu,\varepsilon} e^{-(1-\varepsilon)\xi},\\
|\xi^{n+1} (\Hank_{\nu+2} \aux^\nu)^{(n)}(\xi)| &\lesssim_{n,\nu,\varepsilon} (\xi e^{-\xi})^{1-\varepsilon}
\end{align*}
for all $\xi \in \Rpos$; in particular,
\begin{equation*}
|\xi^n  (\Hank_\nu \aux^\nu)^{(n)}(\xi)-(\xi')^n  (\Hank_\nu \aux^\nu)^{(n)}(\xi')|\lesssim_{n,\nu} |\xi-\xi'|
\end{equation*}
for all $\xi,\xi' \in \Rpos$.
\end{lm}
\begin{proof}
By the definition \eqref{eq:H_transf} of Hankel transform, together with \eqref{eq:aux} and \cite[(10.22.46), (10.27.3) and (10.32.8)]{DLMF},
\begin{align*}
\Hank_\nu \aux^\nu(\xi) &= \frac{1}{\nu} \xi K_{-1}(\xi) = \frac{1}{\nu} \xi^2 \int_1^\infty e^{-\xi t} \sqrt{t^2-1}\,\dd t,\\
\Hank_{\nu+2} \aux^\nu(\xi) &= K_0(\xi) = \int_1^\infty e^{-\xi t}\frac{\dd t}{\sqrt{t^2-1}},
\end{align*}
where $K_s$ denotes the modified Bessel function of the second kind of order $s$.

As a consequence, for all $n\in\NN$,
\begin{align*}
\xi^n (\Hank_\nu\aux^\nu)^{(n)}(\xi) &=\frac{1}{\nu} \sum_{j=0}^{\min\{n,2\}} c_{n,j} F_{n-j}(\xi),\\
\xi^{n} (\Hank_{\nu+2}\aux^\nu)^{(n)} &= (-1)^n G_n(\xi),
\end{align*}
where $c_{n,j}$ are some constants, while
\begin{align*}
F_n(\xi) &= \xi^2 \int_1^{\infty} e^{-t\xi} (t\xi)^n \sqrt{t^2-1}\,\dd t = \int_0^\infty e^{-(t+\xi)} (t+\xi)^n \sqrt{t} \sqrt{t+2\xi}\,\dd t,\\
G_n(\xi) &= \int_1^\infty e^{-\xi t} (\xi t)^n \frac{\dd t}{\sqrt{t^2-1}} =  \int_0^\infty e^{-(t+\xi)} (t+\xi)^n \frac{\dd t}{\sqrt{t}\sqrt{t+2\xi}}. 
\end{align*}
Thus, it suffices to prove the bounds in the statement of the lemma for $F_n(\xi)$ and $\xi G_n(\xi)$ in place of $\xi^n (\Hank_\nu\aux^\nu)^{(n)}(\xi)$ and $\xi^{n+1} (\Hank_{\nu+2}\aux^\nu)^{(n)}(\xi)$.

First, for all $\varepsilon \in (0,1)$,
\begin{equation*}
| F_n(\xi) | \leq \sqrt{2} \int_0^\infty e^{-(t+\xi)} (t+\xi)^{n+1} \,\dd t \lesssim_{n,\varepsilon} e^{-(1-\varepsilon)\xi}.
\end{equation*}
Further,
\[
F_n'(\xi) = -F_n(\xi) + n F_{n-1}(\xi) + \tilde F_n(\xi),
\]
where
\[
\tilde F_n(\xi) = \int_0^\infty e^{-(t+\xi)} (t+\xi)^n \sqrt{\frac{t}{t+2\xi}}\,\dd t \leq \int_0^\infty e^{-(t+\xi)} (t+\xi)^n \,\dd t \lesssim_{n,\varepsilon} e^{-(1-\varepsilon)\xi};
\]
thus, from the bound for $F_n(\xi)$ we also deduce that
\[
|F_n'(\xi)| \lesssim_{n,\varepsilon} e^{(1-\varepsilon)\xi}.
\]
This clearly implies that $F_n$ is Lipschitz-continuous, as required.

Finally, for all $\varepsilon \in (0,1/2]$,
\[
G_n(\xi) \leq (2\xi)^{-\varepsilon} \int_0^\infty e^{-(t+\xi)} (t+\xi)^n \,\frac{\dd t}{t^{1-\varepsilon}} \lesssim_{\varepsilon,n} \xi^{-\varepsilon} e^{-(1-\varepsilon)\xi},
\]
and the estimate $|\xi G_n(\xi)| \lesssim_{\varepsilon,n} (\xi e^{-\xi})^{1-\varepsilon}$ follows.
\end{proof}

\subsection{The Riesz transform \texorpdfstring{$\Riesz_0$}{R0}}

We shall prove that the operators $\tilde\Riesz_0^{(j)}$, $j=1,2,3$, from \eqref{eq:R_0_decomp}, are bounded on $L^p(G_\nu)$, $p \in (1,2]$. This completes the proof of the $L^p$-bounds for $\Riesz_0$ stated in Theorem \ref{thm:main_riesz}.

The operator $\tilde\Riesz_0^{(3)}$ is relatively easy to discuss: indeed, it can be factorised into two operators, which can be easily discussed by looking separately at the factors $X_\nu$ and $\RR$ of the semidirect product $G_\nu$.

\begin{prop}\label{prop:K03_bd}
The operator
\begin{equation*}
\tilde\Riesz_0^{(3)}: f\mapsto f\diamond_\nu K_0^{(3)}
\end{equation*}
is of weak type $(1,1)$ and bounded for all $p\in(1,\infty)$.
\end{prop}
\begin{proof}
Notice that, from \eqref{eq:20} and \eqref{eq:39}, it follows that
\[
f\diamond_\nu K_{0}^{(3)}(x,u)= -\frac{2\nu}{\pi} \int_\RR \big( f^v\ast_\nu \aux^\nu_{(e^v)}\big)(x) \frac{\ind_{\{|u-v|\geq1 \} }}{u-v}\,\dv.
\]
Much as in the proof of \cite[Proposition 4.4]{Ma23}, this can be rewritten as
\[
f\diamond_\nu K_{0}^{(3)}= -\frac{2\nu}{\pi} ABf,
\]
where
\begin{equation*}
(B f)^u = f^u\ast_\nu \aux^\nu_{(e^u)},\qquad (Af)_x = f_x\ast k,
\end{equation*} 
while $f_x(u)=f^u(x) = f(x,u)$ and $k(u)= \frac{\ind_{\{|u|\geq 1 \}}}{u}$. Since $\aux^\nu\in L^1(X_\nu)$ and the scaling $\aux^\nu_{(e^u)}$ does not change the $L^1$-norm, we see that $B$ is bounded on $L^p(G_\nu)$, $p\in[1,\infty]$. On the other hand, $A = \id \otimes \tilde A$, where $\tilde A$ is a convolution operator on $\RR$ with a truncated Calder\'on--Zygmund kernel $k$, so it is of weak type $(1,1)$ and bounded on $L^p(G_\nu)$, $p\in(1,\infty)$. Thus, their composition $AB$ has the claimed boundedness properties.
\end{proof}

We now move to the operator $\tilde\Riesz_0^{(2)}$, which shall be treated via the operator-valued multiplier theorem from Section \ref{s:opval}.

\begin{prop}\label{prop:12}
The operator
\begin{equation*}
\tilde\Riesz_0^{(2)} : f\mapsto f\diamond_\nu K_0^{(2)}
\end{equation*}
is bounded on $L^p(G_\nu)$ for all $p\in(1,\infty)$.
\end{prop}
\begin{proof}
We shall employ the abstract theory proved in Section \ref{s:opval} for $L=L_\nu$.
In other words, we take $X=X_\nu$, $\Omega=X_\nu$, $\Upsilon=(2^{\nu/2-1} \Gamma(\nu/2))^{-1} \Hank_\nu$ and $\ell(\omega)=\omega^2$. Recall from \eqref{eq:hankel_conv_transf} that, if $K\in L^2(X_\nu)$ is bounded, then
\begin{equation*}
f\ast_\nu K = \Hank_\nu K(\sqrt{L_\nu})f.
\end{equation*}
Moreover, $L_\nu$ has a bounded $\hormander{s}$-calculus on $L^p(X_\nu)$ for all $s> \frac{\nu}{2}$ \cite{Kap,KanPre,Stk1}. Thus, we can apply Corollary \ref{cor:3} with $N=\lfloor\frac{\nu+5}{2} \rfloor$, that is, the smallest integer greater than $\frac{\nu+3}{2}$.

Recall that $G_\nu=X_\nu\rtimes\RR$ and $\dmub(x,u) = \dmu(x) \,\du$. We denote
\begin{equation*}
\HaId\defeq \Hank_\nu \otimes \id : L^2(G_\nu)\to L^2(G_\nu).
\end{equation*}
Clearly, $\HaId$ is a multiple of an isometry on $L^2(G_\nu)$. Now let $M : X_\nu\to\LinBnd(L^2(\RR))$ be measurable and uniformly bounded. Then, from \eqref{eq:55} with $Y=\RR$ we have 
\begin{equation}\label{eq:56}
\HaId(M(L_\nu) f)(\xi,u) = (M(\xi^2)\HaId f(\xi,\cdot))(u).
\end{equation}

On the other hand, observe that, by \eqref{eq:HaId_conv},
\begin{equation}\label{eq:57}
\HaId (f\diamond_\nu K)(\xi,u) = \int_\RR \HaId f(\xi,v) \, \HaId K(e^{v}\xi,u-v)\,\dd v,
\end{equation}
Let $B_K : X_\nu\to \LinBnd(L^2(\RR))$ be given by
\begin{equation*}
B_K(\xi) g = \int_\RR B_K^\xi(u,v)\, g(v) \, \dv,\qquad g\in L^2(\RR),
\end{equation*}
where
\begin{equation}\label{eq:kernel_conv_op}
B_K^\xi(u,v) \defeq \HaId K(e^v \xi,u-v).
\end{equation}
With this notation, by comparing \eqref{eq:56} and \eqref{eq:57}, we see that
\begin{equation*}
f\diamond_\nu K = B_K(\sqrt{L_\nu}) f.
\end{equation*}

We shall apply this for $K=K_{0}^{(2)}$.
In this case, by \eqref{eq:kernel_conv_op} and \eqref{eq:20},
\begin{equation*}
B_{K_0^{(2)}}^\xi(u,v) = -\frac{2\nu}{\pi} \ind_{\{u \geq v+1\}} \frac{\Hank_\nu\aux^\nu(e^u\xi) -\Hank_\nu\aux^\nu(e^v \xi)}{u-v}.
\end{equation*}
As a consequence, for any $j \in \NN$, the operator $\xi^j\partial_\xi^j B_{K_0^{(2)}}(\xi)$
has integral kernel 
\begin{equation*}
\xi^j \partial_\xi^j B_{K_0^{(2)}}^\xi(u,v) = -\frac{2\nu}{\pi} \ind_{\{u \geq v+1\}} \frac{(\xi^j \partial_\xi^j \Hank_\nu\aux^\nu)(e^u\xi) - (\xi^j \partial_\xi^j \Hank_\nu\aux^\nu)(e^v \xi)}{u-v}.
\end{equation*}
The estimates of Lemma \ref{lm:16}, together with \cite[Proposition~5.6(i)]{Ma23}, then show that each $\xi^j\partial_\xi^j B_{K_0^{(2)}}(\xi)$ is bounded on $L^2(w)$ for any $w\in A_2(\RR)$, with a bound only depending on $j$ and $[w]_{A_2}$. Thus, by Corollary \ref{cor:3}, we deduce the desired $L^p(G_\nu)$-bound for $f \mapsto f \diamond_\nu K_0^{(2)}$.
\end{proof}

We remain with the operator $\tilde\Riesz_0^{(1)}$, which effectively contains the ``local part'' of $\Riesz_0 - \Riesz_0^*$ and can be treated with the standard singular integral theory on $G_\nu$.

\begin{prop}\label{prop:13}
The operator
$\tilde\Riesz_0^{(1)}$
is of weak type $(1,1)$ and bounded for all $p\in(1,2]$.
\end{prop}
\begin{proof}
By Propositions \ref{prop:K03_bd} and \ref{prop:12} we know that the operators corresponding to $K_0^{(3)}$ and $K_{0}^{(2)}$ are bounded on $L^2(G_\nu)$. Also $\Riesz_0-\Riesz_0^\ast$ has this property. Thus, by the decomposition \eqref{eq:R_0_decomp}, we know that $\tilde\Riesz_0^{(1)}$ is bounded on $L^2(G_\nu)$. 
We shall now apply Theorem \ref{thm:3} (and Corollary \ref{cor:4}) to prove the claimed $L^p$-boundedness of $\tilde\Riesz_0^{(1)}$.
	
For all $\bdx = (x,u) \in G_\nu \setminus \{\bdzero\}$, set
\begin{equation*}
K_n(\bdx) = -\frac{2\nu}{\pi} \frac{u}{r^{\nu+2}}  \chi(2^n r),
\end{equation*}
where $r=|\bdx|_\dist$ is as in \eqref{eq:distance_origin}, while $\chi \in C^\infty_c(\Rnoz)$ is as in \eqref{eq:dyadicpart}. Then, by \eqref{eq:20}, we can decompose the off-diagonal $\diamond_\nu$-convolution kernel $K_0^{(1)}$ of $\tilde\Riesz_0^{(1)}$ as
\begin{equation*}
K_0^{(1)} = \sum_{n\in\NN} K_n + E,
\end{equation*}
where $E$ is an $L^1(G_\nu)$-integrable function.
Thus, to prove the $L^p$-boundedness of $\Riesz_0^{(1)}$ by means of Theorem \ref{thm:3}, it suffices to prove that
\begin{equation}\label{eq:42}
\int_{G_\nu} |K_n(\bdx)| (1+2^n r)\,\dmub(\bdx)\lesssim 1,\qquad \int_{G_\nu} \big| \nabla_\nu K^\ast_n(\bdx)\big| \,\dmub(\bdx)\lesssim 2^n.
\end{equation}
Recall that, since $r\lesssim 1$ on the support of $K_n$, there holds $r^2 \simeq x^2+u^2$ (see \eqref{eq:41}).

Thus, by using that $2^n r\simeq 1$ and $|u|x^{\nu-1} \lesssim r^{\nu}$ on the support of $K_n$,
\begin{equation*}
\int_{G_\nu} |K_n(\bdx)| (1+2^n r)\,\dmub(\bdx) \simeq_\nu  \int_{G_\nu} |K_n(\bdx)| \,\dmub(\bdx)\lesssim  \int_{2^{-(n+1)}}^{2^{-(n-1)}} r^{-1}\, \dd r\simeq 1,
\end{equation*}
which proves the the first bound in \eqref{eq:42}.

For the second bound in \eqref{eq:42} notice that
\begin{equation*}
K^\ast_n(\bdx) = \frac{2\nu}{\pi} \frac{e^{-\nu u} u}{r^{\nu+2}}  \chi(2^n r).
\end{equation*}
Therefore, by \eqref{eq:vfs} and \eqref{eq:44}, if we set $\tilde\chi(t) = t \chi'(t) - (\nu+2)\chi(t)$, then
\begin{align*}
\vfX_0 K^\ast_{n}(\bdx) &= \frac{2\nu}{\pi} e^{-\nu u} \frac{1-\nu u}{r^{\nu+2}}\chi(2^n r) 
- \frac{2\nu}{\pi} e^{-\nu u} u \frac{\tilde\chi(2^n r)}{r^{\nu+3}}\cdot\frac{\sinh u -e^{-u}x^2/2}{\sinh r} ,\\
\vfX_1 K^\ast_{n}(\bdx) &= \frac{2\nu}{\pi} e^{-\nu u} u  \frac{\tilde\chi(2^n r) }{r^{\nu+3}} \cdot \frac{x}{e^u\sinh r}.
\end{align*}
Observing again that $2^n\simeq r^{-1}$ and $r^2 \simeq x^2+u^2 \lesssim 1$ on the support of $K_n^*$, we immediately see that
\begin{equation*}
|\vfX_j K_n^\ast (x,u)| \lesssim_\nu r^{-(\nu+2)},\qquad j=0,1,
\end{equation*}
whence
\begin{align*}
\int_{G_\nu} \big|\nabla_\nu K_n^\ast(x,u)\big|\,\dmu(x)\,\du &\lesssim_\nu \int_{2^{-(n+1)}}^{2^{-(n-1)}} r^{-2}\,\dd r \simeq 2^n.
\end{align*}

Thus, Theorem \ref{thm:3} implies the $L^p(G_\nu)$-boundedness of $\tilde\Riesz_0^{(1)}$ on $L^p(G_\nu)$, $p\in(1,2]$.
\end{proof}

\subsection{The Riesz transform \texorpdfstring{$\Riesz_1$}{R1}}

We shall prove the $L^p$-boundedness for $p\in(1,2]$ of the operators $\tilde\Riesz_1^{(j)}$, $j=1,2$, from \eqref{eq:R_1_decomp}. This completes the proof of the $L^p$-bounds for $\Riesz_1$ stated in Theorem \ref{thm:main_riesz}.

Much as before, we first discuss the ``part at infinity'' $\tilde\Riesz_1^{(2)}$ of $\Riesz_1^*$. This part shall again be treated via the operator-valued multiplier theorem, this time applied to the Dunkl transform instead of the Hankel transform.

\begin{prop}
The operator
\begin{equation*}
f\mapsto f\diamond_\nu^{\rD} (\ext_\odd K_1^{(2)})
\end{equation*}
is bounded on $L^p(G_\nu^\rD)$ for all $p\in(1,\infty)$. In particular, the operator
\[
\tilde\Riesz_1^{(2)} : f\mapsto ((\ext_\odd f)\diamond_\nu^{\rD} (\ext_\odd K_1^{(2)}))|_{G_\nu}
\]
is bounded on $L^p(G_\nu)$ for all $p\in(1,\infty)$.
\end{prop}
\begin{proof}
Much as in the proof of Proposition \ref{prop:12}, we shall employ the abstract theory proved in Section \ref{s:opval}, in this case applied to $L=-iD_\nu$.
In other words, here we can take $X=\Omega=X_\nu^\rD$, $\Upsilon=(2^{\nu/2-1} \Gamma(\nu/2))^{-1}\Hank_\nu^\rD$ and $\ell(\omega)=\omega$. Recall from \eqref{eq:Dtr_conv} and \eqref{eq:Dtr_op} that, if $K\in L^2(X_\nu^\rD)$ is bounded, then
\begin{equation*}
f\ast_\nu^\rD K =(\Hank_\nu^\rD K)(-iD_\nu)f.
\end{equation*}
Observe that $-iD_\nu$ has a bounded $\hormander{s}$-calculus for any $s>\nu/2$ (see, e.g., \cite[Section~8]{DzHe}). Thus, again, we can apply Corollary \ref{cor:3} with $N=\lfloor\frac{\nu+5}{2} \rfloor$.

Recall that $G_\nu^\rD$ is the product space $X_\nu^\rD \times \RR$ with measure $\dmub^\rD(x,u) = \dmu^\rD(x) \,\du$. As in Section \ref{ss:dunkl}, we denote
\begin{equation*}
\HaId^\rD \defeq \Hank_\nu^\rD \otimes \id : L^2(G_\nu^\rD)\to L^2(G_\nu^\rD).
\end{equation*}
Clearly, $\HaId^\rD$ is a multiple of an isometry on $L^2(G_\nu^\rD)$. Now let $M : X_\nu^\rD \to \LinBnd(L^2(\RR))$ be measurable and uniformly bounded. Then, from \eqref{eq:55} with $Y=\RR$ we have 
\begin{equation}\label{eq:56D}
\HaId^\rD(M(-iD_\nu) f)(\xi,u) = (M(\xi)\HaId^\rD f(\xi,\cdot))(u).
\end{equation}

On the other hand, observe that, by \eqref{eq:Dtr_diamond},
\begin{equation}\label{eq:57D}
\HaId^\rD (f\diamond_\nu^\rD K)(\xi,u) = \int_\RR \HaId^\rD f(\xi,v) \, \HaId^\rD K(e^{v}\xi,u-v)\,\dd v.
\end{equation}
Let $B_K : X_\nu^\rD \to \LinBnd(L^2(\RR))$ be given by
\begin{equation*}
B_K(\xi) g = \int_\RR B_K^\xi(u,v)\, g(v) \, \dv,\qquad g\in L^2(\RR),
\end{equation*}
where
\begin{equation}\label{eq:kernel_conv_opD}
B_K^\xi(u,v) \defeq \HaId^\rD K(e^v \xi,u-v).
\end{equation}
With this notation, by comparing \eqref{eq:56D} and \eqref{eq:57D}, we obtain that
\begin{equation*}
f\diamond_\nu^\rD K = B_K(-iD_\nu) f.
\end{equation*}

We shall apply this for $K=\ext_\odd K_{1}^{(2)}$. In this case, by \eqref{eq:kernel_conv_opD}, \eqref{eq:62} and \eqref{eq:58},
\[
B_{\ext_\odd K_{1}^{(2)}}^\xi(u,v)= \frac{2i}{\pi} \ind_{\{ u \leq v - 1 \}} \frac{S_\nu(e^v\xi)}{u-v}, \qquad S_\nu(\xi) \defeq  \xi \Hank_{\nu+2}\aux^\nu(\xi).
\]
As a consequence, for any $j \in \NN$, the operator
 $\xi^j\partial_\xi^j B_{\ext_\odd K_{1}^{(2)}}(\xi)$
has integral kernel 
\[
\xi^j\partial_\xi^j B_{\ext_\odd K_{1}^{(2)}}^\xi(u,v) = \frac{2i}{\pi} \ind_{\{ u \leq v - 1 \}} \frac{(\xi^j \partial_\xi^j S_\nu)(e^v\xi)}{u-v}.
\]
The estimates of Lemma \ref{lm:16}, together with \cite[Proposition~5.6(ii)]{Ma23}, then show that each $\xi^j\partial_\xi^j B_{\ext_\odd K_{1}^{(2)}}(\xi)$ is bounded on $L^2(w)$ for any $w\in A_2(\RR)$, with a bound only depending on $j$ and $[w]_{A_2}$. Thus, by Corollary \ref{cor:3}, we deduce the desired $L^p(G_\nu^\rD)$-bound for $f \mapsto f \diamond_\nu^\rD (\ext_\odd K_{1}^{(2)})$.
\end{proof}

We remain with the ``local part'' $\tilde\Riesz_1^{(1)}$ of $\Riesz_1^*$, which again shall be treated with the standard singular integral theory.

\begin{prop}
The operator $\tilde\Riesz_1^{(1)}$ is of weak type $(1,1)$ and bounded on $L^p(G_\nu)$ for all $p\in(1,2]$.
\end{prop}
\begin{proof}
Much as in the proof of Proposition \ref{prop:13}, from the decomposition \eqref{eq:R_1_decomp}, we deduce by difference that $\tilde\Riesz_1^{(1)}$ is bounded on $L^2(G_\nu)$. We shall now employ the Calder\'on--Zygmund theory of Section \ref{s:CZ}, and specifically Corollary \ref{cor:6}, to prove that $\tilde\Riesz_1^{(1)}$ has the claimed boundedness properties.

Namely, by \eqref{eq:62} we can decompose
\begin{equation*}
\ext_\odd K_1^{(1)} = \sum_{n\in\NN} K_n + E,
\end{equation*}
where $E \in L^1(G_\nu^\rD)$, while, for all $n \in \NN$ and $\bdx \in G_\nu^\rD$,
\[
K_n(\bdx) = \frac{\nu}{\pi} \frac{x}{r^{\nu+2}} \chi(2^n r)
\]
and $r=|\bdx|_\dist$, while $\chi \in C^\infty_c(\Rnoz)$ is as in \eqref{eq:dyadicpart}.

Recall that $r^2 \simeq x^2 + u^2 \lesssim 1$ on the support of $K_n$. Thus,
\begin{equation*}
\int_{G_\nu^\rD} |K_{n}(\bdx)| \, (1+2^n |\bdx|_\dist) \, \dmub^\rD(\bdx) \lesssim_\nu \int_{2^{-(n+1)}}^{2^{-(n-1)}} r^{-1}\,\dd r \simeq 1.
\end{equation*}

Moreover, notice that
\begin{equation*}
K^{\bullet}_n(\bdx) = -\frac{\nu}{\pi} \frac{e^{-(\nu +1)u}x}{r^{\nu+2}} \chi(2^n r).
\end{equation*}
Therefore, by \eqref{eq:vfs} and \eqref{eq:44}, if we set $\tilde\eta(t) = t \eta'(t) - (\nu+2)\eta(t)$, then
\begin{align*}
\vfX_0 K^{\bullet}_{n}(\bdx) &= \frac{\nu}{\pi} e^{-(\nu+1) u} x \frac{\nu+1}{r^{\nu+2}}\eta(2^n r) 
- \frac{\nu}{\pi} e^{-(\nu+1) u} x \frac{\tilde\chi(2^n r)}{r^{\nu+3}}\cdot\frac{\sinh u -e^{-u}x^2/2}{\sinh r} ,\\
\vfX_1 K^{\bullet}_{n}(\bdx) &= -\frac{\nu}{\pi} \frac{e^{-\nu u}}{r^{\nu+2}} \chi(2^n r) -\frac{\nu}{\pi} e^{-\nu u} x  \frac{\tilde\chi(2^n r) }{r^{\nu+3}}\cdot\frac{x}{e^u\sinh r}.
\end{align*}
Moreover, as $(K_n^\bullet)_\odd = K_n^\bullet$, from \eqref{eq:tildeX1} we get
\[
\tilde U (K_n^\bullet)_\odd(\bdx) = -\frac{\nu}{\pi} \frac{e^{-\nu u}}{r^{\nu+2}} \chi(2^n r).
\]
All in all, this gives
\[
|\nabla^\nu K^{\bullet}_{n}(\bdx)| + |\tilde U (K_n^\bullet)_\odd(\bdx)| \lesssim_\nu  r^{-(\nu+2)}
\]
and
\[
\int_{G_\nu^\rD} (|\nabla^\nu K^{\bullet}_{n}(\bdx)| + |\tilde U (K_n^\bullet)_\odd(\bdx)|) \,\dmub^\rD(\bdx) \lesssim_\nu \int_{2^{-(n+1)}}^{2^{-(n-1)}} r^{-2}\,\dd r \simeq 2^{n}.
\]
The desired $L^p$-boundedness for $p \in (1,2]$ then follows by Corollary \ref{cor:6}.
\end{proof}

\end{document}